\documentclass[a4paper,3p]{elsarticle}
\usepackage{amsmath}
\usepackage{amsfonts}
\usepackage{amssymb}
\usepackage{amsthm}
\usepackage[numbers]{natbib}
\usepackage{bbm}
\usepackage{mathrsfs} 
\usepackage{xifthen} 
\usepackage{graphicx}
\usepackage{mathrsfs}
\usepackage{dsfont}

\usepackage{mathrsfs}
\usepackage{scalerel}

\DeclareMathAlphabet{\Mymathbb}{U}{bbold}{m}{n}
\DeclareMathAlphabet{\mathpzc}{OT1}{pzc}{m}{it}

\newcommand{\T}{\mbox{\tiny T}}

\newcommand{\Forall}{\forall \,}

\newcommand{\timeSymbol}{t}
\newcommand{\Time}{\timeSymbol} 
\newcommand{\TimeF}{{\timeSymbol}_{\scriptscriptstyle{f}}} 
\newcommand{\tempo}{\Time} 

\newcommand{\FirstSymbol}{\mathcal{G}}
\newcommand{\First}[1][]{
  \ifthenelse{\equal{#1}{}}
  {\FirstSymbol}
  {\FirstSymbol_{#1}}
}
\newcommand{\STFirstSymbol}{\mathbb{G}}
\newcommand{\STFirst}[1][]{
  \ifthenelse{\equal{#1}{}}
  {\STFirstSymbol}
  {\STFirstSymbol_{#1}}
}

\newcommand{\FirstForm}[1][]{
  \ifthenelse{\equal{#1}{}}       
  {\operatorname{I_{\point}}}             
  {\operatorname{I_{#1}}}             
}
\newcommand{\GradSymbol}{\operatorname{\mathbf{\nabla}}}
\newcommand{\Grad}{\GradSymbol}
\newcommand{\Div}{\Grad\cdot}

\newcommand{\GradSurf}{\GradSymbol_{\scriptscriptstyle\First}}
\newcommand{\GradSurfST}{\GradSymbol_{\scriptscriptstyle\STFirst}}

\newcommand{\DivSurf}{\GradSurf\cdot}
\newcommand{\DivSurfST}{\GradSurfST\cdot}

\newcommand{\Der}[1][]
{
  \ifthenelse{\equal{#1}{}}
  {\partial}
  {\partial_{\scriptscriptstyle{#1}}}
}

\newcommand{\DerPar}[2]{\frac{\Der #1}{\Der #2}}

\newcommand{\DerTot}[2][t]
{
  \ifthenelse{\equal{#2}{}}
  {\frac{d #2}{dt}}
  {\frac{d #2}{d #1}}
}
\newcommand{\Diffsymbol}{\operatorname{d}\!}
\newcommand{\Diff}[2][]{
  \ifthenelse{\equal{#1}{}}
  {\Diffsymbol{#2}}
  {\Diffsymbol{#2}_{\scriptscriptstyle{#1}}}
}
\newcommand{\DirDerSymb}{D}
\newcommand{\DirDer}[2][]
{
  \ifthenelse{\equal{#1}{}}
  {\DirDerSymb^{#2}}
  {\DirDerSymb^{#2}_{#1}}
}

\newcommand{\BigO}[1]{\ensuremath{\operatorname{\mathcal{O}}\!\left(#1\right)}}

\newcommand{\REALsymbol}{\mathbb{R}}
\newcommand{\REAL}[1][]{
  \ifthenelse{\equal{#1}{}}
  {\REALsymbol}
  {{\REALsymbol}^{#1}}
}
\newcommand{\EUCLsymbol}{\mathbb E}
\newcommand{\EUCL}[1][]{
  \ifthenelse{\equal{#1}{}}
  {\EUCLsymbol}
  {{\EUCLsymbol}^{#1}}
}
\newcommand{\NATURALsymbol}{\mathbb N}
\newcommand{\NATURAL}[1][]{
  \ifthenelse{\equal{#1}{}}
  {\NATURALsymbol}
  {{\NATURALsymbol}^{#1}}
}

\newcommand{\ABS}[2][]
{
  \ifthenelse{\equal{#1}{}}
  {\left| #2 \right|}
  {\left| #2 \right|_{#1}}
}

\newcommand{\NORM}[2][]
{
  \ifthenelse{\equal{#1}{}}
  {\left\| #2 \right\|}
  {\left\| #2 \right\|_{#1}}
}

\newcommand{\SCAL}[3][]
{
  \ifthenelse{\equal{#1}{}}
  {\left\langle{#2},{#3}\right\rangle}
  {\left\langle{#2},{#3}\right\rangle_{#1}}
}
\newcommand{\SCALF}[3][]
{
  \ifthenelse{\equal{#1}{}}
  {\left({#2},{#3}\right)}
  {\left({#2},{#3}\right)_{#1}}
}

\newcommand{\scalprodSurf}[3][]
{
  \ifthenelse{\equal{#1}{}}
  {\left\langle {#2},{#3} \right\rangle_{\scriptscriptstyle\First}}
  {\left\langle {#2},{#3} \right\rangle_{{#1}}}
}
\newcommand{\scalprodSurfST}[3][]
{
  \ifthenelse{\equal{#1}{}}
  {\left\langle {#2},{#3} \right\rangle_{\scriptscriptstyle\STFirst}}
  {\left\langle {#2},{#3} \right\rangle_{{#1}}}
}

\newcommand{\DET}[1]{\ensuremath{\operatorname{det}(#1)}}

\newcommand{\Point}{\mathbf{P}}
\newcommand{\point}[1][]
{
  \ifthenelse{\equal{#1}{}}
  {\mathbf{p}}
  {\mathbf{p}_{#1}}
}

\newcommand{\midPoint}{\mathbf{m}}

\newcommand{\Dt}{\Delta\tempo}

\newcommand{\From}{:}
\newcommand{\To}{\rightarrow}

\newcommand{\Closure}[1]{\overline{#1}}

\newcommand{\RegionSymb}{R}
\newcommand{\Region}[1][]
{
  \ifthenelse{\equal{#1}{}}
  {\RegionSymb}
  {\RegionSymb^{#1}}
}
\newcommand{\bRegion}[1][]
{
  \ifthenelse{\equal{#1}{}}
  {\overline{\RegionSymb}}
  {\overline{\RegionSymb}^{#1}}
}

\newcommand{\diam}{\operatorname{diam}}

\newcommand{\IT}[1][] 
{
  \ifthenelse{\equal{#1}{}}
  {[0,\TimeF]}
  {[\Time^{#1},\Time^{#1+1}]}
}
\newcommand{\Interval}[1][]
{
  \ifthenelse{\equal{#1}{}}
  {I}
  {I_{#1}}
}

\newcommand{\SurfDomainsymb}{\Gamma}
\newcommand{\SurfDomain}[1][]{
  \ifthenelse{\equal{#1}{}}
  {\SurfDomainsymb}
  {\SurfDomainsymb_{#1}}
}
\newcommand{\SurfDomainBndsymb}{\partial\Gamma}
\newcommand{\SurfDomainBnd}[1][]{
  \ifthenelse{\equal{#1}{}}
  {\SurfDomainBndsymb}
  {\SurfDomainBndsymb_{#1}}
}
\newcommand{\ClosedSurfDomain}[1][]{
  \ifthenelse{\equal{#1}{}}
  {\Closedsymb{\SurfDomain}}
  {\Closedsymb{\SurfDomain[#1]}}
}

\newcommand{\heightsymb}{\mathcal{H}}
\newcommand{\height}[1][]{
  \ifthenelse{\equal{#1}{}}
  {\heightsymb}
  {\heightsymb_{#1}}
}

\newcommand{\BSMsymbol}{H}

\newcommand{\BSM}[1][]
{
  \ifthenelse{\equal{#1}{}}
  {\BSMsymbol}
  {\BSMsymbol_{#1}}
}

\newcommand{\Surf}{\mathcal{S}}

\newcommand{\SurfBndsymb}{\partial\Surf}
\newcommand{\SurfBnd}[1][]{
  \ifthenelse{\equal{#1}{}}
  {\SurfBndsymb}
  {\SurfBndsymb_{#1}}
}
\newcommand{\Closedsymb}[1]{\bar{#1}}
\newcommand{\ClosedSurf}[1][]{
  \ifthenelse{\equal{#1}{}}
  {\Closedsymb{\Surf}}
  {\Closedsymb{\Surf[#1]}}
}

\newcommand{\Vector}[1]{\mathbf{#1}}

\newcommand{\grav}{g}

\newcommand{\Depth}{\eta}

\newcommand{\MapUsymb}{\phi}
\newcommand{\MapU}[1][]
{
  \ifthenelse{\equal{#1}{}}
    {\MapUsymb}
    {\MapUsymb_{#1}}
}
\newcommand{\MapVsymb}{\psi}
\newcommand{\MapV}[1][]
{
  \ifthenelse{\equal{#1}{}}
    {\MapVsymb}
    {\MapVsymb_{#1}}
}
\newcommand{\Transsymb}{\Phi}
\newcommand{\Trans}[1][]
{
  \ifthenelse{\equal{#1}{}}
    {\Transsymb} 
    {\Transsymb_{\scriptscriptstyle{#1}}}
}
\newcommand{\InvMapsymb}{\Psi}
\newcommand{\InvMap}[1][]
{
  \ifthenelse{\equal{#1}{}}
    {\InvMapsymb}
    {\InvMapsymb_{\scriptscriptstyle{#1}}}
}

\newcommand{\fsymb}{f}
\newcommand{\scalFun}[1][]
{
  \ifthenelse{\equal{#1}{}}
  {\fsymb}
  {\fsymb_{#1}}
}
\newcommand{\tscalFun}[1][]
{
  \ifthenelse{\equal{#1}{}}
  {\tilde{\fsymb}}
  {\tilde{\fsymb}_{#1}}
}
\newcommand{\gsymb}{g}
\newcommand{\scalFung}[1][]
{
  \ifthenelse{\equal{#1}{}}
  {\gsymb}
  {\gsymb_{#1}}
}
\newcommand{\Fsymb}{F}
\newcommand{\FvecFun}[1][]
{
  \ifthenelse{\equal{#1}{}}
  {\Fsymb}
  {\Fsymb_{#1}}
}
\newcommand{\tFvecFun}[1][]
{
  \ifthenelse{\equal{#1}{}}
  {\tilde{\Fsymb}}
  {\tilde{\Fsymb}_{#1}}
}
\newcommand{\Ffunc}[2][]
{
  \ifthenelse{\equal{#1}{}}
  {\ensuremath{\Fsymb_{#2}}}
  {\ensuremath{\Fsymb^{#1}_{#2}}}
}

\newcommand{\PrincipalK}[1][]
{
  \ifthenelse{\equal{#1}{}}
  {k}
  {k_{#1}}
}

\newcommand{\vecsymb}{u}
\newcommand{\vecFun}[1][]{
  \ifthenelse{\equal{#1}{}}
  {\Vector{\vecsymb}}
  {\vecsymb^{#1}}
}
\newcommand{\tvecFun}[1][]{
  \ifthenelse{\equal{#1}{}}
  {\tilde{\vecsymb}}
  {\tilde{\vecsymb}_{#1}}
}

\newcommand{\vvsymb}{v}
\newcommand{\vv}[1][]{
  \ifthenelse{\equal{#1}{}}
  {\mathbf{\vvsymb}}
  {\vvsymb^{#1}}
}
\newcommand{\wwsymb}{w}
\newcommand{\ww}[1][]
{
  \ifthenelse{\equal{#1}{}}
  {\mathbf{\wwsymb}}
  {\wwsymb^{#1}}
}
\newcommand{\uusymb}{u}
\newcommand{\uu}[1][]
{
  \ifthenelse{\equal{#1}{}}
  {\mathbf{\uusymb}}
  {\uusymb^{#1}}
}

\newcommand{\VecFieldSymbol}{X}
\newcommand{\VecField}[1][]
{
  \ifthenelse{\equal{#1}{}}
  {\VecFieldSymbol}
  {\VecFieldSymbol^{#1}}
}
\newcommand{\VecFieldYSymbol}{Y}
\newcommand{\VecFieldY}[1][]
{
  \ifthenelse{\equal{#1}{}}
  {\VecFieldYSymbol}
  {\VecFieldYSymbol^{#1}}
}

\newcommand{\xvsymb}{x}
\newcommand{\xv}[1][]
{
  \ifthenelse{\equal{#1}{}}
  {\mathbf{\xvsymb}}
  {\mathbf{\xvsymb}_{\scriptscriptstyle{#1}}}
}
\newcommand{\xvcomp}[1][]{
  \ifthenelse{\equal{#1}{}}
  {\xvsymb}
  {\xvsymb^{\scriptscriptstyle{#1}}}
}
\newcommand{\xcg}[1][]{
  \ifthenelse{\equal{#1}{}}
  {\xvcomp[1]}
  {\xvcomp[1]_{\scriptscriptstyle{#1}}}
}
\newcommand{\ycg}[1][]{
  \ifthenelse{\equal{#1}{}}
  {\xvcomp[2]}
  {\xvcomp[2]_{\scriptscriptstyle{#1}}}
}
\newcommand{\zcg}[1][]{
  \ifthenelse{\equal{#1}{}}
  {\xvcomp[3]}
  {\xvcomp[3]_{\scriptscriptstyle{#1}}}
}

\newcommand{\svsymb}{s}
\newcommand{\sv}[1][]{
  \ifthenelse{\equal{#1}{}}
  {\mathbf{\svsymb}}
  {\mathbf{\svsymb}_{\scriptscriptstyle{#1}}}
}
\newcommand{\svcomp}[1][]
{
  \ifthenelse{\equal{#1}{}}
  {\svsymb}
  {\svsymb^{\scriptscriptstyle{#1}}}
}
\newcommand{\xcl}[1][]{
  \ifthenelse{\equal{#1}{}}
   {\svcomp[1]}
   {\svcomp[1]_{\scriptscriptstyle{#1}}}
}
\newcommand{\ycl}[1][]{
  \ifthenelse{\equal{#1}{}}
   {\svcomp[2]}
   {\svcomp[2]_{\scriptscriptstyle{#1}}}
}
\newcommand{\zcl}[1][]{
  \ifthenelse{\equal{#1}{}}
   {\svcomp[3]}
   {\svcomp[3]_{\scriptscriptstyle{#1}}}
}

\newcommand{\ProjSymb}{\operatorname{\pi}}
\newcommand{\ProjFun}[2][]{
  \ifthenelse{\equal{#1}{}}
  {\ProjSymb\left(#2\right)}
  {\ProjSymb_{\scriptscriptstyle{#1}}\left(#2\right)}
}
\newcommand{\Prm}[1][]
{
  \ifthenelse{\equal{#1}{}}
  {\operatorname{pr}}
  {\operatorname{pr}_{\scriptscriptstyle{#1}}}
}
\newcommand{\TanPlane}[2][]
{
  \ifthenelse{\equal{#1}{}}
  {T_{\scriptscriptstyle{\point}}#2}
  {T_{\scriptscriptstyle{#1}}#2}
}

\newcommand{\SubsetSymbol}{\mathcal{U}}
\newcommand{\SubsetU}[1][]
{
  \ifthenelse{\equal{#1}{}}
  {{U}}
  {{U}_{#1}}
}
\newcommand{\SubsetV}[1][]
{
  \ifthenelse{\equal{#1}{}}
  {{V}}
  {{V}_{#1}}
}
\newcommand{\SubsetW}[1][]
{
  \ifthenelse{\equal{#1}{}}
  {{W}}
  {{W}_{#1}}
}
\newcommand{\NeighSymbol}{\mathcal{N}}
\newcommand{\Neigh}[1][]
{
  \ifthenelse{\equal{#1}{}}
  {\NeighSymbol_{\point}}
  {\NeighSymbol_{#1}}
}
\newcommand{\NeighSurf}[1][]
{
  \ifthenelse{\equal{#1}{}}
  {\SubsetSymbol_{\point}}
  {\SubsetSymbol_{#1}}
}
\newcommand{\NormSymb}{\mathbf{N}}
\newcommand{\normalvec}{\NormSymb}
\newcommand{\normalSurf}[1][]
{
  \ifthenelse{\equal{#1}{}}
  {\NormSymb}
  {\NormSymb(#1)}
}
\newcommand{\normalInterp}[1][]
{
  \ifthenelse{\equal{#1}{}}
   {\tilde{\NormSymb}}
   {\tilde{\NormSymb}_{\scriptscriptstyle{#1}}}
}

\newcommand{\normalEdge}{\mathbf{\nu}} 
\newcommand{\basisCC}{t}
\newcommand{\basisGC}{e}
\newcommand{\vecBaseGC}[1][]
{
  \ifthenelse{\equal{#1}{}}
  {\mathbf{\basisGC}}
  {\mathbf{\basisGC}_{#1}}
}
\newcommand{\vecBasePhys}[1][]
{
  \ifthenelse{\equal{#1}{}}
  {\mathbf{\basisGC}}
  {\mathbf{\basisGC}_{#1}}
}

\newcommand{\vecBaseCCcv}[1][]
{
  \ifthenelse{\equal{#1}{}}
  {\mathbf{\basisCC}}
  {\mathbf{\basisCC}_{#1}}
}

\newcommand{\tvecBaseCCcv}[1][]
{
  \ifthenelse{\equal{#1}{}}
  {\tilde{\mathbf{\basisCC}}}
  {\tilde{\mathbf{\basisCC}}_{#1}}
}
\newcommand{\hvecBaseCCcv}[1][]
{
  \ifthenelse{\equal{#1}{}}
  {\hat{\mathbf{\basisCC}}}
  {\hat{\mathbf{\basisCC}}_{#1}}
}
\newcommand{\vecBaseCCctrv}[1][]
{
  \ifthenelse{\equal{#1}{}}
  {\mathbf{\basisCC}}
  {\mathbf{\basisCC}^{#1}}
}

\newcommand{\metricsymbol}{g}

\newcommand{\metrcoef}[1]{\metricsymbol_{\mbox{\tiny{#1#1}}}}
\newcommand{\metrcoefsymbol}{h}
\newcommand{\metrcoefH}[1]{\metrcoefsymbol_{\mbox{\tiny{(#1)}}}}

\newcommand{\first}[1]{
  \IfEqCase{#1}{
    {1}{\operatorname{E}}
    {2}{\operatorname{F}}
    {3}{\operatorname{G}}
  }
  [\PackageError{first}{Undefined option to first: #1}{}]%
}
\newcommand{\SecondFormSymbol}{\ensuremath{\operatorname{II}}}
\newcommand{\SecondForm}[1][]
{
  \ifthenelse{\equal{#1}{}}
  {\SecondFormSymbol_{\point}}
  {\SecondFormSymbol_{#1}}
}
\newcommand{\second}[1]{
  \IfEqCase{#1}{
    {1}{\operatorname{e}}
    {2}{\operatorname{f}}
    {3}{\operatorname{g}}
  }
  [\PackageError{first}{Undefined option to first: #1}{}]
}
\newcommand{\WeigSymbol}{\mathcal{W}}
\newcommand{\Weig}[1][]
{
  \ifthenelse{\equal{#1}{}}
  {\WeigSymbol}
  {\WeigSymbol_{#1}}
}
\newcommand{\ChristSymb}[2]{\Gamma_{#1}^{#2}}

\newcommand{\velSymbol}{u}
\newcommand{\vectvel}[1][]
{
   \ifthenelse{\equal{#1}{}}
   {\vec{\velSymbol}}
   {\vec{\velSymbol}(#1)}
}

\newcommand{\velcompContr}[2][i]
{
   \ifthenelse{\equal{#2}{}}
   {\velSymbol^{#1}}
   {\velSymbol^{#1}(#2)}}
\newcommand{\velcompPhys}[2][i]
{
   \ifthenelse{\equal{#2}{}}
   {\velSymbol_{(#1)}}
   {\velSymbol_{(#1)}(#2)}}
\newcommand{\velcomp}{\velcompContr}

\newcommand{\velSymbolRP}{v}
\newcommand{\velcompContrRP}[2][i]
{
   \ifthenelse{\equal{#2}{}}
   {\velSymbolRP^{#1}}
   {\velSymbolRP^{#1}(#2)}}
\newcommand{\velRP}[1][]
{
  \ifthenelse{\equal{#1}{}}
  {\velSymbolRP}
  {\velSymbolRP_{#1}}
}

\newcommand{\velcompApprox}[2][i]
{
   \ifthenelse{\equal{#2}{}}
   {\velSymbol^{#1)}}
   {\velSymbol^{#1}_{(#2)}}
}

\newcommand{\VelSymbol}{U}
\newcommand{\vectVel}[1][]
{
   \ifthenelse{\equal{#1}{}}
   {\vec{\VelSymbol}}
   {\vec{\VelSymbol}(#1)}
}
\newcommand{\Velcomp}[2][i]
{
   \ifthenelse{\equal{#2}{}}
   {\VelSymbol^{#1}}
   {\VelSymbol^{#1}(#2)}
}
\newcommand{\VprimoSymbol}{\tilde{u}}
\newcommand{\Vprimo}[1][]
{
   \ifthenelse{\equal{#1}{}}
   {\VprimoSymbol}
   {\VprimoSymbol(#1)}
}
\newcommand{\VprimoComp}[2][i]
{
   \ifthenelse{\equal{#2}{}}
   {\VprimoSymbol^{#1}}
   {\VprimoSymbol^{#1}(#2)}
}
\newcommand{\ttvelSymbol}{\tilde{\Mymathbb{u}}}
\newcommand{\ttvel}[1][]
{
   \ifthenelse{\equal{#1}{}}
   {\mathbf{\ttvelSymbol}}
   {\mathbf{\ttvelSymbol}(#1)}
}
\newcommand{\ttvelComp}[2][i]
{
   \ifthenelse{\equal{#2}{}}
   {\ttvelSymbol^{#1}}
   {\ttvelSymbol^{#1}(#2)}
}

\newcommand{\MatAlphaSymbol}{\mathbb{A}}
\newcommand{\MatAlpha}[1][]{%
  \ifthenelse{\equal{#1}{}}
  {\MatAlphaSymbol}
  {\MatAlphaSymbol_{#1}}
}

\newcommand{\QSymbol}{q}
\newcommand{\Qdisch}[1][]
{
   \ifthenelse{\equal{#1}{}}
   {\mathbf{\QSymbol}}
   {\mathbf{\QSymbol}(#1)}
}
\newcommand{\Qcomp}[2][i]
{
   \ifthenelse{\equal{#2}{}}
   {\QSymbol^{#1}}
   {\QSymbol^{#1}(#2)}
}
\newcommand{\Qvect}[1][]
{
   \ifthenelse{\equal{#1}{}}
   {\vec{\QSymbol}}
   {\vec{\QSymbol}(#1)}
}
\newcommand{\FricSymbol}{f}
\newcommand{\vectFric}[1][]
{
   \ifthenelse{\equal{#1}{}}
   {\mathbf{\FricSymbol}}
   {\mathbf{\FricSymbol}_{\scriptscriptstyle{#1}}}
}
\newcommand{\Friccomp}[2][i]
{
   \ifthenelse{\equal{#2}{}}
   {\FricSymbol_{#1}}
   {\FricSymbol_{#1}(#2)}}
\newcommand{\BFsymbol}{\tau}
\newcommand{\BottomFriction}[1][]{
  \ifthenelse{\equal{#1}{}}
  {\BFsymbol_{b}}
  {\BFsymbol_{b}^{#1}}
}

\newcommand{\ProjMatSymb}{\mathbb{P}}
\newcommand{\ProjMat}[1][]{
  \ifthenelse{\equal{#1}{}}
  {\ProjMatSymb}
  {\ProjMatSymb_{#1}}
}
\newcommand{\IDSymbol}{\mathbb{I}}
\newcommand{\IDtens}[1][]{
  \ifthenelse{\equal{#1}{}}
  {\IDSymbol}
  {\IDSymbol(#1)}
}

\newcommand{\tensSymbol}{\mathbb{T}}
\newcommand{\tenscompSymbol}{\tau}
\newcommand{\tens}[1][]{
  \ifthenelse{\equal{#1}{}}
  {\tensSymbol}
  {\tensSymbol(#1)}
}
\newcommand{\tenscomp}[2][ij]
{
  \ifthenelse{\equal{#2}{}}
  {\tenscompSymbol^{#1}}
  {\tenscompSymbol^{#1}(#2)}
}
\newcommand{\tensrow}[2][i]
{
  \ifthenelse{\equal{#2}{}}
  {\tensSymbol^{(#1)}}
  {\tensSymbol^{(#1)}(#2)}
}
\newcommand{\TensSymbol}{\mathbf{T}}
\newcommand{\Tens}[1][]{
  \ifthenelse{\equal{#1}{}}
  {\TensSymbol}
  {\TensSymbol_{#1}}
}

\newcommand{\TensCompSymbol}{\TensSymbol}
\newcommand{\TensComp}[2][ij]
{
  \ifthenelse{\equal{#2}{}}
  {\TensCompSymbol^{#1}}
  {\TensCompSymbol^{#1}(#2)}
}

\newcommand{\tensPrimoSymbol}{\tilde{\mathbf{\tau}}}
\newcommand{\tensPrimo}[1][]{
  \ifthenelse{\equal{#1}{}}
  {\tensPrimoSymbol}
  {\tensPrimoSymbol(#1)}
}
\newcommand{\tensPrimoCompSymbol}{\tensPrimoSymbol}
\newcommand{\tensPrimoComp}[2][ij]
{
  \ifthenelse{\equal{#2}{}}
  {\tensPrimoCompSymbol^{#1}}
  {\tensPrimoCompSymbol^{#1}(#2)}
}
\newcommand{\MCxl}[1][]{\ifthenelse{\equal{#1}{}}{h_{(1)}}{h_{(1),#1}}} 
\newcommand{\MCyl}[1][]{\ifthenelse{\equal{#1}{}}{h_{(2)}}{h_{(2),#1}}} 
\newcommand{\MCzl}[1][]{\ifthenelse{\equal{#1}{}}{h_{(3)}}{h_{(3),#1}}}

\newcommand{\MPsymb}{h}

\newcommand{\smallestH}[1][]
{
  \ifthenelse{\equal{#1}{}}
    {{l}}
    {{l}_{#1}}
}
\newcommand{\meshparam}[1][]
{
  \ifthenelse{\equal{#1}{}}
    {\MPsymb}
    {\MPsymb_{#1}}
}
\newcommand{\InradiusSymbol}{r}
\newcommand{\Inradius}[1][]
{
  \ifthenelse{\equal{#1}{}}
  {\InradiusSymbol}
  {\InradiusSymbol_{\scriptscriptstyle{#1}}}
}
\newcommand{\Tsymb}{\mathcal{T}}
\newcommand{\Triang}[1][]
{
  \ifthenelse{\equal{#1}{}}
    {\Tsymb}
    {\Tsymb_{#1}}
}
\newcommand{\TriangH}[1][]
{
  \ifthenelse{\equal{#1}{}}
    {\Tsymb_{\meshparam}}
    {\Tsymb_{#1}}
}
\newcommand{\Edgesymb}{e}
\newcommand{\Edge}[1][]{
  \ifthenelse{\equal{#1}{}}
    {\Edgesymb}
    {\Edgesymb_{#1}}
}
\newcommand{\EdgeH}[1][]{
  \ifthenelse{\equal{#1}{}}
    {\Edgesymb_{\meshparam}}
    {\Edgesymb_{\meshparam,#1}}
}
\newcommand{\NEdge}[1][]{
  \ifthenelse{\equal{#1}{}}
    {N_{\Edgesymb}}
    {N_{\Edgesymb({#1})}}
}
\newcommand{\Cellsymb}{K}
\newcommand{\Cell}[1][]{
  \ifthenelse{\equal{#1}{}}
    {\Cellsymb}
    {\Cellsymb_{#1}}
}
\newcommand{\CellH}[1][]{
  \ifthenelse{\equal{#1}{}}
    {\Cellsymb_{\meshparam}}
    {\Cellsymb_{\meshparam,#1}}
}
\newcommand{\areaSymb}{\mathcal{A}}
\newcommand{\CellArea}[1][]
{
  \ifthenelse{\equal{#1}{}}
    {\areaSymb_{\Cell}}
    {{A({#1})}}
}
\newcommand{\CellHArea}[1][]
{
  \ifthenelse{\equal{#1}{}}
    {\areaSymb_{\CellH}}
    {\areaSymb_{\meshparam,#1}}
}

\newcommand{\NCell}[1][]{
  \ifthenelse{\equal{#1}{}}
    {N_{\Cellsymb}}
    {N_{\Cellsymb({#1})}}
}

\newcommand{\edgeLength}[1][]
{
  \ifthenelse{\equal{#1}{}}
  {\ABS{\Edge}}
  {\ABS{\Edge_{#1}}}
}
 
\newcommand{\edgeHLength}[1][]
{
  \ifthenelse{\equal{#1}{}}
  {\lengthSymb_{\EdgeH}}
  {\lengthSymb_{\meshparam,#1}}
}
\newcommand{\subVolume}[2]{D_{#1}^{#2}}

\newcommand{\Sourcesymb}{\mathbf{S}}
\newcommand{\Source}[1][]
{
  \ifthenelse{\equal{#1}{}}
    {\Sourcesymb}
    {\Sourcesymb_{#1}}
}

\newcommand{\SourceEdge}[1][]{
  \ifthenelse{\equal{#1}{}}
    {\Sourcesymb_{ij}}
    {\Sourcesymb_{#1}}
}

\newcommand{\ConservVar}[2]{\mathbf{U}_{#1}^{{#2}}}

\newcommand{\bConservVar}[2]{\overline{\mathbf{U}}_{#1}^{{#2}}}
\newcommand{\FluxSWE}{\mathbf{F}}
\newcommand{\Flux}{\underline{\underline{F}}}

\newcommand{\FluxEdgesymb}{\mathbf{F}}
\newcommand{\FluxEdge}[1][]{
  \ifthenelse{\equal{#1}{}}
    {\FluxEdgesymb_{ij}}
    {\FluxEdgesymb_{#1}}
}

\newcommand{\numFlux}[1][]
{
  \ifthenelse{\equal{#1}{}}
  {\tilde{\fsymb}}
  {\tilde{\fsymb}_{#1}}
}

\newcommand{\FluxFuncNormSymbol}{\mathbf{F}}
\newcommand{\FluxFuncNorm}[1][]{
  \ifthenelse{\equal{#1}{}}
    {\FluxFuncNormSymbol^{\normalEdge}}
    {\FluxFuncNormSymbol^{\normalEdge}_{#1}}    
}

\newcommand{\JacobianSymbol}{\mathbf{A}}
\newcommand{\Jacobian}[1][]{
  \ifthenelse{\equal{#1}{}}
    {\JacobianSymbol}
    {\JacobianSymbol_{#1}}    
}
\newcommand{\EValSymbol}{\lambda}
\newcommand{\EVal}[1][]{
  \ifthenelse{\equal{#1}{}}
  {\EValSymbol}
  {\EValSymbol_{#1}}
}
\newcommand{\EVecSymbol}{\mathbf{r}}
\newcommand{\EVec}[2][]{
  \ifthenelse{\equal{#1}{}}
  {\EVecSymbol^{(#2)}}
  {\EVecSymbol^{(#2)}_{#1}}
}

\newcommand{\midPointEdge}[1][]
{
  \ifthenelse{\equal{#1}{}}
  {\midPoint_{\scriptscriptstyle\Edge}}
  {\midPoint_{\scriptscriptstyle\Edge[#1]}}
}
\newcommand{\gpPointEdgeDG}[1][]
{
  \ifthenelse{\equal{#1}{}}
  {\point_{\scriptscriptstyle\Edge}}
  {\point_{\scriptscriptstyle\Edge,#1}}
}
\newcommand{\midPointCell}[1][]
{
  \ifthenelse{\equal{#1}{}}
  {\midPoint_{\scriptscriptstyle\Cell}}
  {\midPoint_{\scriptscriptstyle\Cell[#1]}}
}

\newcommand{\speedRS}[2][]
{
  \ifthenelse{\equal{#2}{}}
  {S_{#2}}
  {S_{#2}^{#1}}
}

\newcommand{\ContSymbol}{C}
\newcommand{\Cont}[1][]{
  \ifthenelse{\equal{#1}{}}
  {\ContSymbol^{0}}
  {\ContSymbol^{#1}}
}
\newcommand{\Cinf}[1][]{
  \ifthenelse{\equal{#1}{}}
  {\ContSymbol^{\infty}}
  {\ContSymbol^{\infty}(#1)}
}
\newcommand{\SobSymbol}{W}
\newcommand{\HilbSymbol}{H}
\newcommand{\Sob}[2][]{
  \ifthenelse{\equal{#1}{}}
  {\HilbSymbol^{#2}}
  {\SobSymbol^{#2,#1}}  
}

\newcommand{\LspaceSymb}{L}
\newcommand{\Lspace}[1][]{
  \ifthenelse{\equal{#1}{}}
  {\LspaceSymb^{2}}
  {\LspaceSymb^{#1}}  
}

\newcommand{\TestSpSymbol}{\mathcal{V}}
\newcommand{\TestSpace}[1][]{
  \ifthenelse{\equal{#1}{}}
  {\TestSpSymbol({\SurfDomain})}
  {\TestSpSymbol_{#1}(\SurfDomain)}
}
\newcommand{\TestSpaceTime}[1][]{
  \ifthenelse{\equal{#1}{}}
  {\TestSpSymbol(\SurfDomain{}(\tempo))}
  {\TestSpSymbol_{#1}(\SurfDomain{}(\tempo))}
}

\newcommand{\AngleSymbol}{\theta}
\newcommand{\DevAngle}[1][]
{
  \ifthenelse{\equal{#1}{}}
  {\AngleSymbol}
  {\AngleSymbol_{\scriptscriptstyle{#1}}}
}
\newcommand{\relheightSymb}{\pi}
\newcommand{\relheight}[1][]
{
  \ifthenelse{\equal{#1}{}}
  {\relheightSymb_{\scriptscriptstyle{\SurfDomain}}}
  {\relheightSymb_{\scriptscriptstyle{#1}}}
}

\newcommand{\SolSymbol}{u}
\newcommand{\Sol}{\SolSymbol}

\newcommand{\BilinearStiffSymbol}{a}
\newcommand{\BilinearStiff}[3][]
{
  \ifthenelse{\equal{#1}{}}
  {\BilinearStiffSymbol(#2,#3)}    
  {\BilinearStiffSymbol_{#1}(#2,#3)}    
}
\newcommand{\BilinearAdvSymbol}{b}
\newcommand{\BilinearAdv}[3][]
{
  \ifthenelse{\equal{#1}{}}
  {\BilinearAdvSymbol(#2,#3)}    
  {\BilinearAdvSymbol_{#1}(#2,#3)}    
}
\newcommand{\BilinearMassSymbol}{m}
\newcommand{\BilinearMass}[3][]
{
  \ifthenelse{\equal{#1}{}}
  {\BilinearMassSymbol(#2,#3)}    
  {\BilinearMassSymbol_{#1}(#2,#3)}    
}
\newcommand{\BilinearReactSymbol}{c}
\newcommand{\BilinearReact}[3][]
{
  \ifthenelse{\equal{#1}{}}
  {\BilinearReactSymbol(#2,#3)}    
  {\BilinearReactSymbol_{#1}(#2,#3)}    
}

\newcommand{\TestSymbol}{v}
\newcommand{\Test}[1][]
{
  \ifthenelse{\equal{#1}{}}
  {\TestSymbol}  
  {\TestSymbol_{\scriptscriptstyle{#1}}}
}

\newcommand{\coeffSol}[2][]
{
  \ifthenelse{\equal{#1}{}}
  {\Sol_{#2}}  
  {\Sol_{#2}^{\scriptscriptstyle{#1}}}
}

\newcommand{\nodalBasis}[2][]
{
  \ifthenelse{\equal{#1}{}}
  {\varphi_{\scriptscriptstyle{#2}}}  
  {\varphi_{\scriptscriptstyle{#2}}^{\scriptscriptstyle{#1}}}
}

\newcommand{\nNodes}[1][]
{
  \ifthenelse{\equal{#1}{}}
  {N}
  {N_{\scriptscriptstyle{#1}}}
}

\newcommand{\TestApprox}[1][]
{
  \ifthenelse{\equal{#1}{}}
  {\TestSymbol_{\scriptscriptstyle{\meshparam}}}  
  {\TestSymbol_{\scriptscriptstyle{\meshparam,#1}}}
}

\newcommand{\viscositySymb}{\nu}
\newcommand{\viscosity}[1][]
{
  \ifthenelse{\equal{#1}{}}
  {\viscositySymb}
  {\viscositySymb_{\scriptscriptstyle{#1}}}
}

\newcommand{\ResidualSymbol}{R}
\newcommand{\Residual}[1][]
{
  \ifthenelse{\equal{#1}{}}
  {\ResidualSymbol}
  {\ResidualSymbol_{#1}}
}

\makeatletter
\def\ps@pprintTitle{%
  \let\@oddhead\@empty
  \let\@evenhead\@empty
  \def\@oddfoot{\reset@font\hfil\thepage\hfil}
  \let\@evenfoot\@oddfoot
}
\makeatother

\title{A geometrically intrinsic Lagrangian-Eulerian scheme for 2D
    Shallow Water Equations with variable topography and discontinuous
    data}

  \author[1]{Eduardo Abreu}
  \address[1]{Department of Applied Mathematics,
      University of Campinas, Brazil}
  \author[2]{Elena Bachini}
  \address[2]{Institute of Scientific Computing, TU
      Dresden, Germany}
  \author[3]{John P\'erez}
  \address[3]{ITM-University Institution, Medellin, Colombia}
  \author[4]{Mario Putti}
  \address[4]{Department of Mathematics ``Tullio
      Levi-Civita'', University of Padua, Italy}
  
\begin{document}

\begin{frontmatter}
  
  \begin{abstract}
    We present a Lagrangian-Eulerian scheme to solve the shallow water
    equations in the case of spatially variable bottom geometry.
    Using a local curvilinear reference system anchored on the bottom
    surface, we develop an effective first-order and high-resolution 
    space-time discretization of the no-flow
    surfaces and solve a Lagrangian initial value problem that
    describes the evolution of the balance laws governing the
    geometrically intrinsic shallow water equations. The evolved
    solution set is then projected back to the original surface grid
    to complete the proposed Lagrangian-Eulerian formulation.
    The resulting scheme maintains monotonicity and captures shocks without providing excessive numerical dissipation also in the presence of non-autonomous fluxes such as
    those arising from the geometrically intrinsic shallow water
    equation on variable topographies. We provide a representative set of numerical examples to illustrate the accuracy and robustness of 
    the proposed Lagrangian-Eulerian formulation for two-dimensional
    surfaces with general curvatures and discontinuous initial
    conditions.
  \end{abstract}

  \begin{keyword}
    Balance laws on surface\sep
    shallow water equations \sep
    spatially variable topography \sep
    intrinsic Lagrangian-Eulerian scheme \sep
    no-flow surfaces
  \end{keyword}

\end{frontmatter}

\section{Introduction}
\label{sec:intro}

Partial Differential Equations (PDEs) modeling processes occurring on
surfaces have been the subject of several studies in recent
years. Typically, the equations governing these processes are
developed from their three-dimensional counterparts and letting the
dimension normal to the surface go to zero or, alternatively, by
averaging the three-dimensional equations along the local normal
direction and employing parameter perturbation analysis to impose the
condition that the process occurs prevalently along a direction
parallel to the surface. Recent examples of these two approaches can
be found in~\cite{art:Nestler2019}, where a model for the elastic
equilibrium of nematic liquid crystals was derived using the first
method, and in~\cite{art:BP20} where averaging and perturbation
analysis are combined to derive the shallow water equations from the
Navier-Stokes model for geophysical applications. These approaches are
typically employed at wide ranges of scales, from
molecular~\cite{art:Neilson2011, art:Lowengrub2016, art:nitschke2012}
to continental~\cite{art:Bouchut2004, art:Decoene-et-al2009,
  art:Bouchut2013, art:Fent2017, art:BP20} and even
planetarial~\cite{book:Holton2004,art:Higdon2006}.  These phenomena
typically act on a surface that can be a material interface such as
the cell membrane~\cite{art:Neilson2011} or the soil-water-atmosphere
interface in rivers or lakes~\cite{book:Vreugdenhil1994} or the
atmosphere~\cite{art:Decoene-et-al2009}, or modeled interfaces such as
those arising in density-dependent multi-layer shallow
water~\cite{art:bonaventura18,art:higdon22}.  In this setting, the
curvature of the interface over which the flow occurs, i.e., the bed
profile, plays a fundamental role and the governing PDEs must be
adapted to the geometrical characteristics of this interface.

In this paper we address Shallow Water (SW) models as used within the
context of geophysical flows to simulate fluid flows in natural
systems such as rivers, lakes, oceans, and the
atmosphere~\cite{art:Lanzoni2006, art:IversonGeorge2014,
  art:bonaventura18, art:higdon22}. Only few attempts at including the
bed geometry into the SW equations are
available~\cite{art:SavageHutter1991, art:Rossmanith2004,
  art:Bouchut2004, art:Bouchut2013}.  Typically, these approaches
yield systems of balance laws that are characterized by non-autonomous
flux functions, due to the presence of spatially dependent geometric
information~\cite{art:Rossmanith2004,art:Bouchut2004}, and to source
terms that contain flux variables. These two characteristics
separately contribute to the difficulties in designing accurate and
efficient numerical solvers.
The presence of non-autonomous flux functions yields discontinuous
local Riemann problems and complicate the identification of the
correct wave structure and interaction that is needed in Godunov-type
Finite Volume (FV) or Discontinuous Galerkin (DG) discretization
approaches (see, e.g.,~\cite{art:bale2003, art:andreianov2010,
  art:andreianov2011, art:xia20, art:qiao22}). Most of the approaches
used in the presence of discontinuous fluxes are tailored to the
specific cases of study and difficulties have been reported to adapt
these approaches to general first order monotone
schemes~\cite{art:qiao22}.  On the other hand, the presence of flux
variables in the source term also requires careful treatment, in
particular to ensure the well-balance of the discrete
equations~\cite{book:Bouchut2004,art:higdon22}.

Recently, starting from the three-dimensional Navier-Stokes equation,
\cite{art:BP20} developed a geometrical shallow-water model that is
defined on a local reference system anchored on the bottom surface and
that uses only intrinsic geometric quantities. The resulting covariant
formulation contains intrinsic differential operators, such as
divergence and gradients that encapsulate the relevant geometric
information.  Among the merits of this approach is that it yields a
well-posed system of hyperbolic balance laws with forcings that do not
contain flux variables.  This allows the use of time-splitting to
incorporate source terms in the numerical discretization by
geometrically intrinsic FV~\cite{art:BP20} or DG
methods~\cite{phd:Bachini20}. However, the intrinsic flux functions
are non-autonomous due to the presence of the metric tensor.  This may
cause oscillations as reported in~\cite{art:BP20} where an intrinsic
first-order Godunov FV method with HLLC Riemann solver adapted to the
geometric setting is proposed.  Other attempts to overcome this
problem in a simpler geometric setting were reported
in~\cite{art:Fent2017} where a FORCE-type scheme~\cite{art:toro09} was
use to avoid altogether the use of Riemann Solvers.

The problem of non-autonomous fluxes can be addressed by means of the
Lagrangian-Eulerian approach~\cite{art:abreu19,art:abreu17,
  art:abreu18a, art:abreu18b, art:abreu20, art:abreu21, art:abreu20b,
  phd:perez15}. In fact, in the past decades several papers have
introduced different Lagrangian-type methods, for instance
semi-Lagrangian methods~\cite{huang2016semi}, and arbitrary
Lagrangian--Eulerian methods~\cite{loubere2010reale}.
An important notion on which we base our developments is contained in
the work of~\cite{douglas2000locally, douglas2000locallyB,
  douglas2001locally}, later considered also
by~\cite{aquino2010lagrangian, huang2016semi}. The original
work~\cite{douglas2000locally} introduces the Locally Conservative
Eulerian--Lagrangian Method (LCELM), which improves the modified
method of characteristics introducing a \emph{space-time control
  volume} where local conservation is enforced. This control volume is
essentially an integral tube whose boundaries are no-flow surfaces.

In this work, we start from the geometrically intrinsic approach
of~\cite{art:BP20} and extend the Lagrangian-Eulerian method for
hyperbolic systems and balance laws developed
in~\cite{art:abreu19,art:abreu17, art:abreu18a, art:abreu18b,
  art:abreu20, art:abreu21, art:abreu20b, phd:perez15} to include
geometric effects. In short, the proposed numerical scheme is
constructed in two steps. A first Lagrangian evolution step evolves
the state from time $t^n$ to $t^{n+1}$ along the space-time control
volume appropriately defined to take into consideration also the
presence of source terms.
The second (Eulerian projection) step averages the evolved solution
over the original grid by maintaining local conservation. The fully
discrete Lagrangian-Eulerian scheme formulation discussed here is
based on the new interpretation of the integral tube, as a Lagrangian
no-flow surface~\cite{art:abreu19, art:abreu20b, art:abreu20}.
A new aspect of our Lagrangian-Eulerian method that we would like to
mention is the dynamic
forward tracking of the no-flow surfaces at each time step. This is
the key to avoid altogether the use of Riemann solvers and maintain
accuracy in the presence of non-autonomous fluxes.
Roughly speaking, our method can be viewed as some sort of combination
of Godunov-upwinding and Central differencing.  An effective
improvement over Godunov-upwind schemes is the fact that no Riemann
solvers are needed, and thus the time-consuming calculation of the
linearized characteristic directions can be avoided. Indeed, as
typical of central differencing, the Jacobian matrix of the associated
flux functions is neither constructed nor evaluated.  The forward
tracking reduces numerical dissipation and makes non-oscillatory and
high-order reconstructions redundant, thus reducing complexity of the
algorithm.
In addition, stabilizing numerical viscosity in the proposed approach
can be shown to assume the well-known $Q-$form
(see~\cite{art:abreu19}) associated to the underlying numerical flux
and satisfies the property that the numerical flux is proportional to
the ratio between the time step size and cell diameter.
This scalar value is used as an approximation of the local speed
of propagation into the space-time control volumes and can be used to
define the CFL stability constraint.

The paper is organized as follows. In Section~\ref{sec:sec2}, we
develop the proposed intrinsic Lagrangian-Eulerian approach.  First,
the geometrically Intrinsic Shallow Water Equations (ISWE) are
re-written in time-space divergence form to adapt to the development
of the intrinsic Lagrangian-Eulerian scheme taking into account the
geometric properties of the bottom surface. Next, the
Lagrangian-Eulerian scheme is developed for general cell shapes
discretizing the surface and it is then further specialized to square
cells in Section~\ref{sec:cartesian}. In this case monotonicity and
stability of the scheme can be easily
proved. Section~\ref{sec:results} reports numerical results on 
surfaces of different shapes. Finally, Section \ref{sec:conclusion}
contains the concluding remarks. 

\section{Intrinsic Lagrangian-Eulerian Finite Volumes for SWEs}
\label{sec:sec2}

\subsection{The geometrically intrinsic SWE flow system}
Consider a surface $\SurfDomain\subset\REAL^3$ equipped with a metric
$\First$ where the system of shallow water equations (SWEs) is defined.
The geometrical framework presented in~\cite{art:BP20} is employed
here and extended when needed to handle the system in a
Lagrangian-Eulerian approach.
The geometrically Intrinsic SWE (ISWE) system can be defined as:
\begin{equation}
  \DerPar{\ConservVar{}{}}{\Time}
  + \DivSurf\FluxSWE(\sv,\ConservVar{}{}) =-  \Source(\sv,\ConservVar{}{}) \,,
\label{eq1swe}
\end{equation}
subject to the Cauchy data:
\begin{equation*}
\ConservVar{}{}(0,\xcl,\ycl)=\ConservVar{0}{}\,, 
\label{eq1sweIC}
\end{equation*}
where
$\ConservVar{}{} = \{\ConservVar{}{\beta}\}_{\beta=1}^3=\left[ \Depth,\,\Qcomp[1]{},\,\Qcomp[2]{}
\right]^T$ is the conservative variable, with $\Depth$ the depth of
the normal fluid and $\Qcomp[i]{}$, $i=1,2$, the averaged velocity in
the two tangential directions. Denoting by $\sv=(\xcl,\ycl)$ the
coordinates with respect to a Local Reference System (LCS) anchored on
the bottom surface $\SurfDomain$, the differential operators ae
defined intrinsically with respect to the LCS by careful use of the
metric tensor $\First$ related to the $\sv$ local coordinates. This is
the meaning, for example, of the symbol $\DivSurf$ in~\eqref{eq1swe}
that denotes the intrinsic divergence operator (see appendix~\ref{app1}.
The hyperbolic flux function $\FluxSWE$ and the source term $\Source$ 
can be written intrinsically as:
\begin{equation*}
  \label{eq:fluxsource}
\!\!\!\FluxSWE(\sv,\ConservVar{}{}) \!=\!\! \begin{bmatrix}
  \Qcomp[1]{}   &  \Qcomp[2]{} \\[0.7em]
     \displaystyle\frac{(q^1)^2}{\Depth} \!+\!  \displaystyle\frac{g
       \Depth^2}{2 h_1^2} \displaystyle\frac{\partial x^3}{\partial
       s^3} \!\!\! & \!\!\! \frac{(q^1 q^2)}{\Depth}  \\[0.7em]
    \displaystyle\frac{q^1 q^2}{\Depth}  \!\!\!\! & \!\!\!\! \displaystyle\frac{(q^2)^2}{\Depth} \!+\!  \displaystyle\frac{g \Depth^2}{2 h_2^2} \displaystyle\frac{\partial x^3}{\partial s^3} 
\end{bmatrix}
\quad \!\!\! \mbox{ and } \!\!\! \quad
\Source(\sv,\ConservVar{}{}) \!=\!\! \begin{bmatrix}
           0 \\[0.7em]
           \displaystyle\frac{g \Depth^2}{2 h_1^2}\frac{\partial}{\partial s^1}\!\!\left( \!\displaystyle\frac{\partial x^3}{\partial s^3} \!\right)  \!+\!
							\displaystyle\frac{g \Depth}{h_1^2}\displaystyle\frac{\partial x^3}{\partial s^1} \!-\! 
							\displaystyle\frac{1}{\rho}[\DivSurf  T_{sw}]^{(1,\cdot)} \!-\! \displaystyle\frac{\tau_b^1}{\rho} \\[0.7em]
           \displaystyle\frac{g \Depth^2}{2 h_2^2}\frac{\partial}{\partial s^2}\!\!\left( \!\displaystyle\frac{\partial x^3}{\partial s^3} \!\right)  \!+\!
							\frac{g \Depth}{h_2^2}\frac{\partial x^3}{\partial s^2} \!-\! 
							\frac{1}{\rho}[\DivSurf  T_{sw}]^{(2,\cdot)} \!-\! \frac{\tau_b^2}{\rho} 
\end{bmatrix},
\end{equation*}
where $x^3$ and $s^3$ in $\tfrac{\partial x^3}{\partial s^3}$ indicate
the coordinates along the vertical and normal directions,
respectively. The tensor $T_{sw}$ is the stress tensor and the
functions $\rho$ and $\tau_b$ are the density of the fluid and the
friction to the bottom, respectively.

Our goal is to develop an Intrinsic Lagrangian-Eulerian approach for
this system of intrinsic SWE, starting by defining a \emph{generalized
  time-space divergence form} for the balance law (\ref{eq1swe}).
First, we need to extend the metric tensor in time-space to handle
time as an ``intrinsic'' variable. Thus we define:
\begin{equation*}
  \mathbb{G} :=
  \left(
    \begin{array}{rr}
      1\; \vline&  \\
      \hline
      \vline& \First 
    \end{array}
  \right)=
  \left(
    \begin{array}{ccc}
      1 & 0&0\\
      0&\NORM{\vecBaseCCcv[1](\Point)}^2 & 0\\
      0&0 &\NORM{\vecBaseCCcv[2](\Point)}^2 \\
    \end{array}
  \right)=
  \left(
    \begin{array}{ccc}
      1& 0 & 0 \\
      0&\metrcoefH{1}^2 & 0\\
      0&0 &\metrcoefH{2}^2  \\
    \end{array}
  \right)=
  \left(
    \begin{array}{ccc}
      \metrcoef{1}& 0 & 0 \\
      0&\metrcoef{2} & 0\\
      0&0 &\metrcoef{3}  \\
    \end{array}
  \right)\,,
\end{equation*}
with the associated reference frame
$\langle\normalvec, \vecBaseCCcv[1],\vecBaseCCcv[2]\rangle$, where
$\vecBaseCCcv[1],\vecBaseCCcv[2]$ span the tangent plane, while
$\normalvec$ is the vector normal to the surface and is used to
extend the tangent set to a basis of $\REAL^3$.
Then~\eqref{eq1swe} can be re-written in divergence form as (see
appendix~\ref{app1}):
\begin{align}
  \label{eq:eqSTdiv}
  \DivSurfST\Flux=-\Source,
\end{align}
where $\DivSurfST$ is the intrinsic divergence with respect to the
extended metric $\mathbb{G}$, and $\Flux$ is an extended (symmetric)
flux function defined as: 
\begin{equation}
  \label{eq:fluxtensor}
  \Flux(\Time,\sv,\ConservVar{}{}) = 
  \begin{bmatrix} 
    \Depth &\Qcomp[1]{} & \Qcomp[2]{} \\[0.5em]
    \Qcomp[1]{}& \frac{(\Qcomp[1]{})^2}{\Depth} + \frac{\grav
      \Depth^2}{2\metrcoefH{1}^2}\DerPar{\zcg}{\zcl} &
    \frac{\Qcomp[1]{}\Qcomp[2]{}}{\Depth}  \\[0.5em]
    \Qcomp[2]{} & \frac{\Qcomp[1]{}\Qcomp[2]{}}{\Depth} &
    \frac{(\Qcomp[2]{})^2}{\Depth} + \frac{\grav
      \Depth^2}{2\metrcoefH{2}^2}\DerPar{\zcg}{\zcl}
  \end{bmatrix} =
  \begin{bmatrix} 
    \Depth &\Qcomp[1]{} & \Qcomp[2]{} \\[0.5em]
    \Qcomp[1]{}& F^{21} & F^{22}  \\[0.5em]
    \Qcomp[2]{} & F^{31} & F^{32}
  \end{bmatrix}.
\end{equation}

\subsection{The Intrinsic Lagrangian-Eulerian Scheme} \quad 

\begin{figure}
  \centering
  \includegraphics[width=0.5\textwidth]{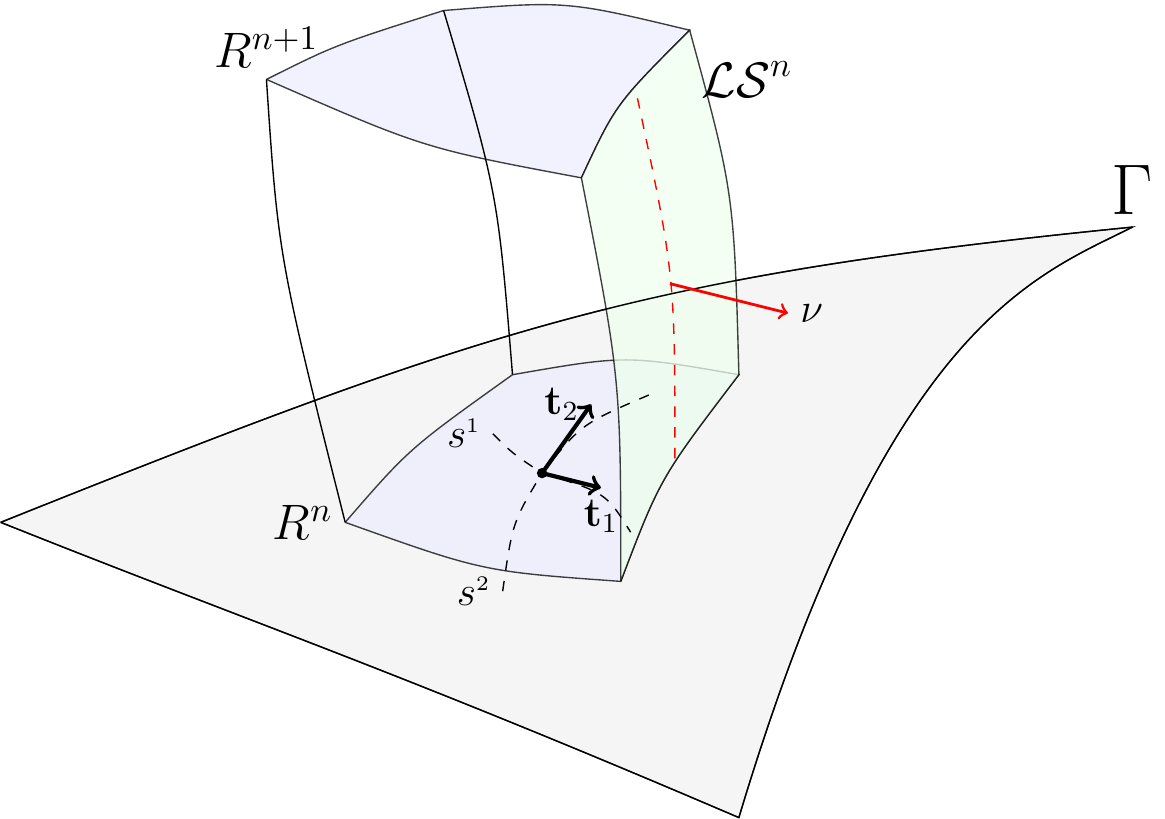}
  \caption{%
    Illustration of the control volume $\subVolume{}{n}$ in the
    time-space intrinsic local coordinates, with ${\partial \subVolume{}{n} =
    \Region[n] \cup \mathcal{LS}^n\cup \bRegion[n+1]}$,
    where
    $\Region[n]$ is the inflow region, $\bRegion[n+1]$ the outflow region and
    $\mathcal{LS}^n$ the no-flow surface.
  }
    \label{fig:FIVb}
\end{figure}

We now turn our attention to the discretization of the ISWE system
\eqref{eq:eqSTdiv} and develop our extension of the
scalar multidimensional Lagragian-Eulerian method, introduced in
\cite{art:abreu20}, to the ISWE  system written in covariant form.
As standard in a finite volume framework, we work in the extended
time-space domain consider a \emph{Discrete Control Volume} in the
Lagrangian-Eulerian local coordinate system defined by the region
$\subVolume{}{n}:=\IT[n]\times\Region\subset\REAL\times\SurfDomain$.
We will be interested in applying the divergence theorem in this
framework, and thus we need to properly define the boundary of this
control volume. In particular, with reference to
Figure~\ref{fig:FIVb},
$\partial \subVolume{}{} = \Region[n] \cup \mathcal{LS}^n\cup
\bRegion[n+1]$ is the boundary of the control volume $\subVolume{}{}$
formed by $\Region[n]$ and $\bRegion[n+1]$, the \emph{Discrete Entry}
and \emph{Discrete Exit} of the flux at times $t^n$ and $t^{n+1}$, and
the lateral surface $\mathcal{LS}^n$, called \emph{Discrete no-flow
  surface}.
Testing equation~\eqref{eq:eqSTdiv} using the characteristic function
of this discrete control volume we obtain:
\begin{align*}
  \int_{\subVolume{}{n}} \DivSurfST\Flux 
  =-\int_{\subVolume{}{n}} \Source \,.
\end{align*}
Applying the divergence theorem yields:
\begin{align}
  \int_{\partial \subVolume{}{n}} \scalprodSurfST{\Flux}{\vec{n}}
  &=-\int_{\subVolume{}{n}} \Source
    =: I(\Source)\,,
    \label{eq4}
\end{align}
where
$\scalprodSurfST{\Flux}{\vec{n}}=\left\{\Flux^{\beta}\STFirst\vec{n}\right\}^{\beta}$,
$\beta=1,2,3$, with $\Flux^{\beta}$ the $\beta$-th row of the flux function~\eqref{eq:fluxtensor}.
The vector $\vec{n}$ denotes outward normal to the boundary
$\partial \subVolume{}{}$ (see central panel in
Figure~\ref{fig:FIVb}). It corresponds to the vector normal to the
surface with $-$ sign in the entry region $\Region[n]$,
i.e. $-\normalvec$, and analogously but with the opposite sign for the
normal to the region $\bRegion[n+1]$. Note that, the vector
$\normalvec$ is equal to $[1,0,0]^{\T}$, when written in the
time-space local coordinate system. As a consequence, the
left-hand-side of equation~\eqref{eq4} can be written as:
\begin{align*}
  \int_{\Region[n]}
  \scalprodSurfST{\Flux}{\begin{bmatrix}-1\\0\\0\end{bmatrix}} 
  + \int_{\mathcal{LS}^n} \scalprodSurfST{\Flux}{\normalEdge}
  + \int_{\bRegion[n+1]}
  \scalprodSurfST{\Flux}{\begin{bmatrix}1\\0\\0\end{bmatrix}} &= I(\Source)\,,\\[1em]
  -\int_{\Region[n]} \ConservVar{}{}(t^n,\xcl,\ycl)
  + \int_{\mathcal{LS}^n} \scalprodSurfST{\Flux}{\normalEdge} 
  + \int_{\bRegion[n+1]}\ConservVar{}{}(t^{n+1},\xcl,\ycl) &= I(\Source)\,,
\end{align*}
where $\normalEdge$ is the normal to $\mathcal{LS}^n$ written in
the LCS (at each corresponding time $t$).

By virtue of the conceptual philosophy of the Lagrangian-Eulerian
approach (see, e.g., \cite{art:douglas00,art:abreu20,art:abreu21}), we
suppose that there is no-flow through the lateral surface
$\mathcal{LS}^n$, and thus we impose the no-flow condition:
\begin{equation}
  \int_{\mathcal{LS}^n} \scalprodSurfST{\Flux}{\normalEdge} = 0\,.
  \label{hypnoflow}
\end{equation}
Hence, the evolution relation for the state variable is given by:
\begin{equation}
  \int_{\bRegion[n+1]} 
  \ConservVar{}{}(t^{n+1},\xcl,\ycl) 
  \, = \int_{\Region[n]} 
  \ConservVar{}{}(t^{n},\xcl,\ycl) 
  \, + I(\Source)\,.
  \label{eq6}
\end{equation}

\paragraph{Lagrangian-Eulerian model}
With respect to the underlying geometrically intrinsic
Lagrangian-Eulerian approach, the following two problems are
equivalent in the weak (distributional) sense:
\begin{equation*}
{\bf (I)}
\begin{cases}
  & \displaystyle\DerPar{\ConservVar{}{}}{\Time}
  + \DivSurf\FluxSWE(\sv,\ConservVar{}{}) + \Source(\sv,\ConservVar{}{}) = 0 \\[0.7em] 
  &\ConservVar{}{}(0,s^1,s^2)=\ConservVar{0}{}
\end{cases}
\end{equation*}
and
\begin{equation*}
  {\bf (II)}
  \begin{cases}
    &  \displaystyle \int_{\mathcal{LS}^n} \scalprodSurfST{\Flux}{\normalEdge} = 0 \\[1em]
    & \displaystyle  \frac{1}{\CellArea[{\Region^{n+1}}]}
    \int_{\bRegion[n+1]} 
  \ConservVar{}{}(t^{n+1},\xcl,\ycl) 
  \, = \frac{1}{\CellArea[{\Region^{n+1}}]} \int_{\Region[n]} 
  \ConservVar{}{}(t^{n},\xcl,\ycl) 
  \, + \frac{1}{\CellArea[{\Region^{n+1}}]}\int_{\subVolume{}{n}} \Source\\[0.7em]
    &\ConservVar{}{}(0,s^1,s^2)=\ConservVar{0}{},
  \end{cases}
\end{equation*}
where $\CellArea[{\bRegion[n+1]}]$ denotes the area of the evolved
region $\bRegion[n+1]$.
Problem \textbf{(II)} represents the geometrically intrinsic
Lagrangian-Eulerian building block, which will form the basis for our
numerical scheme.
Recalling that the regions $\Region[n],\bRegion[n+1]$ and
$\mathcal{LS}^n$ are defined in a time-space reference frame, we can
understand  the Cauchy
problem \textbf{(II)} to represent the Lagrangian evolution
step followed by the projection onto the original region
$\Region[n]$ obtained by imposing in an Eulerian sense the local mass
conservation (Eulerian projection step).


\subsubsection{Time-space discretization}
\label{sec:approx}
We assume we are given a surface triangulation $\Triang(\SurfDomain)$
of the bottom surface $\SurfDomain$, fixed in time, formed by the
union of non-intersecting geodesic triangles with vertices on
$\SurfDomain$ (edges are geodesics), and such that
$\Triang(\SurfDomain)=\cup_{\ell=1}^{\NCell}\Cell[\ell]=\Closure{\SurfDomain}$
and $\Edge[\ell m]=\Cell[\ell]\cap\Cell[m]$ is an internal geodesic edge.
The triangulation is characterized by a mesh parameter
$\meshparam=\max_\ell \diam(\Cell[\ell])$.
Together with the space discretization, we consider a time
discretization $\{t^n\}$ of the time interval $\IT$, where we define a
non-uniform time-step $\Dt=t^{n+1}-t^n$.
At each time $t^n$, the cell $\Cell[\ell]$ coincides with the inflow
region $\Region[n]_\ell$, namely the triangulation $\Triang(\SurfDomain)$
is our working grid.
We assume that each point of $\Cell[\ell]$ can be expressed in the
tangent plane $\TanPlane[\midPoint_\ell]\SurfDomain$ passing through the
curved-cell mid-point $\midPoint_i$ with an error of
$\BigO{\meshparam^2}$, or more precisely, the distance between every
point of $\Cell[i]$ and its projection on
$\TanPlane[\midPoint_\ell]\SurfDomain$ is proportional to $\meshparam^2$.
This enables us to work on the flat cell defined by projecting the
vertices of the curved cell on $\TanPlane[\midPoint_\ell]\SurfDomain$.

We can define in each cell $\Region[n]_\ell$ the approximate 
solution by:
\begin{equation*}
 \ConservVar{\ell}{n} := \frac{1}{\CellArea[{\Region[n]_\ell}]} 
                \int_{\Region[n]_\ell} \ConservVar{}{}(t^n,\xcl,\ycl)\,,
\end{equation*}
while, the approximate solution at time $\Time^{n+1}$, i.e., the mean
value $\ConservVar{\ell}{n+1}$ over the region $\bRegion[n+1]_\ell$, is
obtained after a projection procedure (Eulerian projection step) over
the original grid $\Cell[\ell]=\Region[n]_\ell$ in the corresponding control
volume.
Denoting by $\bConservVar{\ell}{n+1}$ the approximate value at time
$\Time^{n+1}$ obtained from the Lagrangian step \eqref{eq6}
along with the desired conservation properties, we can write:
\begin{multline}
  \bConservVar{\ell}{n+1} = \frac{1}{\CellArea[{\bRegion[n+1]_\ell}]}
  \int_{\bRegion[n+1]_\ell} \ConservVar{}{}(t^{n+1},\sv) \\
  = 
	\frac{\CellArea[{\Region[n]_\ell}]}{\CellArea[{\bRegion[n+1]_\ell}]}   
\left[ \frac{1}{\CellArea[{\Region[n]_\ell}]} \int_{\Region[n]_\ell}
  \ConservVar{}{}(t^{n},\sv)
  + \frac{1}{\CellArea[{\Region[n]_\ell}]} I(\Source) \right] = 
 \frac{\CellArea[{\Region[n]_\ell}]}{\CellArea[{\bRegion[n+1]_\ell}]} 
\left[\ConservVar{\ell}{n} + \frac{1}{\CellArea[{\Region[n]_\ell}]} I(\Source)\right]\,.
  \label{eq7}
\end{multline}

\paragraph{Approximation of the no-flow condition~\eqref{hypnoflow}}

\begin{figure}
\centering
\includegraphics[width=0.48\textwidth]{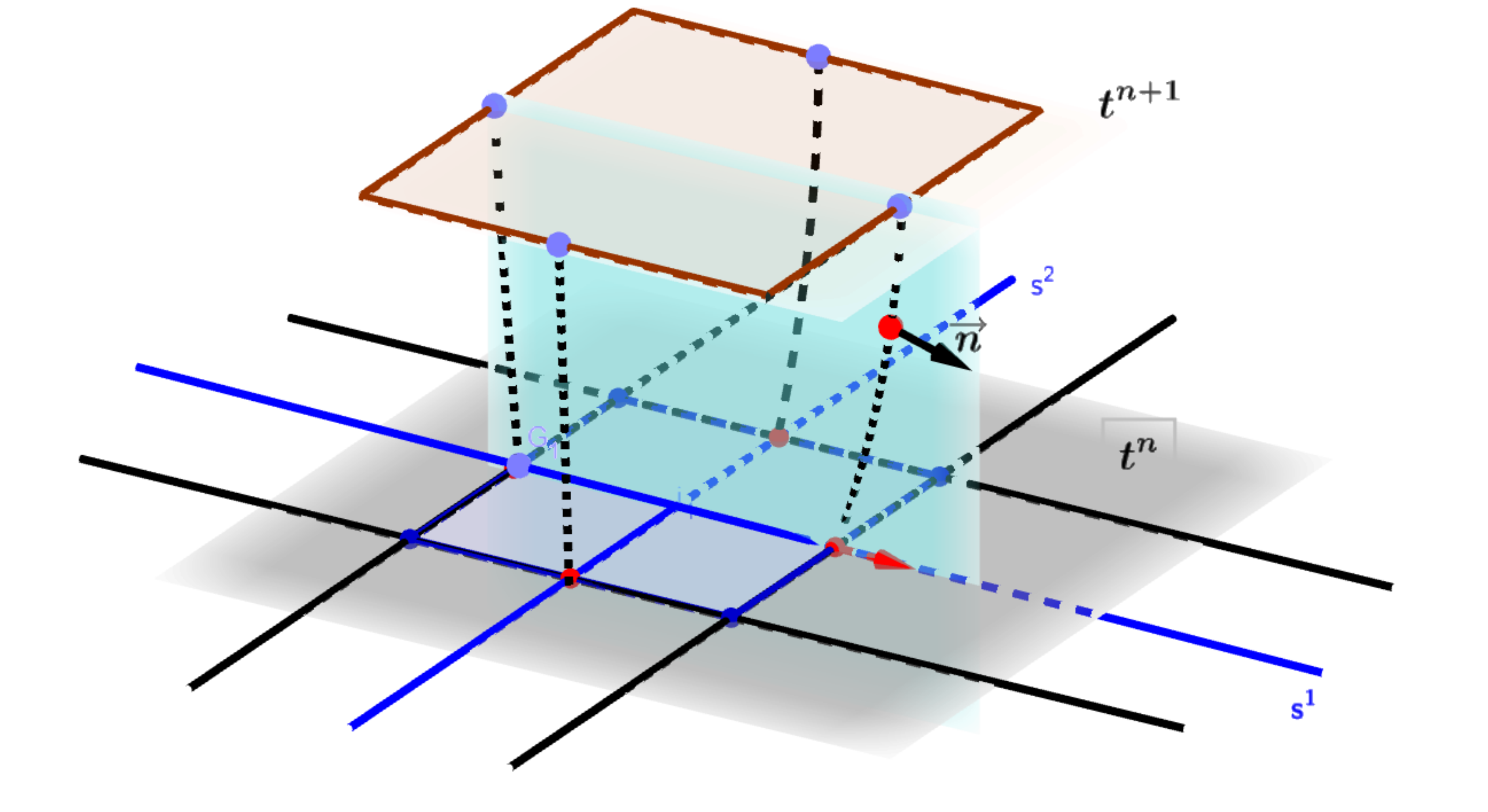}
\hspace{10pt}
\includegraphics[width=0.48\textwidth]{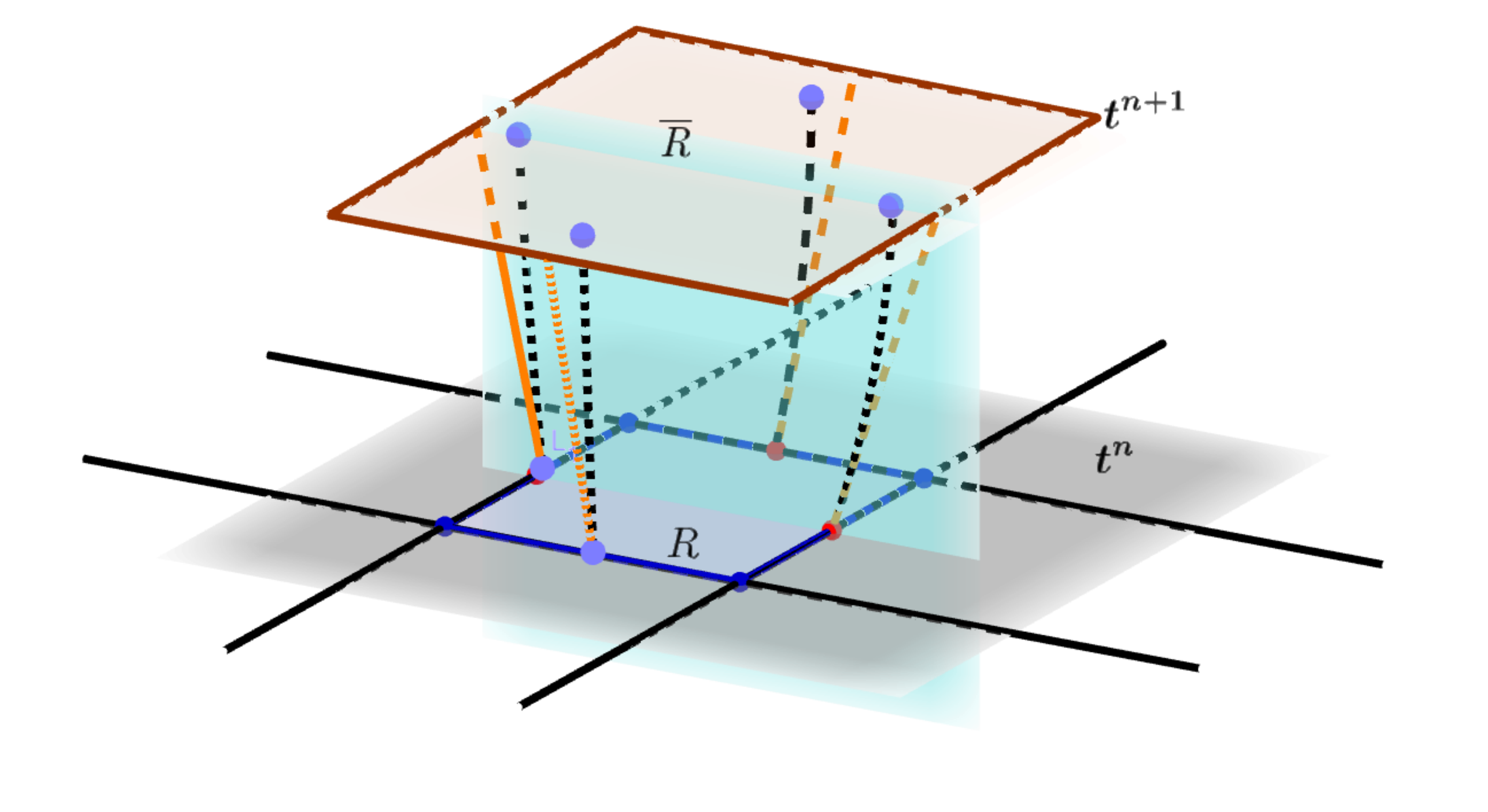}
\caption{Approximation of the no-flow surface.}
\label{fig:shortGLE}
\end{figure}

In this paragraph we describe how to approximate
equation~\eqref{hypnoflow} on each face $\mathcal{LS}_{\ell m}^n$ forming
the lateral no-flow surface $\mathcal{LS}_\ell^n$.  For any time
$\tau\in\IT[n]$, the equation can be re-written separating space and
time as:
\begin{equation*}
  0=
  \int_{t^n}^{\tau} \int_{\Edge[\ell m](t)}
  \scalprodSurfST{\Flux(t,\sv)}{\normalEdge(t,\sv)} \Diff \sv\;\Diff t
  \approx \int_{t^n}^{\tau} \edgeLength[\ell m](t)
  \scalprodSurfST{\Flux(t,\midPointEdge)}{\normalEdge(t,\midPointEdge)}
  \Diff t\,,
\end{equation*}
where $\edgeLength[\ell m]$ and $\midPointEdge$ are the length and the
mid-point of $\Edge[\ell m](t)$, respectively.  Note that we used the
mid-point rule for the space integral suggesting that the consistent
characterization of the no-flow surface can be obtained by
imposing:
\begin{equation*}
  \label{eq:noflow-argument}
  \scalprodSurfST{\Flux(t,\midPointEdge)}%
  {\normalEdge(t,\midPointEdge)}=\mathbf{0}
  \qquad\qquad \Forall t\in\IT[n]\,,
\end{equation*}
Now, we can identify the no-flow curve
$\sigma_{\ell m}\subset \mathcal{LS}_{\ell m}^n$ that starts from
$\midPointEdge$ and write its parametrization as (see
figure~\ref{fig:shortGLE}, left):
\begin{equation*}
  \label{eq:sigmat}
  \sigma_{\ell m}(t) = [t,\sigma^1_{\ell m}(t),\sigma^2_{\ell m}(t)]\,.
\end{equation*}
The normal $\normalEdge_{\ell
  m}(\Time)=\normalEdge(t,\midPointEdge)$
can be written implicitly in terms of $\sigma'_{\ell m}(t)$ in
the time-space local reference system. This yields a nonlinear system of
ODEs that describes $\sigma'_{\ell m}$ as a function of the SWE flux, with
initial conditions given by $\sigma_{\ell m}(t^n) =
[t^n,\midPointEdge^1(t^n),\midPointEdge^2(t^n)]^{\T}$
and $\scalprodSurfST{\normalEdge_{\ell m}(t^n)}{\sigma_{\ell m}'(t^n)}=0$.
The system simplifies considerably if we linearize the equation around
$t^n$, which is equivalent to apply explicit Euler to solve the above
ODE (see figure~\ref{fig:shortGLE}, right).  Following this approach,
the position of the evolved mid-point at $t\in\IT[n]$ can be
calculated as:
\begin{equation*}
 \sigma_{\ell m}(t) = \sigma_{\ell m}(t^n) +
 (t-t^n)\sigma_{\ell m}'(t^n)\,.
\end{equation*}
For example, in the case of an edge with the unit tangent vector in
$\midPointEdge(t^n)$ that is aligned with $\vecBaseCCcv[2]$, we 
obtain explicitly $\sigma_{\ell m}'=[1,(\sigma^1_{\ell m})',0]^{\T}$ by
solving the following linear system:
\begin{equation}
0 = \scalprodSurfST{\Flux(t^n,\midPointEdge)}{\normalEdge(t^n,\midPointEdge)}\Big|^{\beta}=
\Flux^{\beta}(t^n,\midPointEdge)\, \STFirst(\midPointEdge)\,\begin{bmatrix}-1\\
\displaystyle\frac{1}{(\sigma_{\ell m}^{1\beta})'(t^n)}\\
0
\end{bmatrix}
\,,
\label{odecp1}
\end{equation}
for each component $\beta=1,2,3$.

\paragraph{Discrete Lagrangian-Eulerian scheme}
We can summarize the fully discrete version of problem \textbf{(II)}
forming our proposed Lagrangian-Eulerian finite volume scheme
as follows: 
\begin{equation*}
  \begin{cases}
    & \displaystyle\sigma_{\ell m}(t^{n+1}) = \sigma_{\ell m}(t^n) +
 \Delta t\; \sigma_{\ell m}'(t^n)\,, \qquad 
    \mbox{where} \\[0.5em]
&\displaystyle\scalprodSurfST{\normalEdge_{\ell m}(t^n)}{\sigma_{\ell m}'(t^n)}=0 \quad
    \mbox{ and }\quad
    \scalprodSurfST{\Flux(t^n,\midPointEdge)}{\normalEdge(t^n,\midPointEdge)}=\mathbf{0}
\quad \mbox{ and } \quad
    \sigma_{\Edge[\ell m]}(t^n) =
      [t^n,\midPointEdge^1,\midPointEdge^2]^{\T}
      \\[0.7em]
      &\displaystyle\bConservVar{\ell}{n+1} = \frac{\CellArea[{\Region[n]_\ell}]}{\CellArea[{\bRegion[n+1]_\ell}]} 
\left[\ConservVar{\ell}{n} + \frac{1}{\CellArea[{\Region[n]_\ell}]} \int_{\Region[n]_\ell}\Source\right]\,,
  \\[1em]
    &\ConservVar{\ell}{}(0,\midPoint_\ell^1,\midPoint_\ell^2)=\ConservVar{0}{}.
  \end{cases}
\end{equation*}
The numerical formulation is completed by the Eulerian projection
step, which, depending on the geometry of the cells, is
described in the next section.
The approach is simple and robust in terms of accuracy and
computational cost, since simple quadrature formulae as such
trapezoidal and midpoint rules can be employed  for the
approximation of the source term $\Source$ (see,
e.g.,~\cite{art:abreu19, art:abreu17, art:abreu21, art:abreu22b}).

\section{Special case: Cartesian grid}
\label{sec:cartesian}

\begin{figure}
\centering
\includegraphics[height=0.1972\textheight,
  width=0.2472\textwidth]{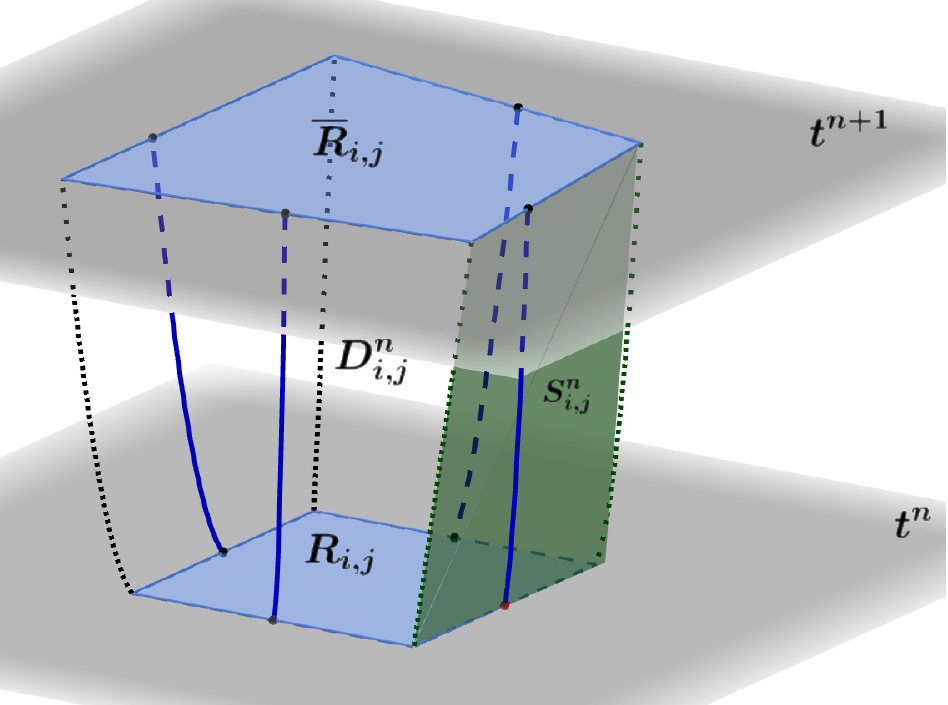}
\hspace{5pt}
\includegraphics[scale=0.40]{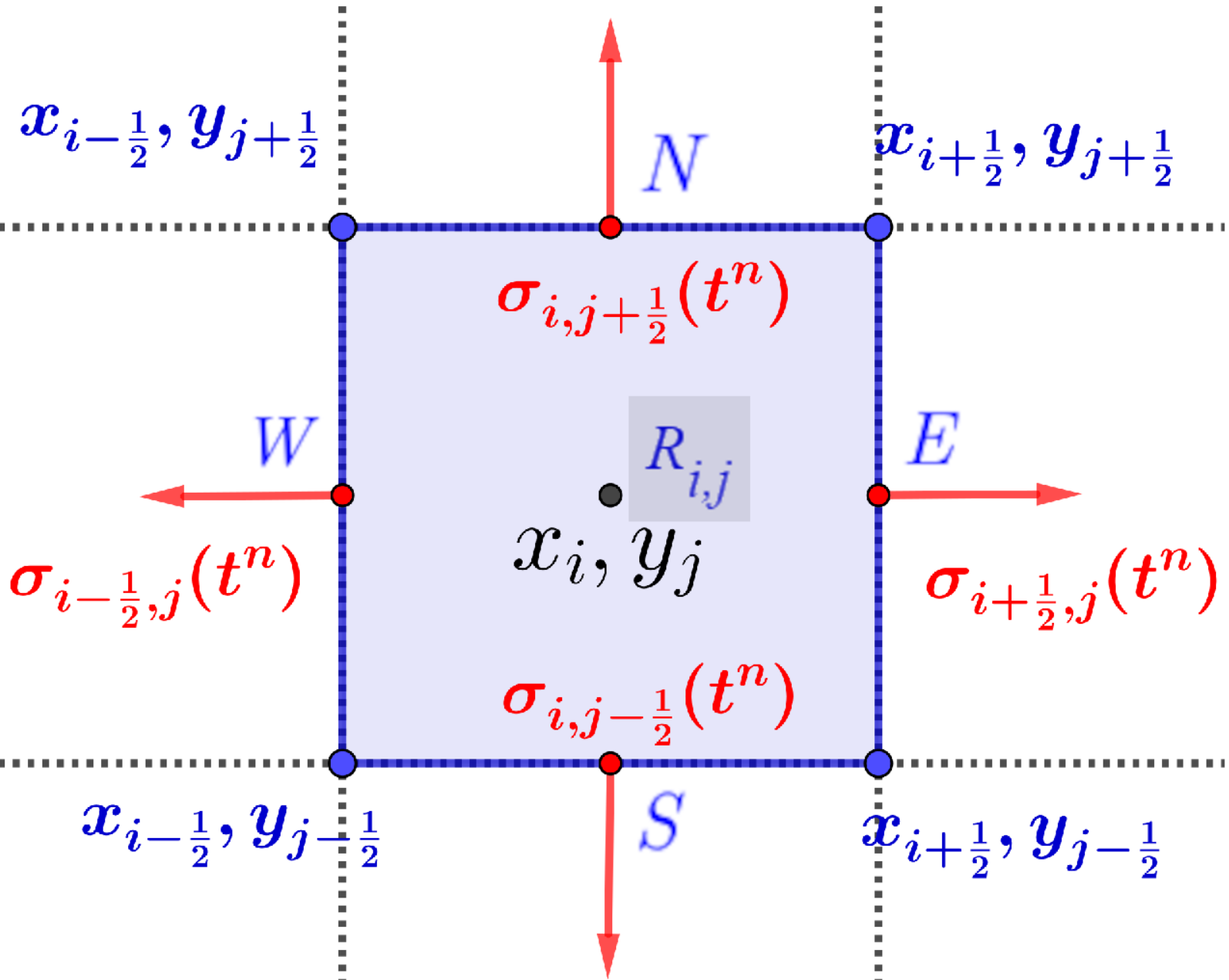} 
\caption{Notation in the Cartesian case. No-flow curves
  ${\sigma_{i-\frac{1}{2},j}(t)}$, ${\sigma_{i+\frac{1}{2},j}(t)}$,
  ${\sigma_{i,j-\frac{1}{2}}(t)}$ and ${\sigma_{i,j+\frac{1}{2}}(t)}$
  for each edge.}
\label{fig:LSnoFLOW}
\end{figure}

In this section, we consider a special tringulation
$\Triang(\SurfDomain)$ of $\SurfDomain$ defined from a regular square
Cartesian grid of the $\REAL[2]$ chart by vertical projection of the
vertices onto the surface. Thus, we can use local flat coordinates
defined on $\TanPlane[\midPoint_\ell]\SurfDomain$ for each cell
$\Cell[\ell]$. We will denote by $(x,y)$ the local flat coordinates.
We use the lattice $\mathds{N}\times\left(\mathds{Z} \times
\mathds{Z}\right) = \{(n,i,j) : n = 0,1,2,\ldots, \,\, i,j = 0, \pm 1,
\pm 2,\ldots \}$ to describe the Cartesian grid of the chart for each
time $t^n$.
We denote by $h_{x}^n = \Delta x^n = x_{i + \frac{1}{2}}^n - x_{i
  - \frac{1}{2}}^n$, and $h_{y}^n = \Delta y^n = y_{j +
  \frac{1}{2}}^n - y_{j - \frac{1}{2}}^n$, where $(x_{i \pm
  \frac{1}{2}}^n,y_{j \pm \frac{1}{2}}^n)$ are the corners of the
$(i,j)$-cell defined on $\TanPlane[\midPoint_i]\SurfDomain$.
The cell centers and vertices are then given by:
\begin{equation*}
(x_i^n,y_j^n) = {(i h_x^n,jh_y^n)} \qquad
\mbox{ and } \qquad
\left(x_{i \pm \frac{1}{2}}^n,y_{j \pm \frac{1}{2}}^n\right) = 
\left(i h_x^n  \pm \frac{h_x^n}{2},j h_y^n\pm \frac{h_y^n}{2}\right)\,.
\end{equation*}
The Lagrangian-Eulerian control volume $\subVolume{i,j}{n}$ (figure
\ref{fig:LSnoFLOW}) is the volume comprised by the union of the 
surface
$\partial\subVolume{i,j}{n}=\Region[n]_{i,j}
\cup\mathcal{LS}_{i,j}^n\cup\bRegion[n+1]_{i,j}$
described in the LCS mesh grid by: 
\begin{itemize}
\item $\Region[n]_{i,j} = \left[x_{i-\frac{1}{2}}^n,x_{i+\frac{1}{2}}^n\right]
  \times 
  \left[y_{j-\frac{1}{2}}^n, y_{j+\frac{1}{2}}^n\right]$ 
  the inflow of $\partial \subVolume{i,j}{n}$;
\item $\bRegion[n+1]_{i,j} = 
  \left[\overline{x}_{i-\frac{1}{2}}^{n+1},\overline{x}_{i+\frac{1}{2}}^{n+1}\right]
  \times\left[\overline{y}_{j-\frac{1}{2}}^{n+1},\overline{y}_{j+\frac{1}{2}}^{n+1}\right]$
  the outflow of $\partial \subVolume{i,j}{n}$;
\item  $\mathcal{LS}_{i,j}^n$ the lateral no-flow surface, formed
  by North (\emph{N}), East (\emph{E}), South (\emph{S}), and West
  (\emph{W}) edges as indicated in figure~\ref{fig:LSnoFLOW}. 
\end{itemize}
The area of every computational square cell
$\left[x_{i - \frac{1}{2}}^n,x_{i + \frac{1}{2}}^n \right] \times 
\left[y_{j - \frac{1}{2}}^n,y_{j + \frac{1}{2}}^n\right]$ (see figure
\ref{fig:LSnoFLOW}, right) is $A(\Region[n]_{i,j})=h_x^nh_y^n$, so that
the approximate solution is:
\begin{equation*}
  \ConservVar{i,j}{n} = \frac{1}{A(\Region[n]_{i,j})} 
                \int_{x_{i - \frac{1}{2}}^n}^{x_{i 
               + \frac{1}{2}}^n} \int_{y_{j - \frac{1}{2}}^n}^{y_{j 
               + \frac{1}{2}}^n}  \ConservVar{}{}(t^n,x,y)\,dx\,dy, 
 \label{2DAp1}
\end{equation*}
with analogous expression for the average solution
$\bConservVar{}{}(t^{n+1},x_i,y_j)$ in the evolved region
$\bRegion[n+1]$.  We recall that the approximate solution $
\ConservVar{i,j}{n+1}$ at time
$t^{n+1}$ is obtained after the projection procedure (Eulerian
projection step) over the original grid (see figure \ref{fig:LSnoFLOW},
left).
Using this notation, equation~\eqref{eq7} provides the expression of
our finite volume scheme in this Cartesian setting.

\paragraph{Approximation of the no-flow surfaces $\mathcal{LS}_{i,j}^n$}
The approximation of the integrals on the no-flow surfaces follows the
procedure described in section~\ref{sec:approx}.
For instance, the parametrization of the no-flow curve for the edge
\emph{N} is indicated by
$\sigma_{N}(t)=\sigma_{i,j+\frac{1}{2}}(t)=
[t,\sigma_{N}^{1}(t),\sigma_{N}^{2}(t)]$,
and the normal vector at $t^n$ is given by $\normalEdge =
\left[-1,0,\displaystyle\frac{1}{(\sigma_N^{2})'(t^n)}\right]^{\T}$.
Then, equation~\eqref{odecp1} becomes:
\begin{equation*}
0 = \scalprodSurfST{\Flux}{\normalEdge}\Big|^{\beta} = \Flux^{\beta}\, \STFirst\,
\begin{bmatrix}-1\\0\\\displaystyle\frac{1}{(\sigma_{N}^{2\beta})'(t^n)}\end{bmatrix},
\end{equation*}
resulting in:
\begin{equation*}
  (\sigma_{N}^{2\beta})'(t^n) =
  \Flux^{(\beta,3)}\metrcoefH{2}^2/\ConservVar{}{\beta,n} \,,
  \qquad \beta = 1,2,3.
\end{equation*}
On the other hand, fo edge \emph{W}, the normal vector at $t^n$ is
given by
$\normalEdge=\left[-1,\displaystyle\frac{1}{(\sigma_{W}^{2})'(t^n)},0\right]^{\T}$,
and equation~\eqref{odecp1} reads:
\begin{equation*}
  0 = \scalprodSurfST{\Flux}{\normalEdge}\Big|^{\beta}
  = \Flux^{\beta}\, \STFirst\,
\begin{bmatrix}-1\\\displaystyle\frac{1}{(\sigma_{W}^{1\beta})'(t^n)}\\0\end{bmatrix}
\,,
\end{equation*}
thus:
\begin{equation*}
  (\sigma_{W}^{1\beta})'(t^n) =
  \Flux^{(\beta,2)}\metrcoefH{1}^2/\ConservVar{}{\beta, n} \,,
  \qquad \beta = 1,2,3\,.
\end{equation*}
In summary, we have the following four expressions for each equation
of the discrete ISWE system for $\beta = 1,2,3$ associated to each
sides \emph{N}, \emph{E}, \emph{S} and \emph{W}, respectively:
\begin{equation*}
(\sigma_{N}^{2\beta})'(t^n) = {\frac{\Flux^{(\beta,3)}\metrcoefH{2}^2}{\ConservVar{}{\beta,n}}}, 
\qquad
(\sigma_{E}^{1\beta})'(t) = {\frac{\Flux^{(\beta,2)}\metrcoefH{1}^2}{\ConservVar{}{\beta,n}}},
\qquad
(\sigma_{S}^{2\beta})'(t) = {\frac{\Flux^{(\beta,3)}\metrcoefH{2}^2}{\ConservVar{}{\beta,n}}},
\qquad
(\sigma_{W}^{1\beta})'(t) = {\frac{\Flux^{(\beta,2)}\metrcoefH{1}^2}{\ConservVar{}{\beta,n}}}\,.
\end{equation*}
Finally, for each edge, the linearization in time reads as follows:
\begin{equation}
  \sigma_{i-\frac{1}{2},j}^{\beta}(t)
  \approx {x_{i-\frac{1}{2}}} + (t - t^n){\frac{\Flux_{W}^{(\beta,2)}\metrcoefH{1}^2}{\ConservVar{}{\beta,n}}}
\quad \text{ and } \quad 
\sigma_{i+\frac{1}{2},j}^{\beta}(t) \approx x_{i+\frac{1}{2}} + (t - t^n){\frac{\Flux_{E}^{(\beta,2)}\metrcoefH{1}^2}{\ConservVar{}{\beta,n}}},
\label{WEodesNF}
\end{equation}
\begin{equation}
\sigma_{i,j-\frac{1}{2}}^{\beta}(t) \approx y_{j-\frac{1}{2}} + (t - t^n){\frac{\Flux_{S}^{(\beta,3)}\metrcoefH{2}^2}{\ConservVar{}{\beta,n}}}
\quad \text{ and } \quad 
\sigma_{i,j+\frac{1}{2}}^{\beta}(t) \approx y_{j+\frac{1}{2}} + (t - t^n){\frac{\Flux_{N}^{(\beta,3)}\metrcoefH{2}^2}{\ConservVar{}{\beta,n}}},
\label{WEodesNS}
\end{equation}
for $t^n<t\le t^{n+1}$.

\begin{figure}
\centering
\includegraphics[height=0.203125\textheight, width=0.253125\textwidth]{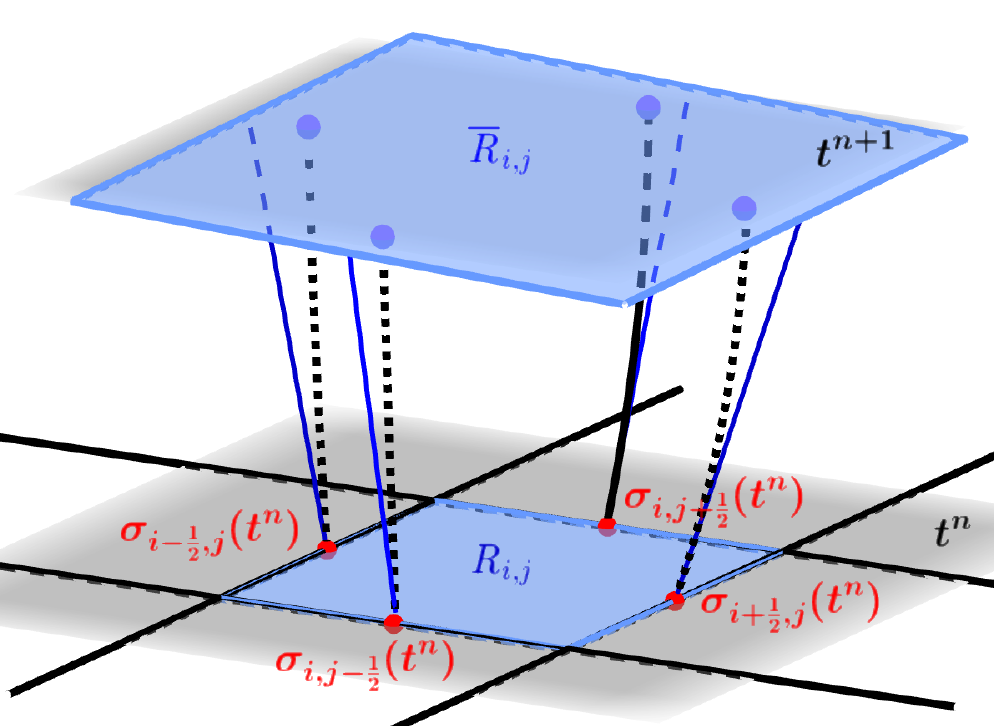}
\hspace{10pt}
\includegraphics[width=0.30\textwidth]{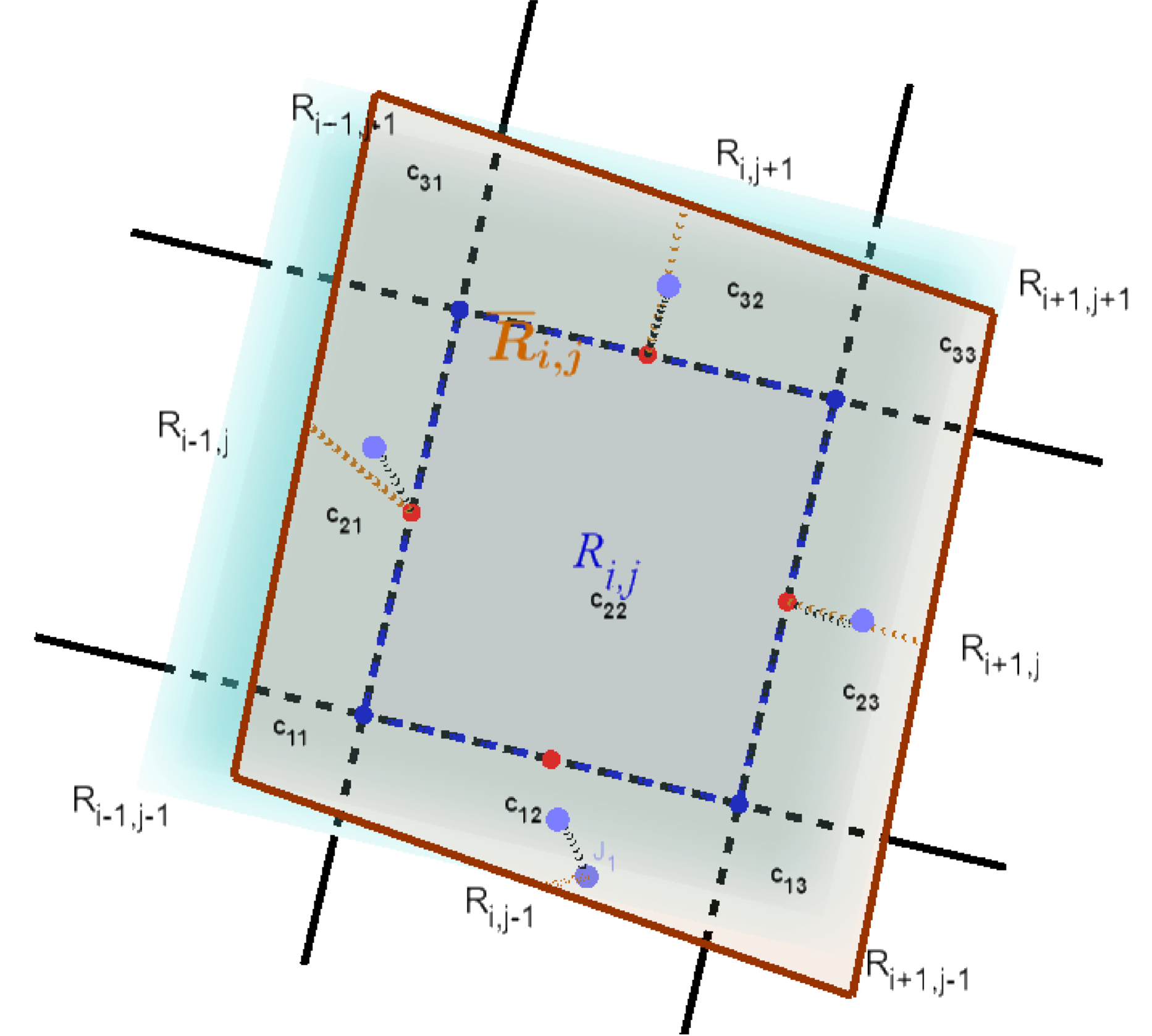}
\caption{Summary of the Geometrically Intrinsic Lagrangian-Eulerian
  construction in the Cartesian case.}
\label{fig:NOFLOWPROJ}
\end{figure}

The geometrically intrinsic Lagrangian-Eulerian building block
on a regular Cartesian grid is summarized as follows:
\begin{itemize}
\item
  {\bf STEP I} (Lagrangian evolution step)
  \begin{equation}
    \begin{cases}
  &\displaystyle \sigma_{i\pm\frac{1}{2},j}^{\beta}(t^{n+1})
  = {x_{i\pm\frac{1}{2}}} + \Delta t \;{\frac{\Flux_{E/W}^{(\beta,2)}\metrcoefH{1}^2}{\ConservVar{}{\beta,n}}}
  \qquad\mbox{
    and } \qquad
\displaystyle\sigma_{i,j\pm\frac{1}{2}}^{\beta}(t^{n+1}) = y_{j\pm\frac{1}{2}} +
  \Delta t\;{\frac{\Flux_{N/S}^{(\beta,3)}\metrcoefH{2}^2}{\ConservVar{}{\beta,n}}} 
\\[1.5em]
  &\displaystyle \bConservVar{i,j}{n+1}=
  \frac{A(\Region[n]_{i,j})}{A(\bRegion[{n+1}]_{i,j})}
  \left[ \ConservVar{i,j}{n}
    + \frac{1}{A(\Region[n]_{i,j})}
    \int_{\Region[n]_{i,j}}\Source^{\beta}\right]
  \\[1.5em] 
  & \displaystyle \ConservVar{i,j}{\beta}(0)=\ConservVar{0}{\beta}
    \end{cases}
  \label{ReltnTnp2}
\end{equation}
for all equations $\beta=1,2,3$, with
\begin{equation*}
A(\bRegion[n+1]_{ij})=
\left[h^n_x - \left({\frac{\Flux_{E}^{(\beta,2)}\metrcoefH{1}^2}{\ConservVar{}{\beta,n}}}
  - {\frac{\Flux_{W}^{(\beta,2)}\metrcoefH{1}^2}{\ConservVar{}{\beta,n}}}\right)\Delta t\right]\,
\left[h^n_y - \left({\frac{\Flux_{N}^{(\beta,3)}\metrcoefH{2}^2}{\ConservVar{}{\beta,n}}} 
  - {\frac{\Flux_{S}^{(\beta,3)}\metrcoefH{2}^2}{\ConservVar{}{\beta,n}}}\right)\Delta t\right]\,.
\end{equation*}
      \item
{\bf STEP II} (Eulerian projection step) 
\begin{equation}
\ConservVar{i,j}{n+1} = \frac{1}{A(\bRegion[{n+1}]_{i,j})} (C1+C2+C3),
\label{proy}	
\end{equation}
with
\begin{align*}
C1= c_{11}\bConservVar{i-1,j-1}{n+1} +
c_{12}\bConservVar{i,j-1}{n+1} +
c_{13}\bConservVar{i+1,j-1}{n+1}\,,\qquad
C2=c_{21}\bConservVar{i-1,j}{n+1} + c_{22}\bConservVar{i,j}{n+1} +
c_{23}\bConservVar{i+1,j}{n+1}\,,
\end{align*}
\begin{align*}
C3=
c_{31}\bConservVar{i-1,j+1}{n+1} + c_{32}\bConservVar{i,j+1}{n+1}
+ c_{33}\bConservVar{i+1,j+1}{n+1}\,,
\end{align*}
where
the coefficients $c_{ij}$ are defined 
as the entries of the matrix (see figure \ref{fig:NOFLOWPROJ})
\begin{equation*}
C = (c_{ij}) = C_x^T  C_y, \quad i,j \in \{1,2,3\}, 
\label{matrixcfl}
\end{equation*}
under the CFL-condition:
\begin{equation}
\max_{i,j, \beta} 
\left\{
\left|{\frac{\Flux_{N}^{(\beta,3)}\metrcoefH{2}^2}{\ConservVar{}{\beta,n}}} \right|,
\left|{\frac{\Flux_{E}^{(\beta,2)}\metrcoefH{1}^2}{\ConservVar{}{\beta,n}}} \right|,
\left|{\frac{\Flux_{S}^{(\beta,3)}\metrcoefH{2}^2}{\ConservVar{}{\beta,n}}} \right|,
\left|{\frac{\Flux_{W}^{(\beta,2)}\metrcoefH{1}^2}{\ConservVar{}{\beta,n}}} \right|
\right\} 
\frac{\Delta t}{h} 
< \frac{1}{2},
\label{matrixcfl2}
\end{equation}
 and the
assumption $h =\Delta x = \Delta y$.
The expression for $C_x$ and $C_y$ takes on the form:
\begin{equation}
  C_x = \left[C_{xl}, {\Delta x} - C_{xl} - C_{xr},C_{xr}\right], \qquad
  \mbox{ and } \qquad
  C_y = \left[C_{yr}, {\Delta y} - C_{yr} - C_{yr},C_{yr}\right],
\label{coefXY}
\end{equation}
assuming:
\begin{align*}
C_{xl} &= 0.5\left[1 + \mbox{sign}\left({\frac{\Flux_{W}^{(\beta,2)}\metrcoefH{1}^2}{\ConservVar{}{\beta,n}}}\right)\right]{\frac{\Flux_{W}^{(\beta,2)}\metrcoefH{1}^2}{\ConservVar{}{\beta,n}}}\Delta t\,, \qquad
C_{xr} = 0.5\left[1 - \mbox{sign}\left({\frac{\Flux_{E}^{(\beta,2)}\metrcoefH{1}^2}{\ConservVar{}{\beta,n}}}\right)\right]{\frac{\Flux_{E}^{(\beta,2)}\metrcoefH{1}^2}{\ConservVar{}{\beta,n}}}\Delta t\,, \\
C_{yl} &= 0.5\left[1 + \mbox{sign}\left({\frac{\Flux_{S}^{(\beta,3)}\metrcoefH{2}^2}{\ConservVar{}{\beta,n}}}\right)\right]{\frac{\Flux_{S}^{(\beta,3)}\metrcoefH{2}^2}{\ConservVar{}{\beta,n}}}\Delta t\,, \qquad
C_{yr} = 0.5\left[1 -
  \mbox{sign}\left({\frac{\Flux_{N}^{(\beta,3)}\metrcoefH{2}^2}{\ConservVar{}{\beta,n}}}\right)\right]{\frac{\Flux_{N}^{(\beta,3)}\metrcoefH{2}^2}{\ConservVar{}{\beta,n}}}\Delta
t\,.
\end{align*}
\end{itemize}

\subsection{Improved spatial resolution of the fully-discrete
  Lagrangian-Eulerian scheme}
\label{ImODE}

Improved spatial resolution can be obtained by means of a linear
reconstruction of the quantities $\ConservVar{i,j}{n}$ in
\eqref{ReltnTnp2} and $\ConservVar{i,j}{n+1}$ in \eqref{proy}, which
translates in a better approximation of the discrete no-flow surfaces.
Following \cite{art:abreu16,art:abreu20,art:abreu20b}, the piecewise
constant numerical data is reconstructed into a piecewise linear
approximation through the use of MUSCL-type interpolants:
$\mathbf{L}_{i,j}(t,x) = \mathbf{U}_{i,j}(t)+(x-x_j) \mathbf{U}'_{i,j}/h$, in the
$x$-direction. The interpolation in the
$y$-direction is similar.
For the
numerical derivative $\mathbf{U}'_{i,j}/h$, there are several
choices of slope limiters for scalar case, and \cite{book:Leveque2002}
contains a good compilation of possible options.
Let us consider the parametric no-flow curve
$\sigma_{i-\frac{1}{2},j}(t)$, the other cases are similar. A
straightforward calculation leads to:
\begin{equation}
\mathbf{U}_{i-\frac{1}{2},j} =
\frac{1}{h}\left[ \int_{x_{i-1,j}^n}^{x_{i-\frac{1}{2},j}^n} 
\mathbf{L}_{i-1,j}(t,x) dx
+ \int_{x_{i-\frac{1}{2},j}^n}^{x_{i,j}^n} 
\mathbf{L}_{i,j}(t,x) dx \right] =
\displaystyle\frac{1}{2} 
(\mathbf{U}_{i-1,j} + \mathbf{U}_{i,j}) + \frac{1}{8} (\mathbf{U}'_{i,j} - \mathbf{U}'_{i-1,j}).
\end{equation}

Finally, in order to show the flexibility 
of the reconstruction, we use the nonlinear Lagrange 
polynomial in $\mathbf{U}_{i-1,j}^n$, $\mathbf{U}_{i,j-1}^n$, $\mathbf{U}_{i,j}^n$, 
$\mathbf{U}_{i,j+1}^n$ and $\mathbf{U}_{i+1,j}^n$. Therefore, equation 
(\ref{ReltnTnp2}) reads:
\begin{equation*}
\overline{\mathbf{U}}_{i,j}^{n+1}  =  
\displaystyle\frac{1}{h}
\displaystyle\int_{w_{j-\frac{1}{2}}^{n}}^{w_{j+\frac{1}{2}}^{n}}\mathbf{P}_2(x,y)dw,
\quad 
\text{where $w = x,y$}
\label{Reco}
\end{equation*}
and
\begin{multline*}
\mathbf{P}_2(x,y) \!=\! 
\mathbf{U}_{i-1,j}^n  \mathbf{L}_{-1}(x - x_{i}) + 
\mathbf{U}_{i,j-1}^n  \mathbf{L}_{-1}(y - y_{j})\\ + 
\mathbf{U}_{i,j}^n \mathbf{L}_{0}(x - x_{i}) 
+\mathbf{U}_{i,j}^n  \mathbf{L}_{0}(y - y_{j}) + 
\mathbf{U}_{i+1,j}^n \mathbf{L}_1(x - x_{j}) + 
\mathbf{U}_{i,j+1}^n \mathbf{L}_1(y - y_{j})\,,
\end{multline*}
with
\begin{equation*}
\mathbf{L}_{\pm \,1}(x) = \frac{1}{2}\left[\left(\frac{x}{h} 
\pm \frac{1}{2}\right)^2
- \frac{1}{4}\right], \quad \quad \quad  
\mathbf{L}_0(x) = 1- \left(\frac{x}{h}\right)^2.
\label{Lagran}
\end{equation*}

\subsection{A connection with monotone finite difference (volume)
  schemes}
\label{monoLE} 

The fact that each working cell is assumed to be a subset of the
tangent plane $\TanPlane[\midPoint_{\ell}]\SurfDomain$ allows us to
recast our proposed scheme within the framework of monotone finite
difference schemes (see
\cite{art:chainais99,art:eymard95a,art:eymard95b}). This theory is
developed only for the case of
scalar hyperbolic equations, thus in the following we consider $\frac{\partial U}{\partial t} + \Div
\FluxSWE(U) = 0$, where $\FluxSWE(U)=[\mathcal{F}_1,\mathcal{F}_2]$.
Under the CFL stability constraint~\eqref{matrixcfl2} and
switching-off the linear reconstruction of the quantities
$U_{i,j}^{n},U_{i,j}^{n+1}$, the fully-discrete Lagrangian-Eulerian
scheme (\ref{ReltnTnp2})-(\ref{proy}) can be written as:
\begin{equation*}	
U_{i,j}^{n+1} = C_x A_{i,j}^{T} C_y^T,
\qquad
A_{i,j} = \left[\begin{array}{lll}
\overline{U}_{i-1,j-1}^{n+1} & \overline{U}_{i,j-1}^{n+1} 
& \overline{U}_{i+1,j-1}^{n+1} \\
\overline{U}_{i-1,j}^{n+1}   & \overline{U}_{i,j}^{n+1}   
& \overline{U}_{i+1,j}^{n+1} \\
\overline{U}_{i-1,j+1}^{n+1} & \overline{U}_{i,j+1}^{n+1} 
& \overline{U}_{i+1,j+1}^{n+1}
\end{array}\right],
\label{proyF}
\end{equation*}
or, in the classical form of conservative monotone
scheme~\cite{art:abreu20, art:crandall80}, as:
\begin{equation}	
U_{i,j}^{n+1} = U_{i,j}^n - \lambda ^x \Delta _{+}^{x} 
     F(U_{i-1,j-1}^n,...,U_{i+1,j+1}^n) 
     - \lambda ^y \Delta _{+}^{y} 
      G(U_{i-1,j-1}^n,...,U_{i+1,j+1}^n),.
\label{proyF1}
\end{equation}
Here, using the standard notations, we define
 $h = \Delta x = \Delta y$, and
$\lambda^{x} = \Delta t/ \Delta x$, 
$\lambda^{y} = \Delta t/ \Delta y$, and
$(\Delta^{x}_{+}U)_{j,k} = U_{j+1,k}-U_{j,k}$, $(\Delta^{y}_{+}U)_{j,k} = U_{j,k+1}-U_{j,k}$,
and we can write:
\begin{multline*}
F(U_{i-1,j-1}^n,...,U_{i+1,j+1}^n) =F_R(U_{i,j-1}^n,U_{i-1,j-1}^n,
U_{i-1,j}^n,U_{i,j}^n,U_{i,j-1}^n) \\- 
F_L(U_{i-1,j+1}^n,U_{i-1,j}^n,U_{i,j-1}^n,U_{i,j}^n,U_{i,j+1}^n)\,,
\end{multline*}
with
\begin{align*}
F_R = h \, C_{xl} &\left( \overline{U}_{i+1,j}^{n+1} - 
\overline{U}_{i,j}^{n+1} \right)\\
&- C_{xl} C_{yr} 
\left(\overline{U}_{i-1,j+1}^{n+1} - 
\overline{U}_{i-1,j}^{n+1} - ( \overline{U}_{i,j+1}^{n+1} 
- \overline{U}_{i,j}^n ) \right) + 
\frac{1}{2}\left(f(U_{i-1,j}^n)+2f(U_{i,j}^n)+f(U_{i+1,j}^n)\right)\,,\\
F_L = h \, C_{xl} &\left( \overline{U}_{i-1,j}^{n+1} - 
\overline{U}_{i,j}^{n+1} \right)\\
&- C_{xl} C_{yl} 
\left(\overline{U}_{i-1,j-1}^{n+1} - 
\overline{U}_{i,j-1}^{n+1} - ( \overline{U}_{i-1,j}^{n+1} - 
\overline{U}_{i,j}^{n+1} )\right)  + 
\frac{1}{2}\left(f(U_{i-1,j}^n)+2f(U_{i,j}^n)+f(U_{i+1,j}^n)\right)\,,
\end{align*}
and
\begin{multline*}
G(U_{i-1,j-1}^n,...,U_{i+1,j+1}^n) =G_R(U_{i-1,j}^n,U_{i,j}^n,
U_{i+1,j+1}^n,U_{i,j+1}^n,U_{i+1,j}^n) \\- 
G_L(U_{i+1,j}^n,U_{i,j}^n,U_{i,j-1}^n,U_{i+1,j}^n,U_{i,j+1}^n)\,,
\end{multline*}
with
\begin{align*}
G_R = h \, C_{yl}  
&\left( \overline{U}_{i,j+1}^{n+1} - 
\overline{U}_{i,j}^{n+1} \right)\\
&-C_{yl}C_{xr} 
\left(\overline{U}_{i+1,j-1}^{n+1} - 
\overline{U}_{i+1,j}^{n+1} - ( 
\overline{U}_{i,j-1}^{n+1} - 
\overline{U}_{i,j}^{n+1} ) \right) + 
\frac{1}{2}\left(g(U_{i,j-1}^n)+2g(U_{i,j}^n)+g(U_{i,j+1}^n)\right)\,,\\
G_L = h \, C_{yl}  
&\left( \overline{U}_{i,j-1}^{n+1} - 
\overline{U}_{i,j}^{n+1} \right)\\
&-C_{yr}C_{xl}  
\left(\overline{U}_{i+1,j+1}^{n+1} - 
\overline{U}_{i+1,j}^{n+1} - ( 
\overline{U}_{i,j+1}^{n+1} - 
\overline{U}_{i,j}^{n+1} )\right)  +  
\frac{1}{2}(g(U_{i,j-1}^n)+2g(U_{i,j}^n)+g(U_{i,j+1}^n))\,.
\end{align*}

In order for the scheme (\ref{proyF1}) be consistent for the 2D scalar case, 
we must have:
\begin{equation}
 F(u,u,...,u) = \mathcal{F}_1(u) \qquad \text{ and } \qquad
 G(u,u,...,u) = \mathcal{F}_2(u), \quad u \in \REAL.
 \label{sclpcons2d}
\end{equation}
We can notice that $F$ and $G$ satisfy condition 
(\ref{sclpcons2d}). This implies consistency with 
$\frac{\partial U}{\partial t} + \Div \FluxSWE(U) = 0$.
Indeed, it turns out that conservative 
monotone schemes as in (\ref{proyF1}) converge to 
entropy solutions \cite{chap:barth17,art:crandall80,book:Leveque2002}. 
It is possible to write the Lagrangian-Eulerian scheme applied to
scalar equations in the form of standard monotone finite difference
schemes (see, e.g.,
\cite{art:abreu22b,art:abreu22,art:abreu20,art:abreu21,art:abreu20b,art:abreu19,art:abreu18a,art:abreu17,art:abreu18b,phd:perez15}),
showing that no Riemann problems are needed to obtain a correct
entropy solution. Generalization to systems can be done only for
special cases, and its identification for the ISWE system is subject
of current studies.

\section{Numerical simulations}
\label{sec:results}
The numerical experiments presented in this section are devised in
order to show that the geometrically intrinsic Lagrangian-Eulerian
scheme is able to handle non-autonomous fluxes coming from the spatial
variability of the bottom topography and preserve the well-balance
property of the system in the discrete setting.
We define verify experimentally the applicability of the scheme on
different surfaces by qualitatively comparing the numerical results
against the benchmarks published in~\cite{art:BP20, phd:Bachini20}.
In addition, we use for all the test cases long-time simulations to
test the well-balance property of the proposed scheme by running the
code until $t=50$~s.

All the test cases consider a gravity-driven flow in a dam-break setting, with
$T_{sw}=0$ and $\tau_b=0$, and no-flow boundary conditions.
Initial condition is set to
describe the dam-break, with zero velocity and depth defined
differently on each surface, according to the specific geometry.
Initial  conditions will be specified for each test case.
The CFL condition~\eqref{matrixcfl2} is used to define an appropriate
time-step size at every time step.
Verification of the well-balance property is performed by looking at
the time-behavior of the $L_1$ discrete norm:
\begin{equation*}
||\Source||_1 =  A(\Region[n]_{i,j}) \sum_{i,j} |\Source_{i,j}|\,.
\end{equation*}

We consider four different bottom surfaces with an increasing
complexity in the geometry: i) a simple sloping plane (i.e., constant
metric), ii) a parabola (i.e., one-dimensional curvature effects),
iii) a hyperboloid with a central bump (i.e., two-dimensional
dependence of the geometry with a central radial symmetry), and iv) a
fully three-dimensional surface.
In all test cases we define our surface by the use of a global
parametrization, where the height function is given in input as
$\zcg=\BSM(x,y)$. The working grid is obtained from a regular
Cartesian grid of a subset $\SubsetU\subset\REAL^2$ expressing the
chart by moving the vertices of the squares vertically to the surface
$\SurfDomain$.

\paragraph{Sloping plane} 
The chart domain is the rectangle
$\SubsetU=[0,1]\times[0,10]\subset\REAL^2$ and the height function is
given by $\BSM(x,y)=0.1y+1.0$. The initial value for the water depth is set
to 2~m for cells located upstream of the dam position, located at
$y=2.0$~m, and a depth of 1~m downstream. A CFL value of $0.20$ is considered to run the
simulation.

\begin{figure}
	\centering
	\includegraphics[width=0.38\textwidth]{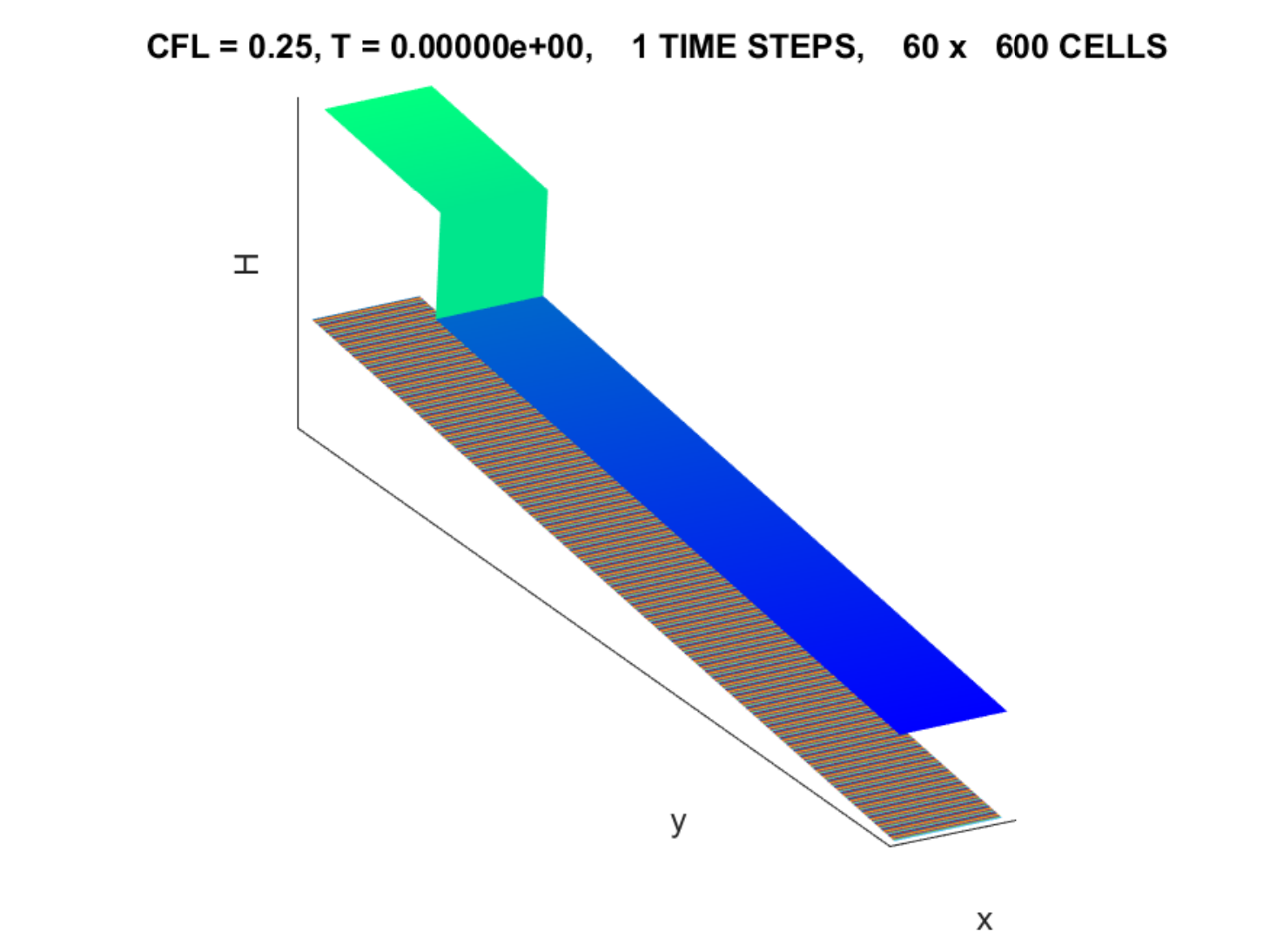}
	\includegraphics[width=0.38\textwidth]{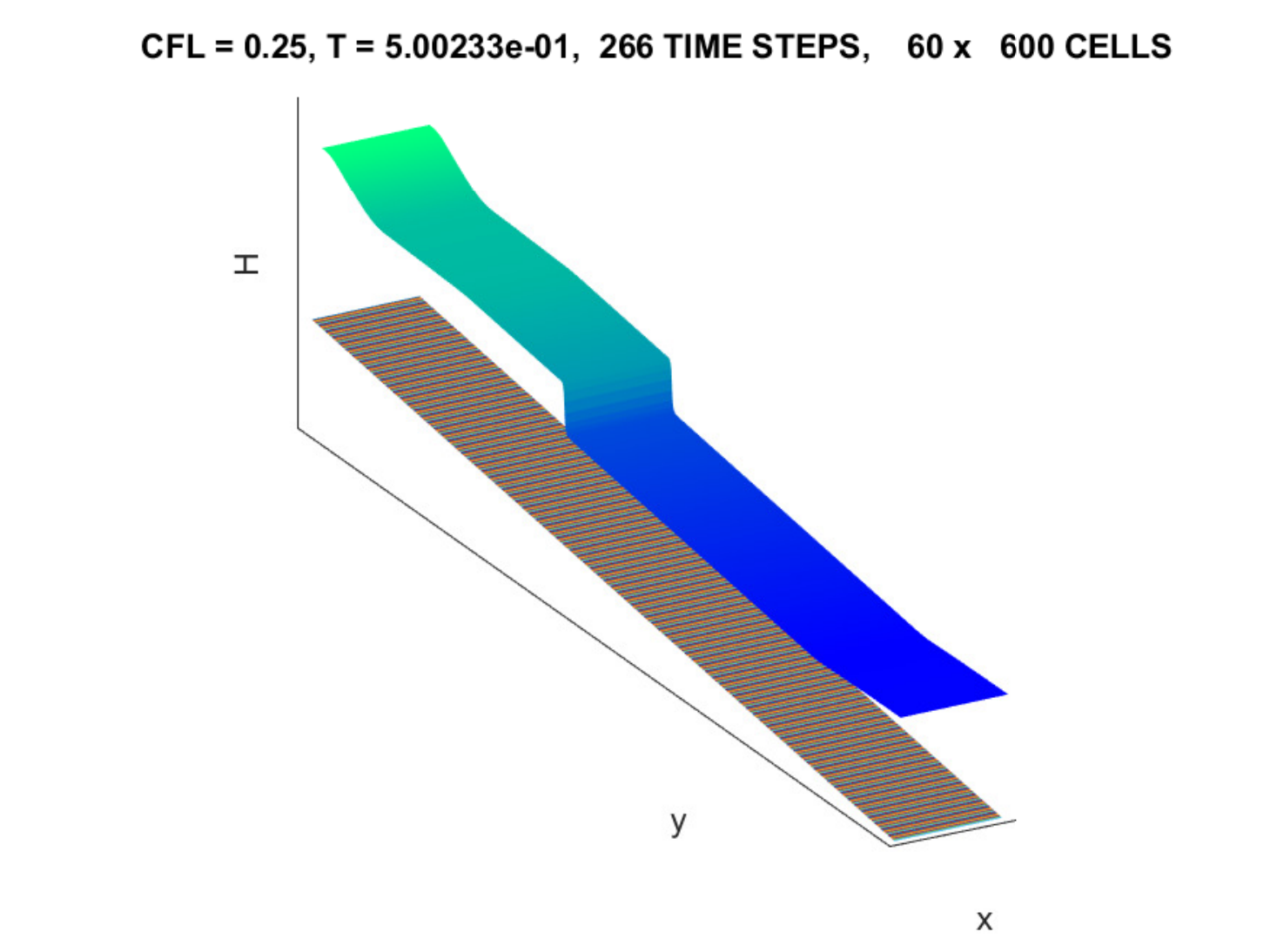}
	\includegraphics[width=0.38\textwidth]{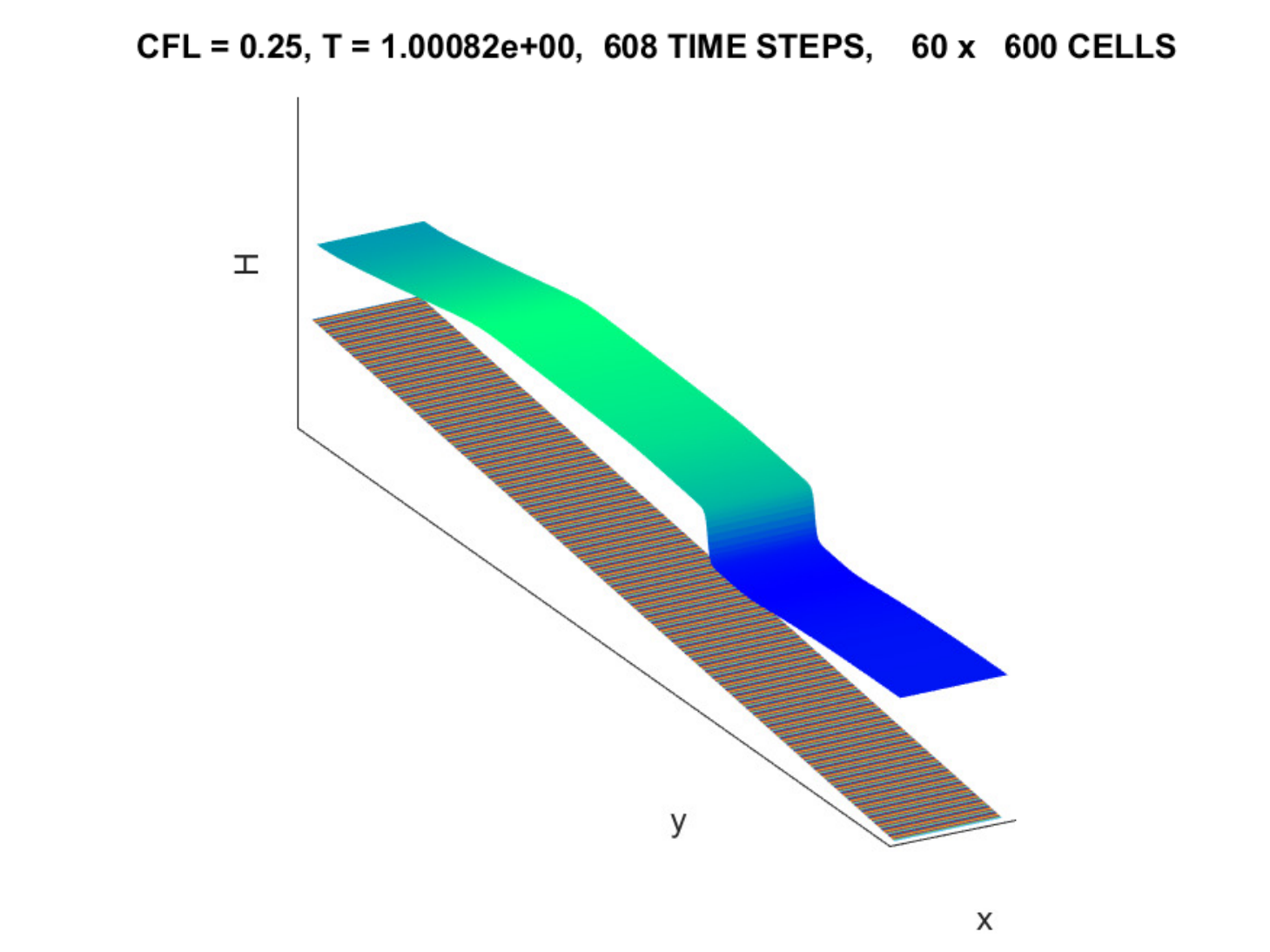}
	\includegraphics[width=0.38\textwidth]{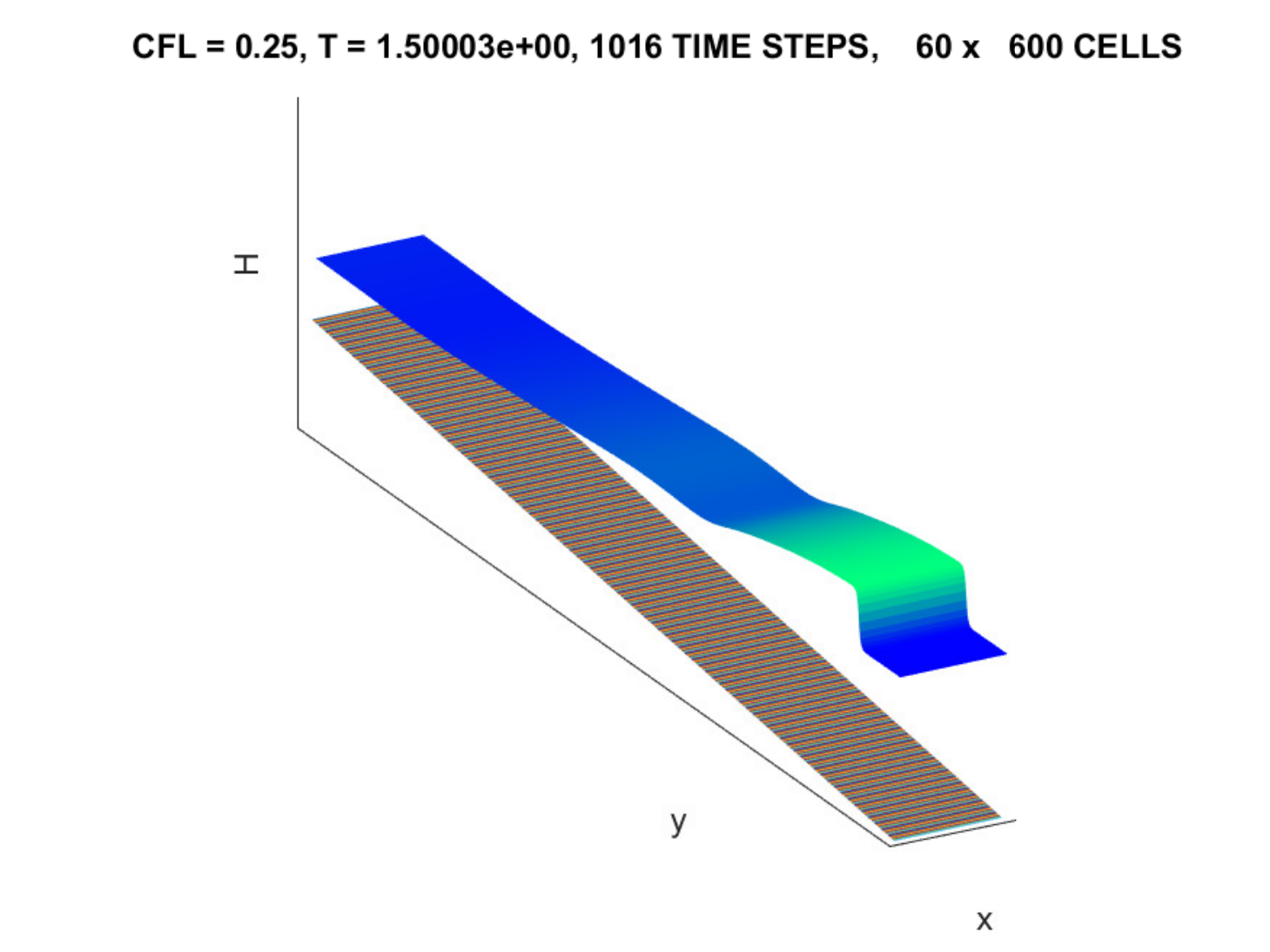}
	\caption{Sloping plane: evolution of the water depth at times $t=$0, 0.5, 1, 1.5~s.}
	\label{fig:Slopeplane}
\end{figure}

\begin{figure}
	\centering
	\includegraphics[width=0.38\textwidth]{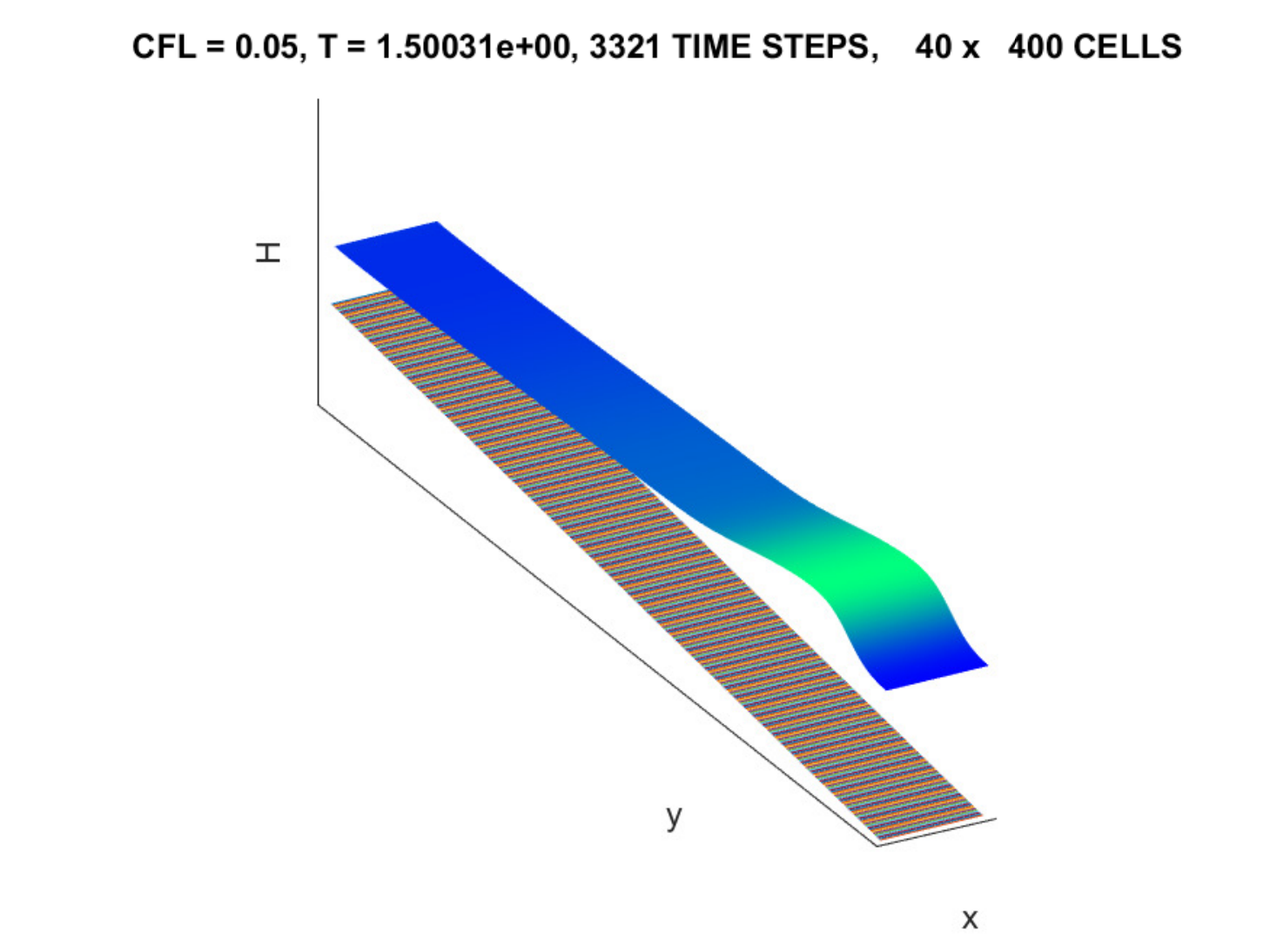}
	\includegraphics[width=0.38\textwidth]{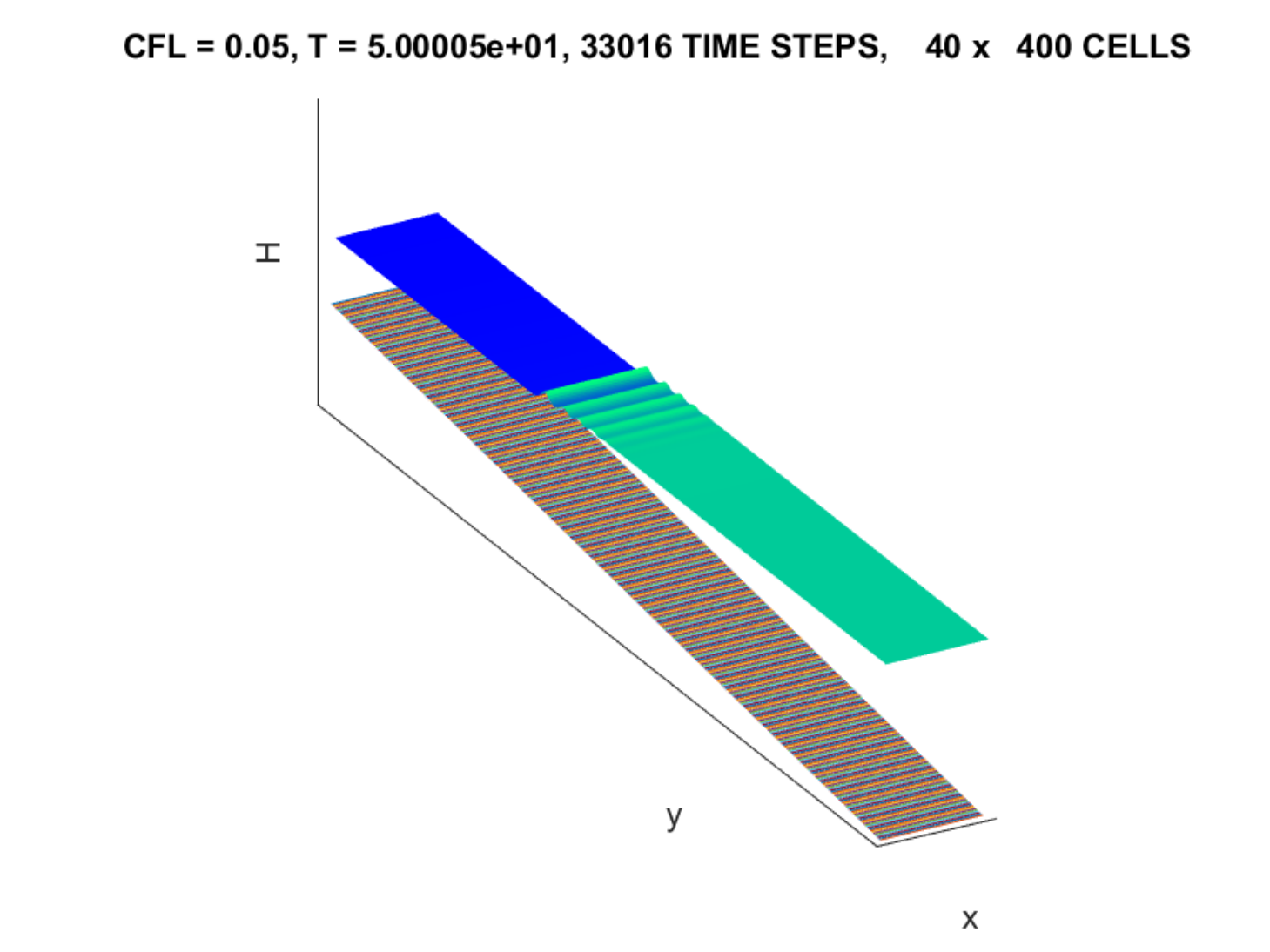}
	\includegraphics[width=0.38\textwidth]{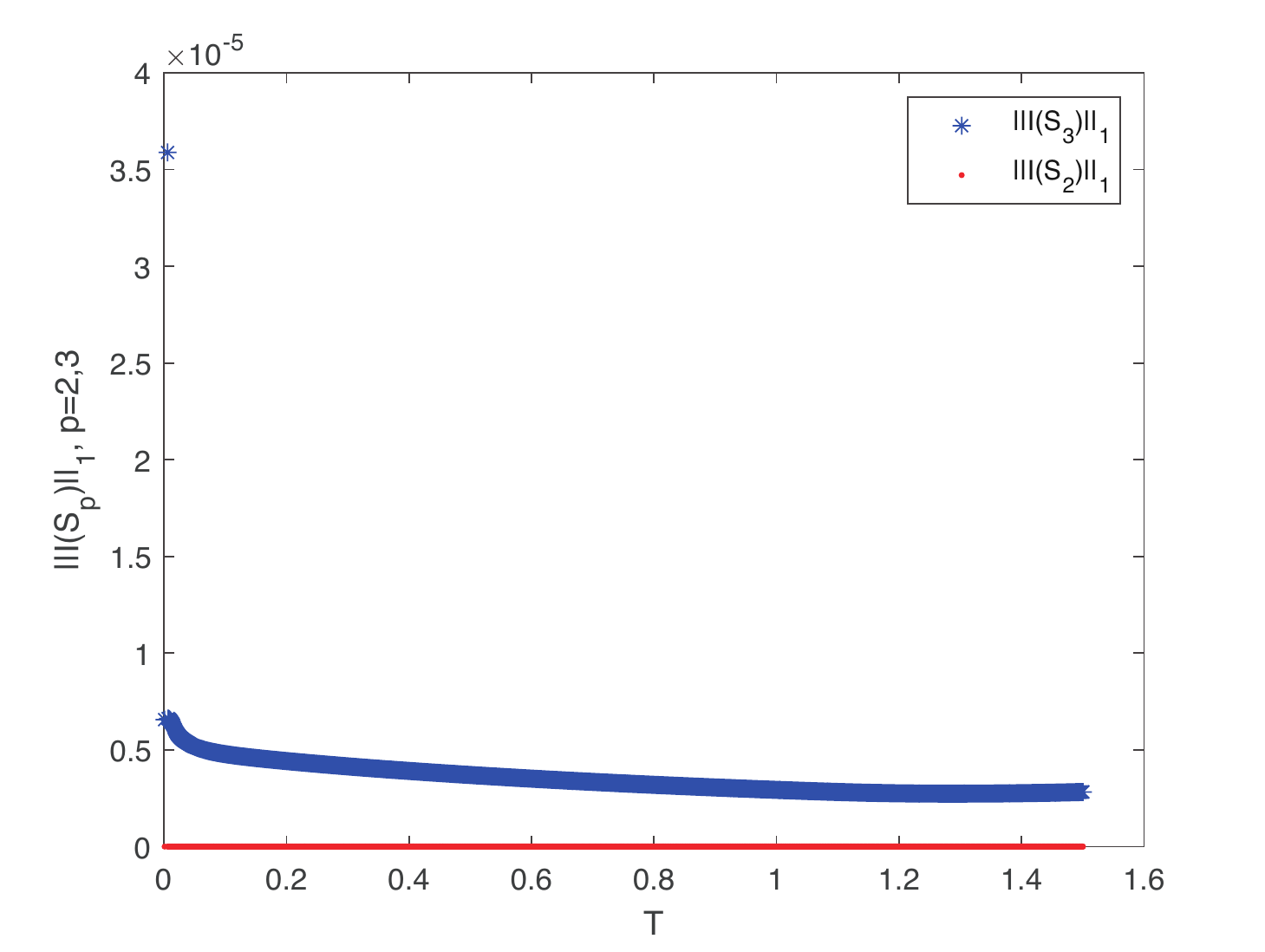}
	\includegraphics[width=0.38\textwidth]{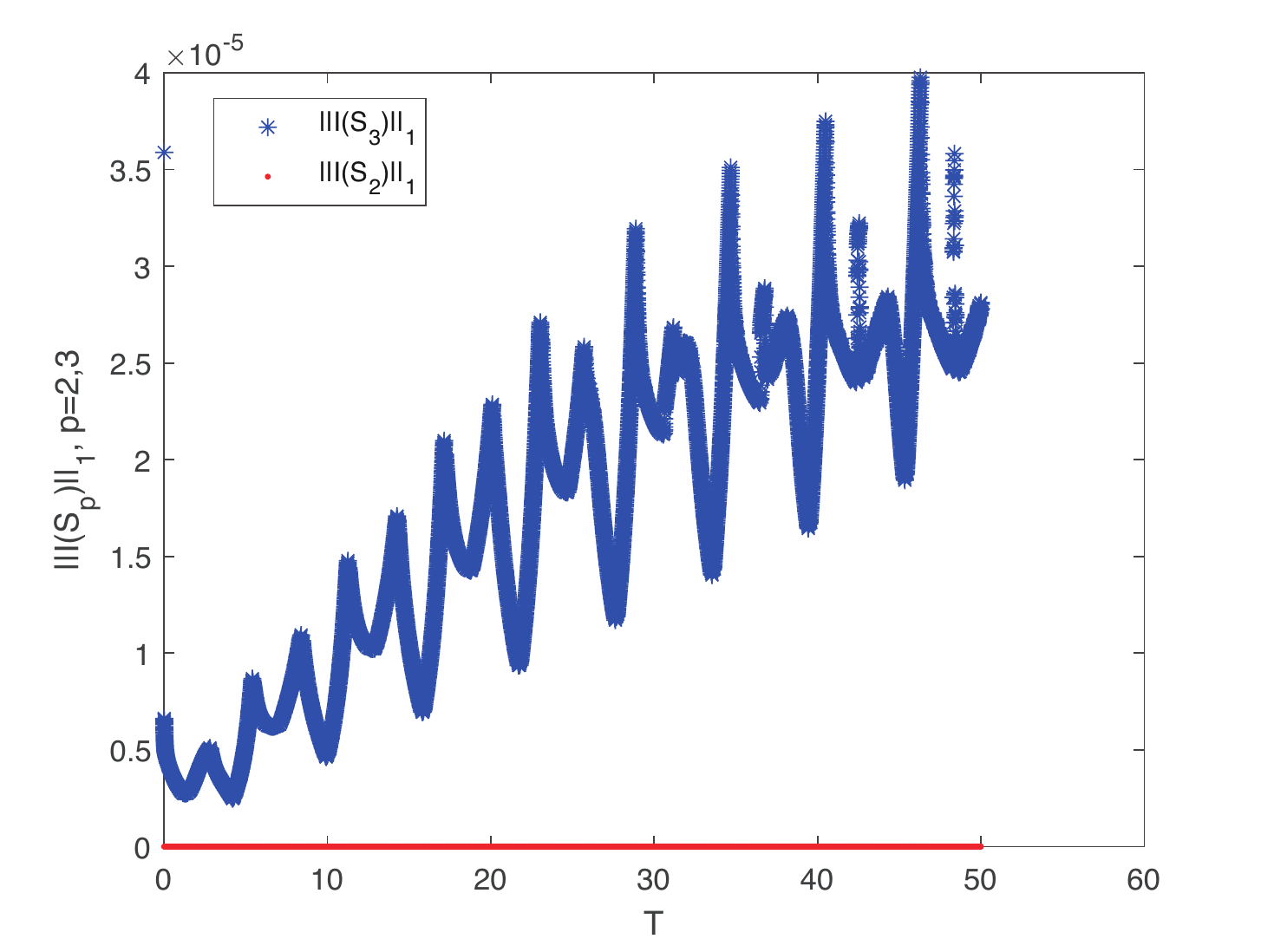}
	\caption{Sloping plane, 40$\times$400 cells: water depth (top panels) and
          plots of the $L^1$-norm of $I(\Source^2)$ and $I(\Source^3)$ (bottom panels), at times
          $t=1.5$~s (left) and $t=50$~s (right). Computational time for $t=50$: 490.75~s.}
	\label{CI1}
\end{figure}

\begin{figure}
	\centering
	\includegraphics[width=0.38\textwidth]{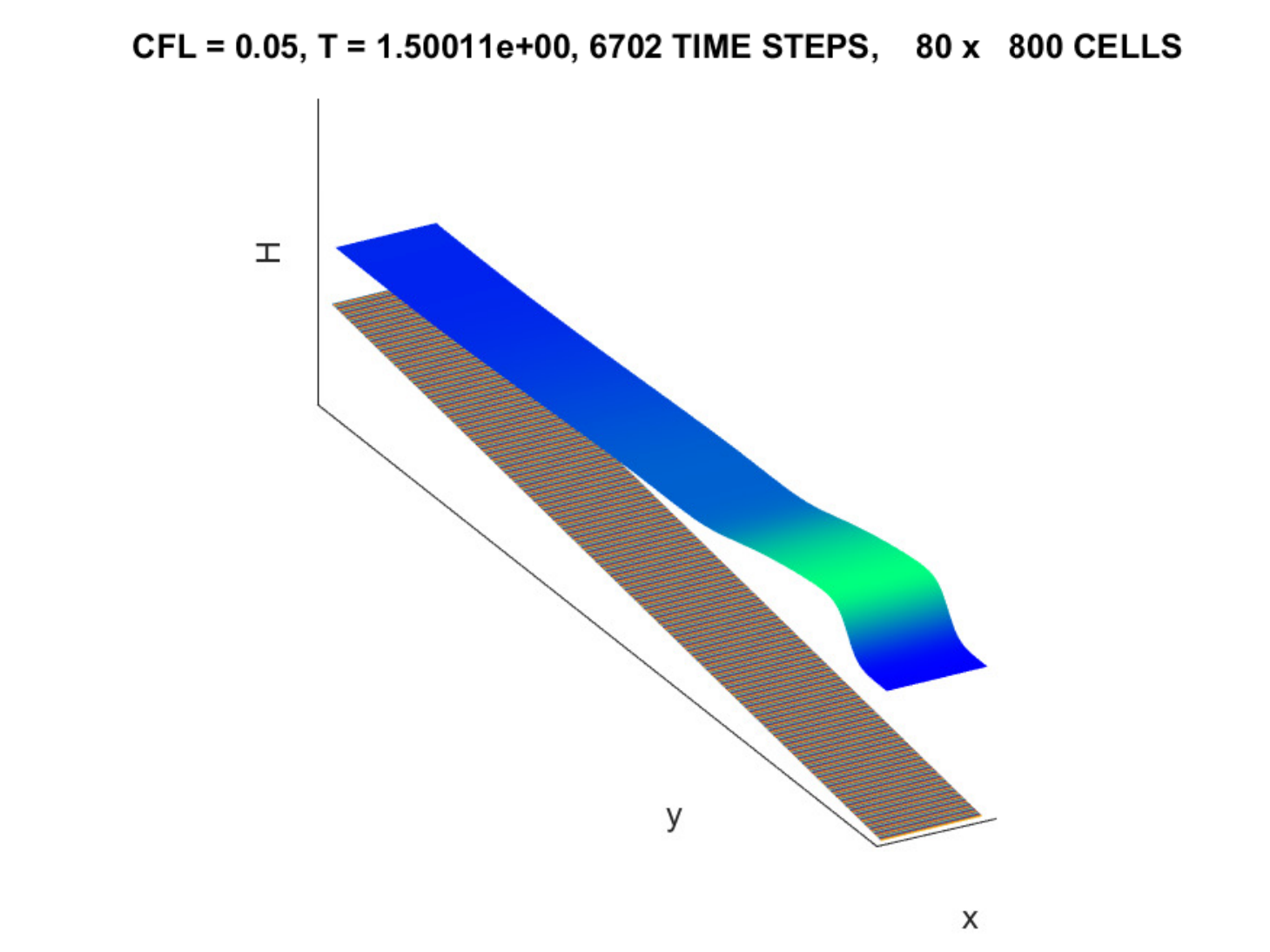}
	\includegraphics[width=0.38\textwidth]{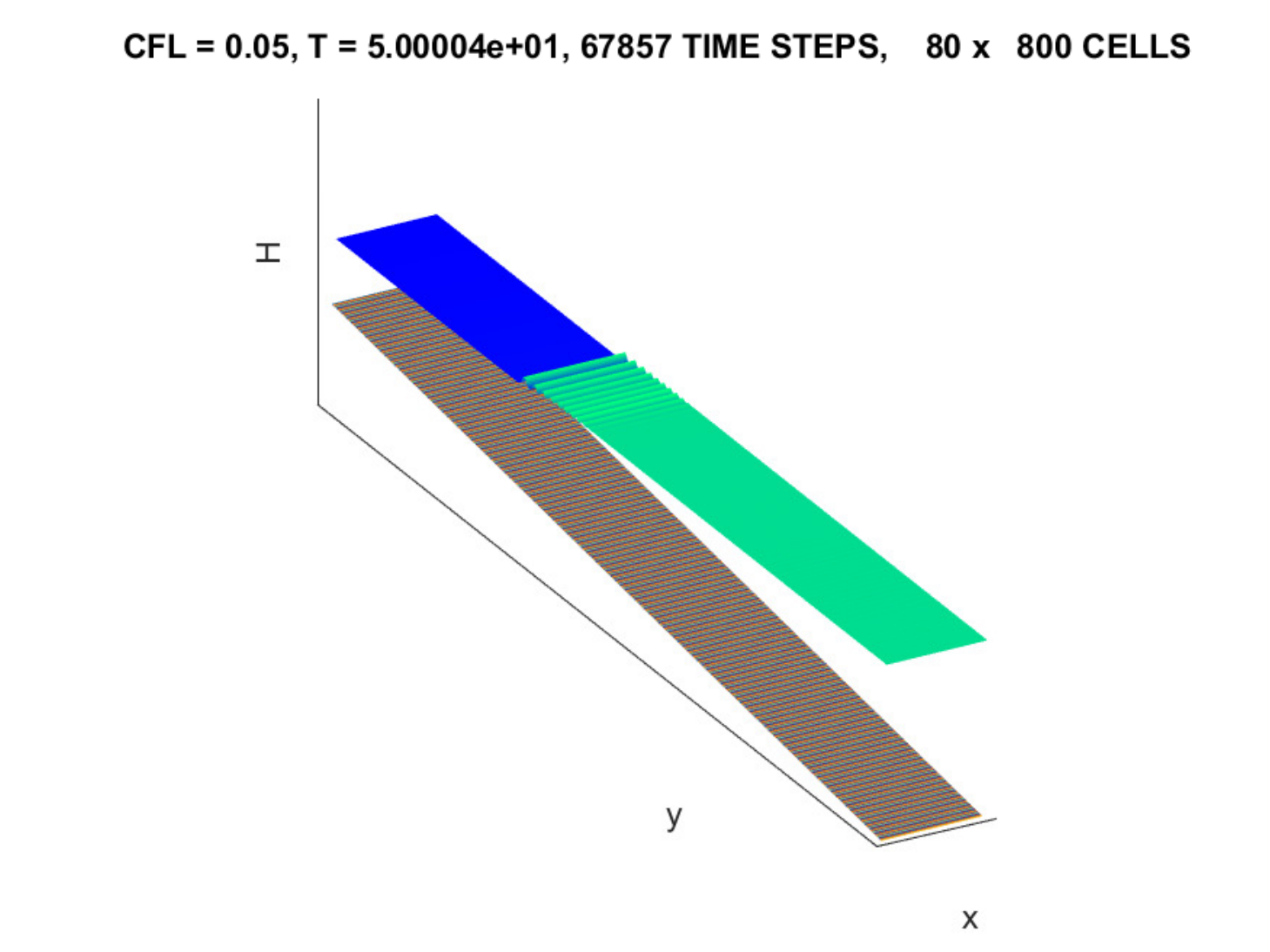}
	\includegraphics[width=0.38\textwidth]{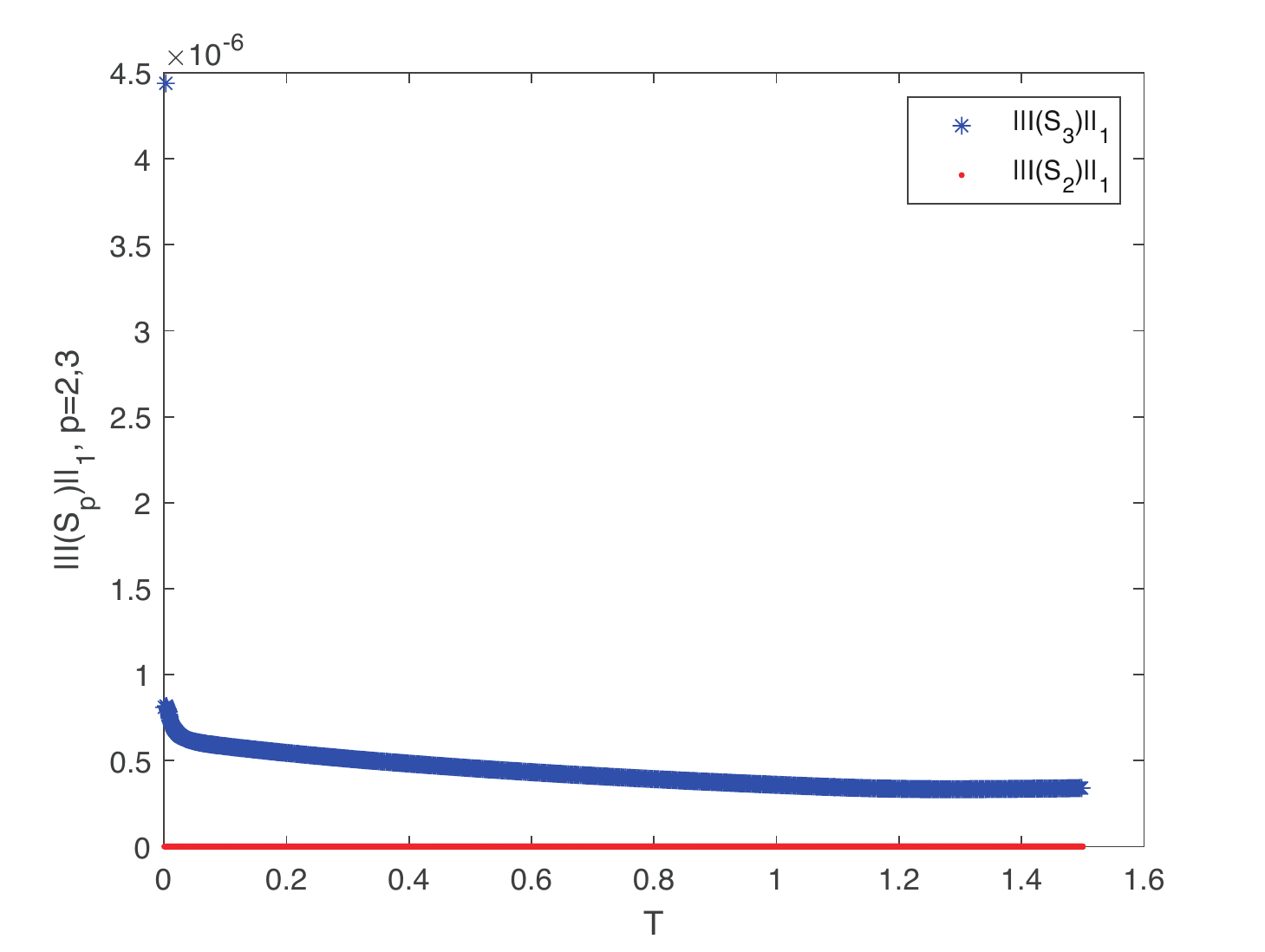}
	\includegraphics[width=0.38\textwidth]{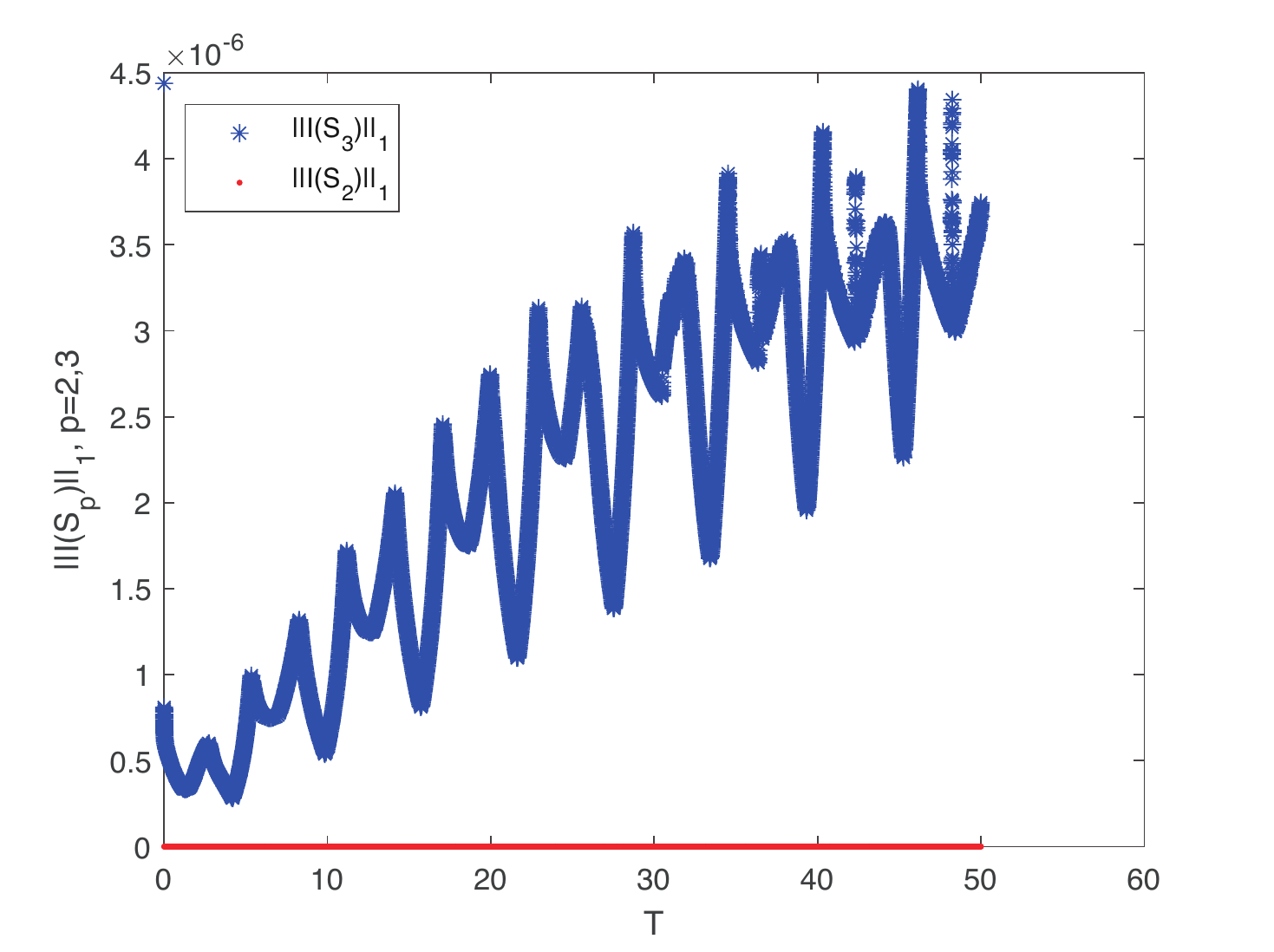}
	\caption{Sloping plane, 80$\times$800 cells: water depth (top panels) and
          plots of the $L^1$-norm of $I(\Source^2)$ and $I(\Source^3)$ (bottom panels) at time
          $t=1.5$~s (left) and $t=50$~s (right). Computational time for $t=50$: 2726~s.}
	\label{CI2}
\end{figure}
Figure~\ref{fig:Slopeplane} shows the evolution of the normal depth at times $t=$0, 0.5, 1, 1.5~s on a grid
formed by $60\times600$ cells along the $x$ and $y$-directions,
respectively.
Figures~\ref{CI1} and \ref{CI2} show the water depth at $t=1.5$~s and
$t=50$~s (top panels) over a mesh with $40\times 400$ cells and
$80\times 800$ cells, respectively. The $\Lspace[1]$-norm of $\Source$
is plotted in the bottom panels.
The results are very similar to those of~\cite{art:BP20, phd:Bachini20}.

\paragraph{Parabola}
In this test case we consider a simple one-dimensional flow in which
the effect of the curvature is present. The parabola is defined by the
height function $\BSM(x,y)=0.04(y-10)^2$ on the rectangular domain
$\SubsetU=[0,1]\times[0,10]\subset\REAL^2$. The initial
depth is 2~m of water upstream and 1~m downstream of the dam located
at $y=2.0$~m and
a CFL value of $0.30$
is considered.

\begin{figure}
	\centering
	\includegraphics[width=0.38\textwidth]{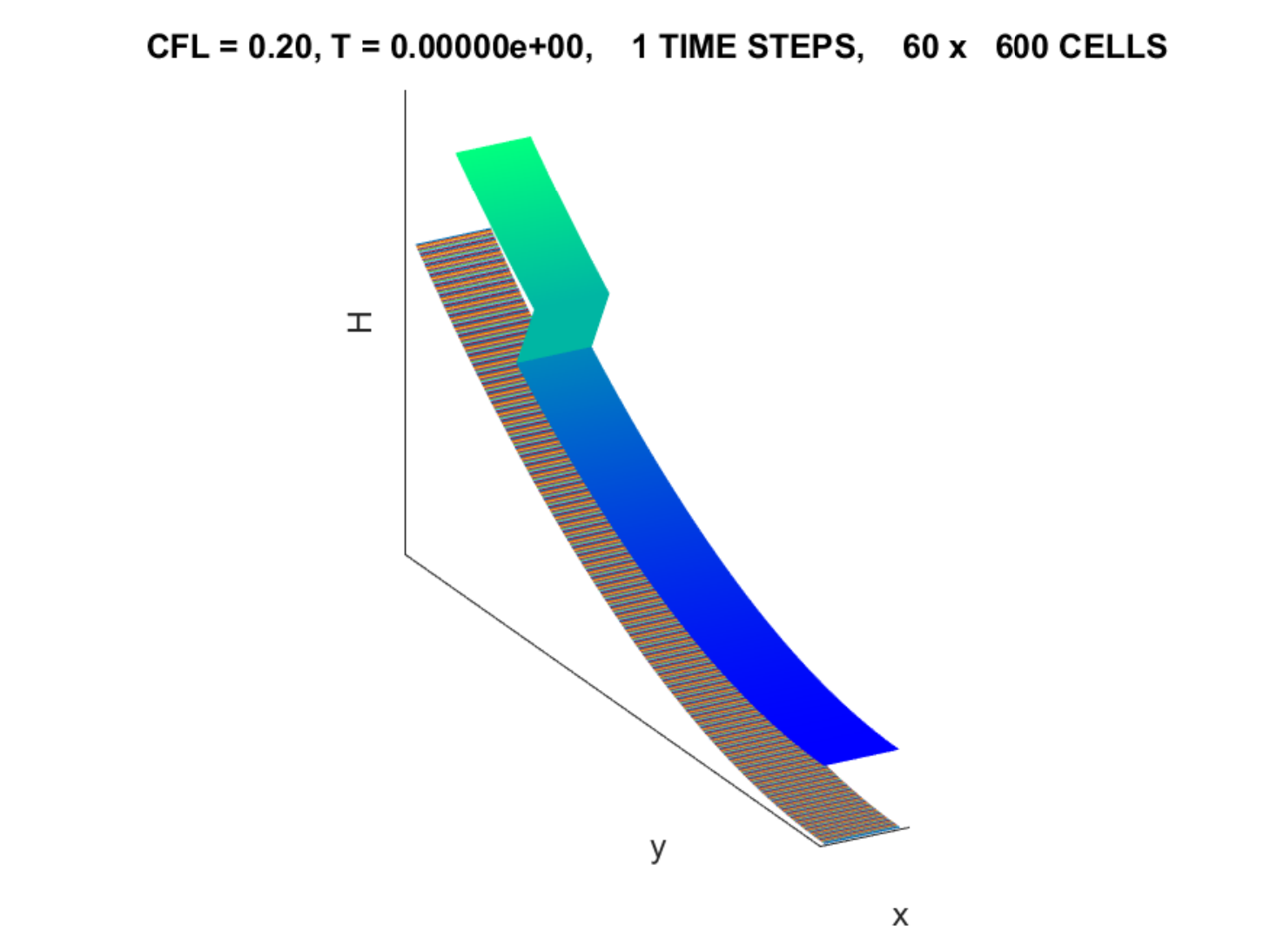}
	\includegraphics[width=0.38\textwidth]{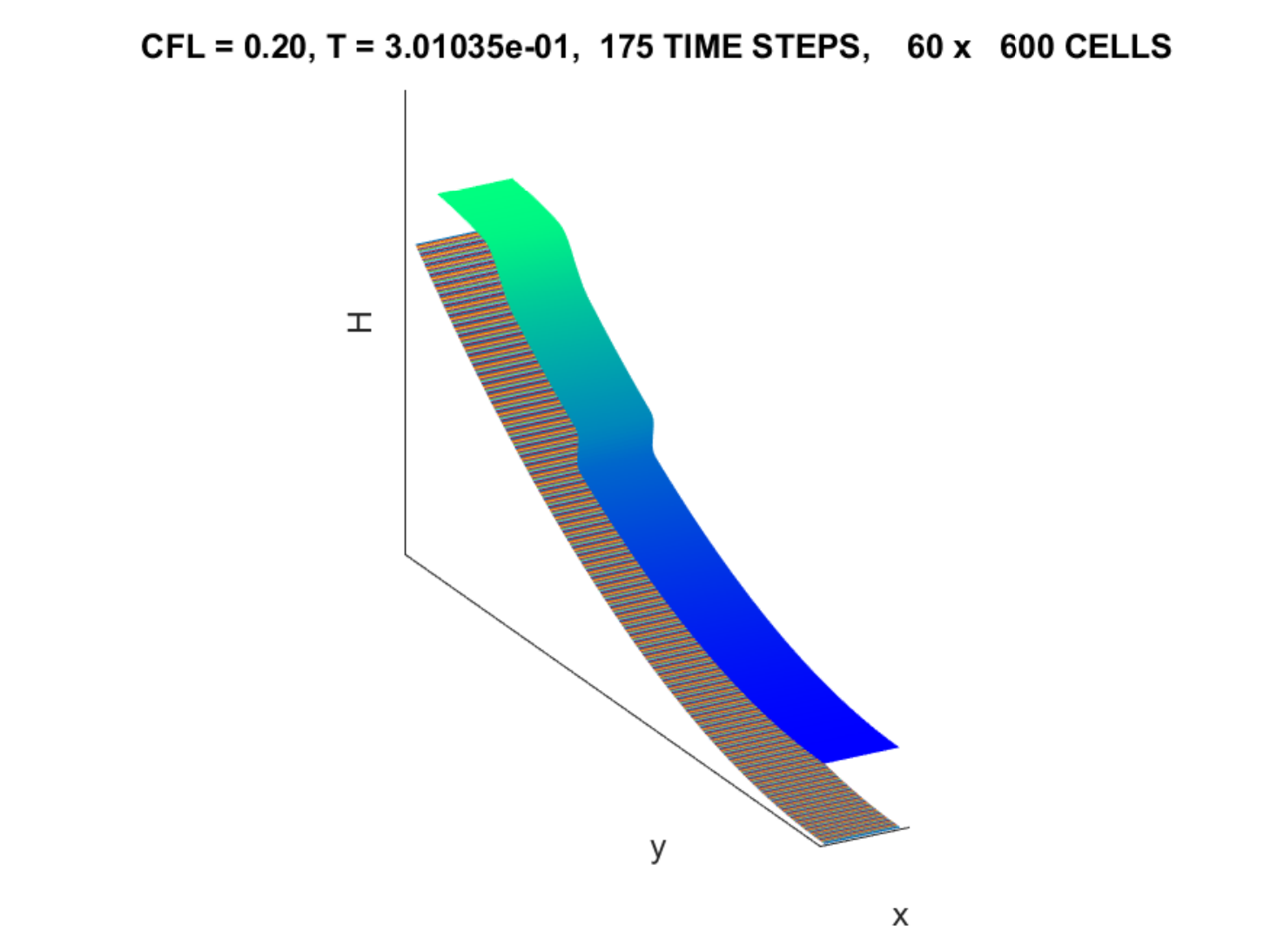}
	\includegraphics[width=0.38\textwidth]{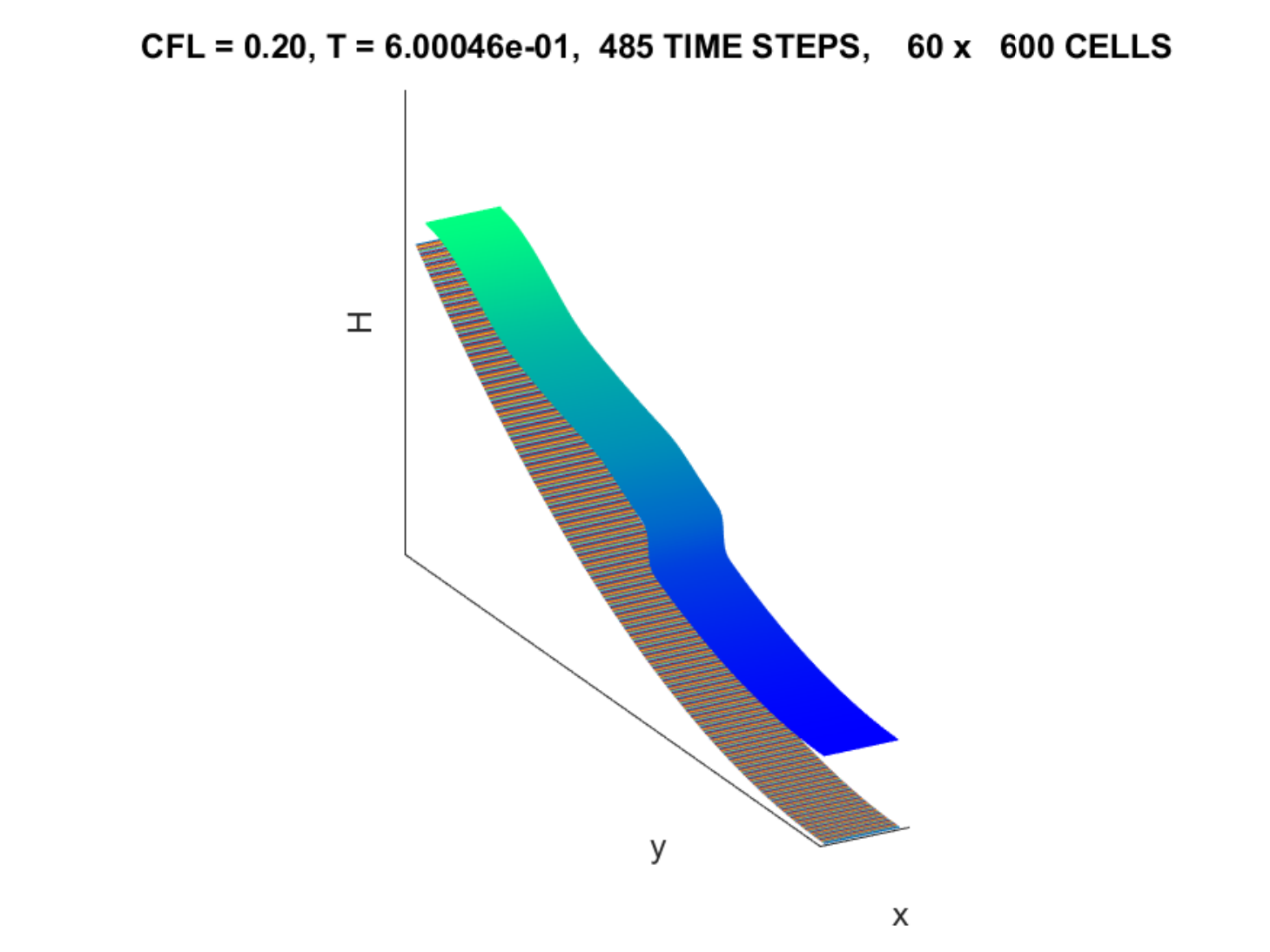}
	\includegraphics[width=0.38\textwidth]{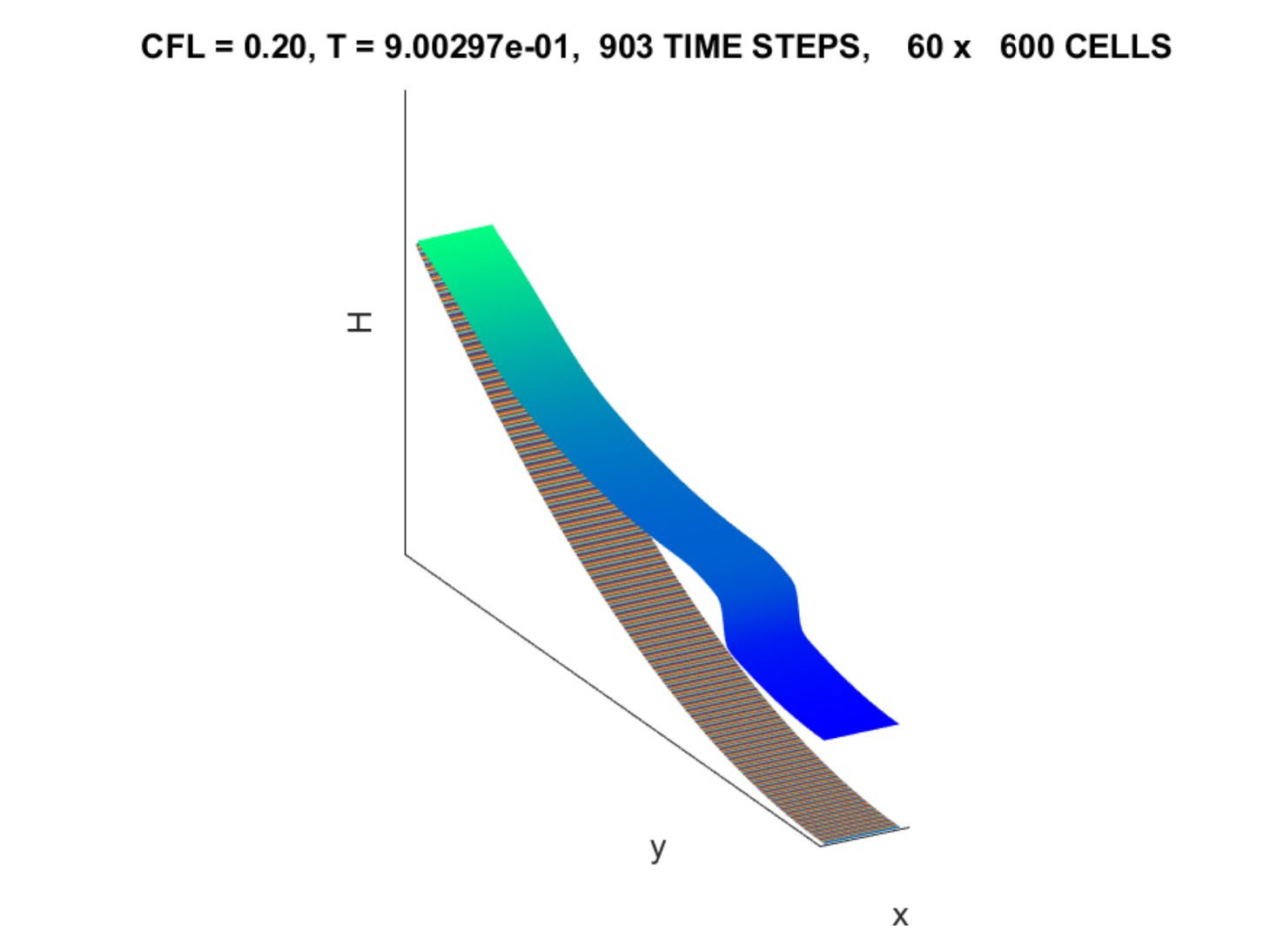}
	\caption{Parabola case: evolution of the water depth at times $t=$0.0, 0.3, 0.6, 0.9~s.}
	\label{parabolacase}
\end{figure}

\begin{figure}
	\centering
	\includegraphics[width=0.38\textwidth]{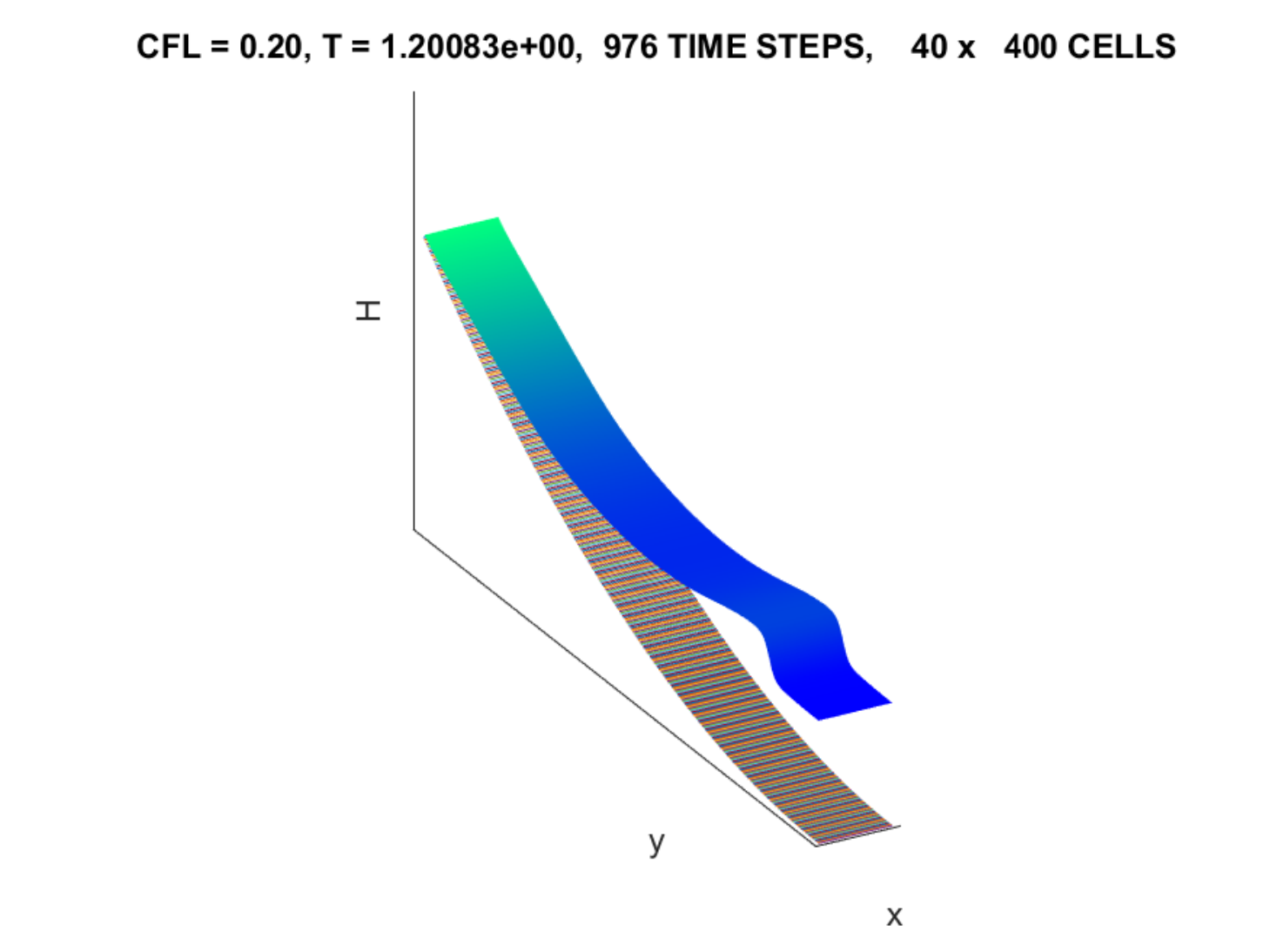}
	\includegraphics[width=0.38\textwidth]{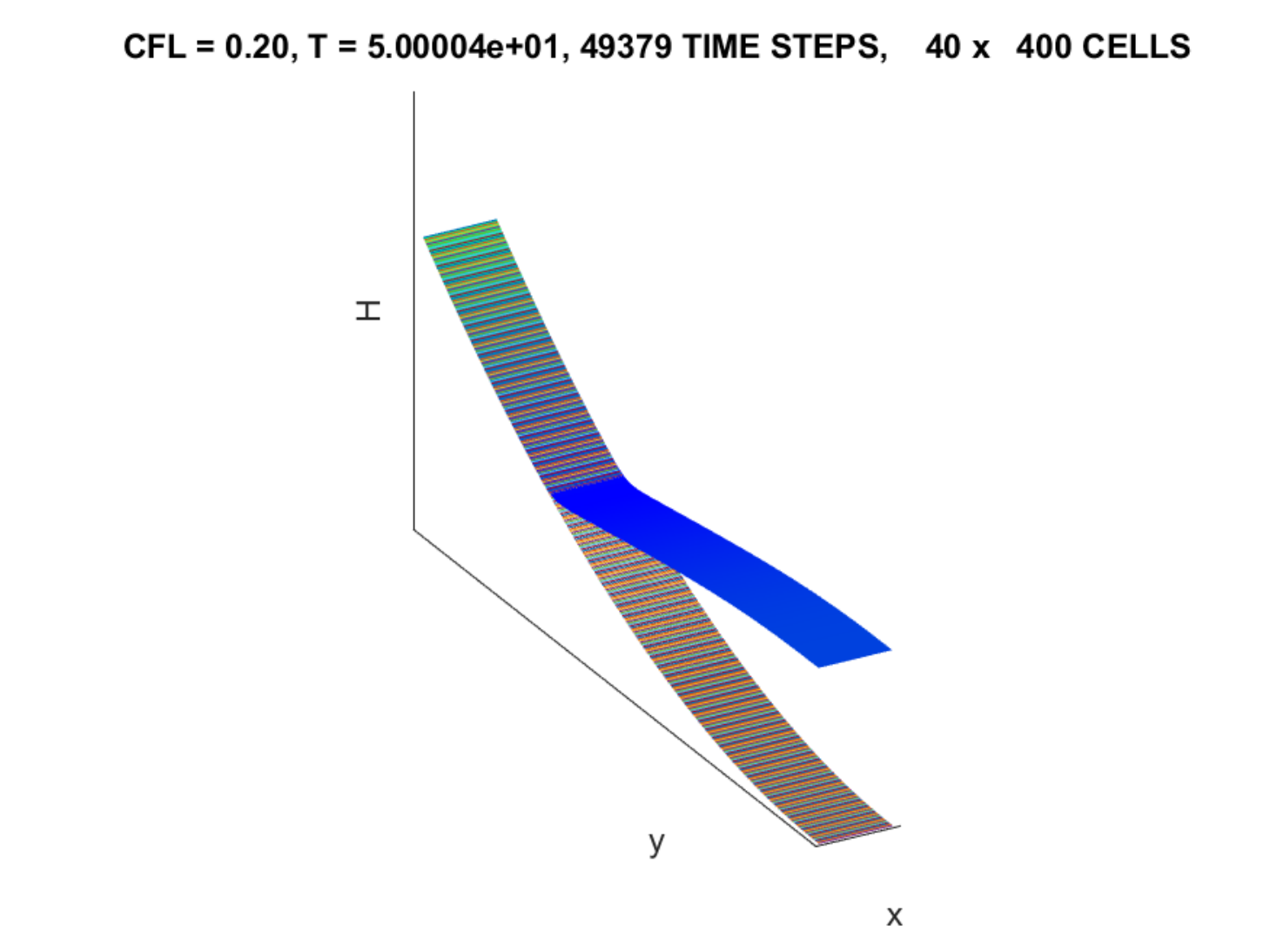}
	\includegraphics[width=0.38\textwidth]{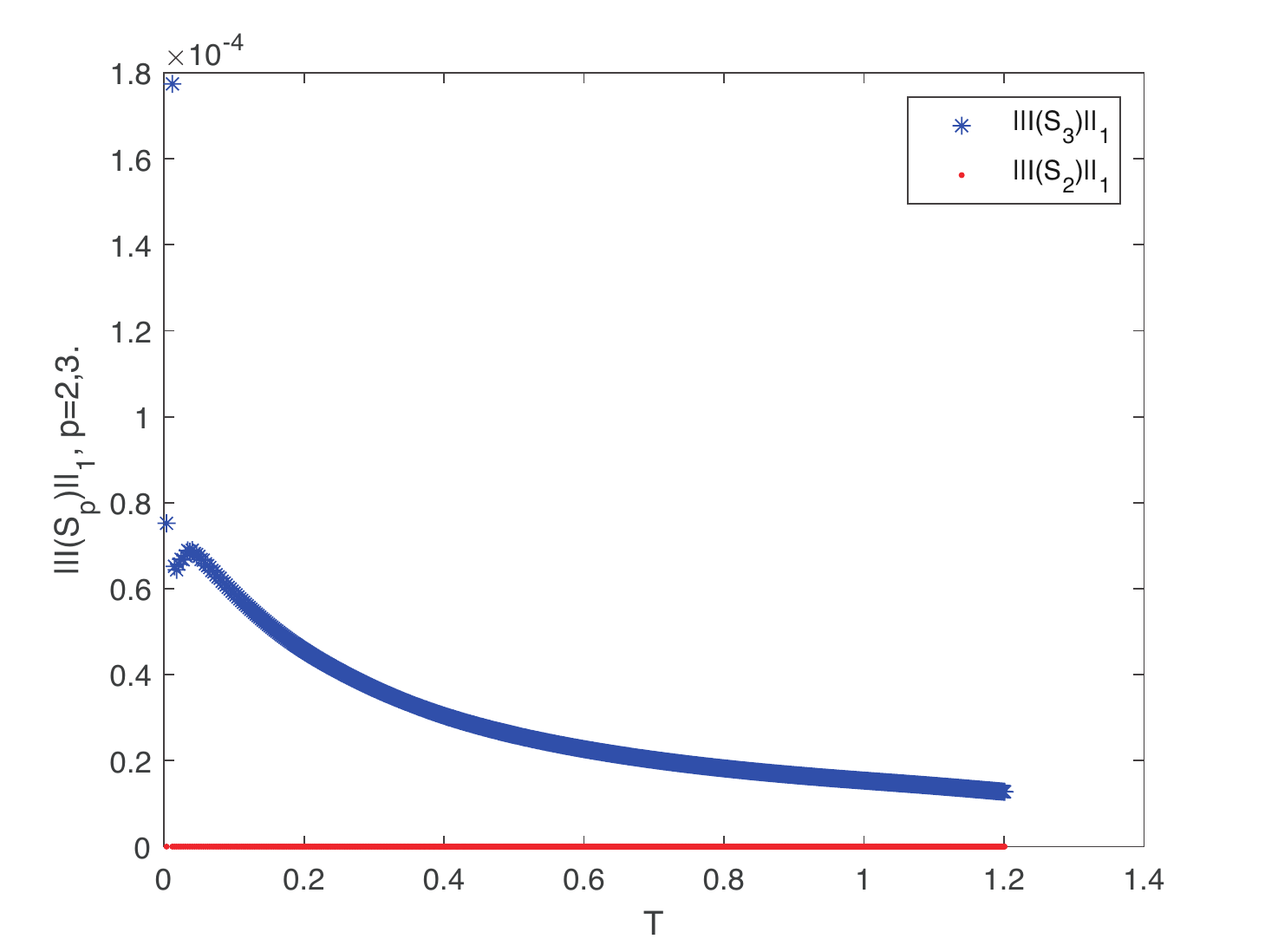}
	\includegraphics[width=0.38\textwidth]{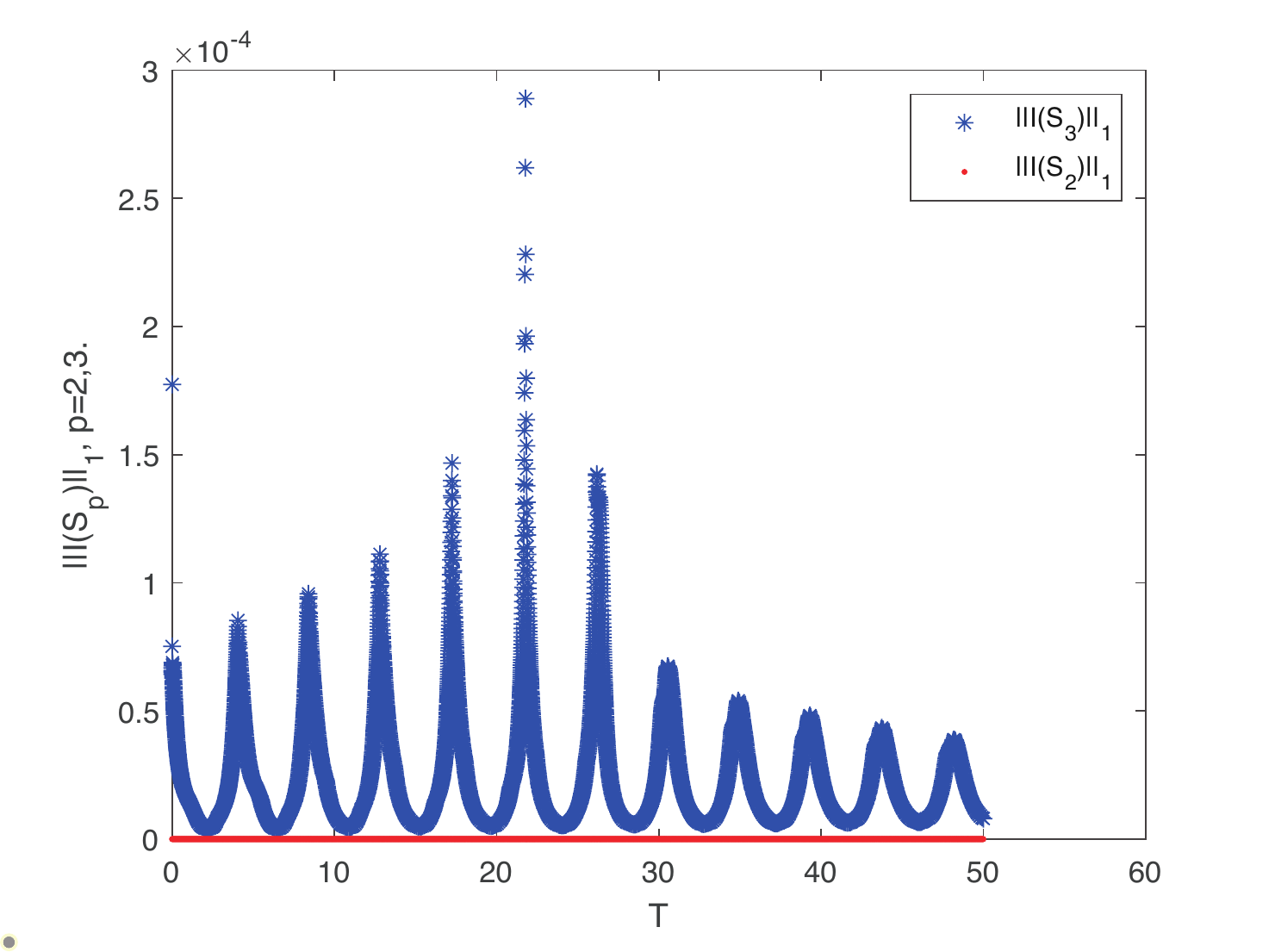}
	\caption{Parabola case, 40$\times$400 cells: water depth (top panels) and
          plots of the $L^1$-norm of $I(\Source^2)$ and $I(\Source^3)$ (bottom panels), at times
          $t=1.2$~s (left) and $t=50$~s (right). Computational time for $t=50$: 676.11~s.} 
	\label{CI3}
\end{figure}

Figure~\ref{parabolacase} shows the normal water depth evolution at
times $t=$0.0, 0.30, 0.6, 0.9~s on a grid formed by 60$\times$600
cells. Figure~\ref{CI3} presents the numerical results to verify the
discrete well-balanced property on a 40$\times$400 grid.
The long time well-balance error
tends to decrease maintaining a maximum value below $10^{-4}$ at all
times. The same simulation was run on a 80$\times$800 grid with
results consistent with the finer resolution.

\paragraph{Hyperboloid-central-bump}
The third example considers a surface defined by the height function:
\begin{equation*}
  \BSM(x,y)=-\frac{4}{5}\sqrt{(x)^2+(y)^2+1}
\end{equation*}
on a two-dimensional square domain
$\SubsetU=[-3,3]\times[-3,3]\subset\REAL^2$. In this case the bottom
surface presents curvatures in both directions.  The initial condition
is defined in order to capture the radial symmetric feature of the
bottom surface in the behavior of the solution: 2~m of water are
considered in the central part of the domain, i.e. inside of a circle
of radius $0.5$~m, and 1~m outside of this central radius.
A CFL
number of $0.30$ is chosen in the condition~\eqref{matrixcfl2}.

\begin{figure}
	\centering
	\includegraphics[width=0.38\textwidth]{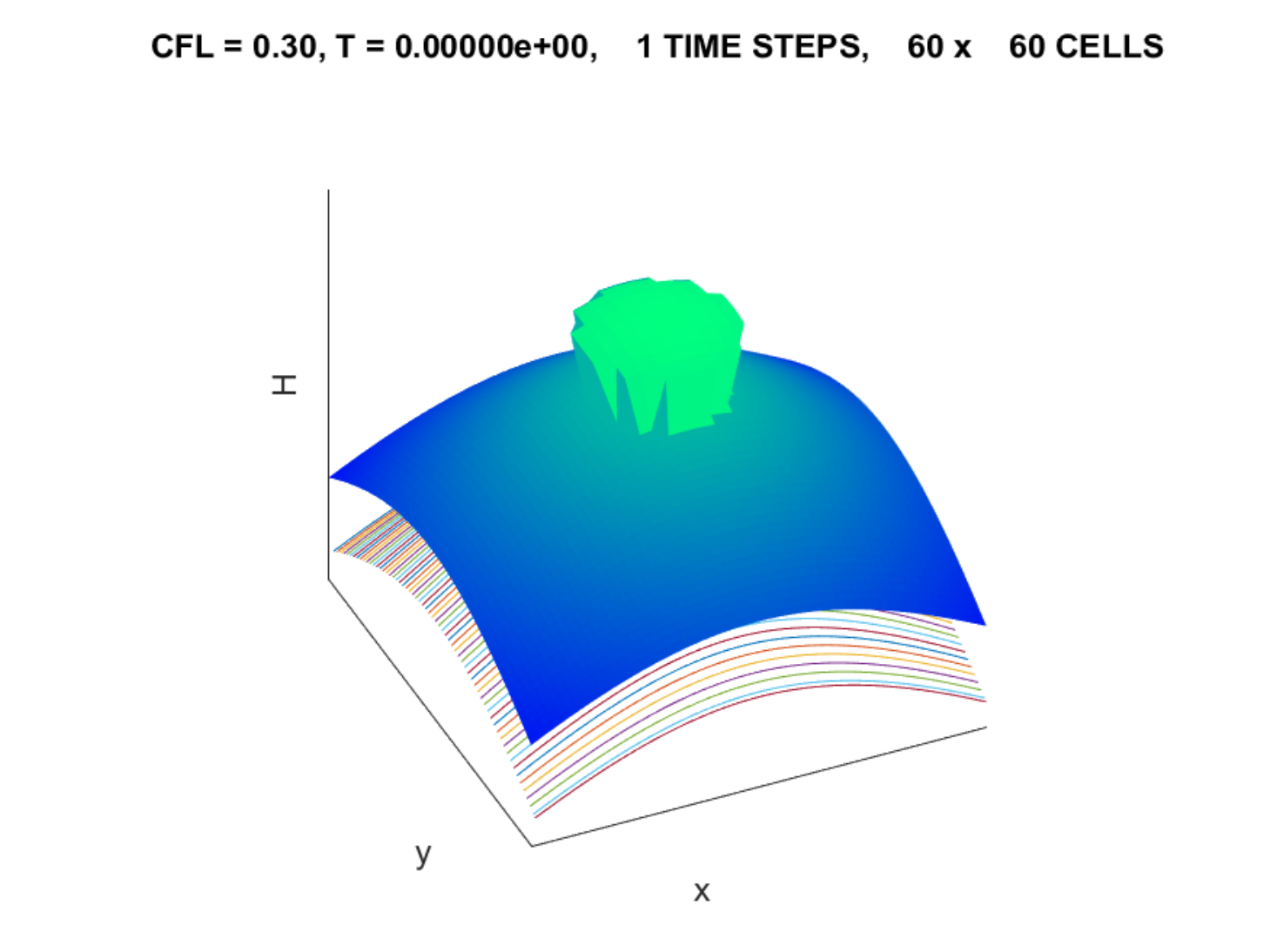}
	\includegraphics[width=0.38\textwidth]{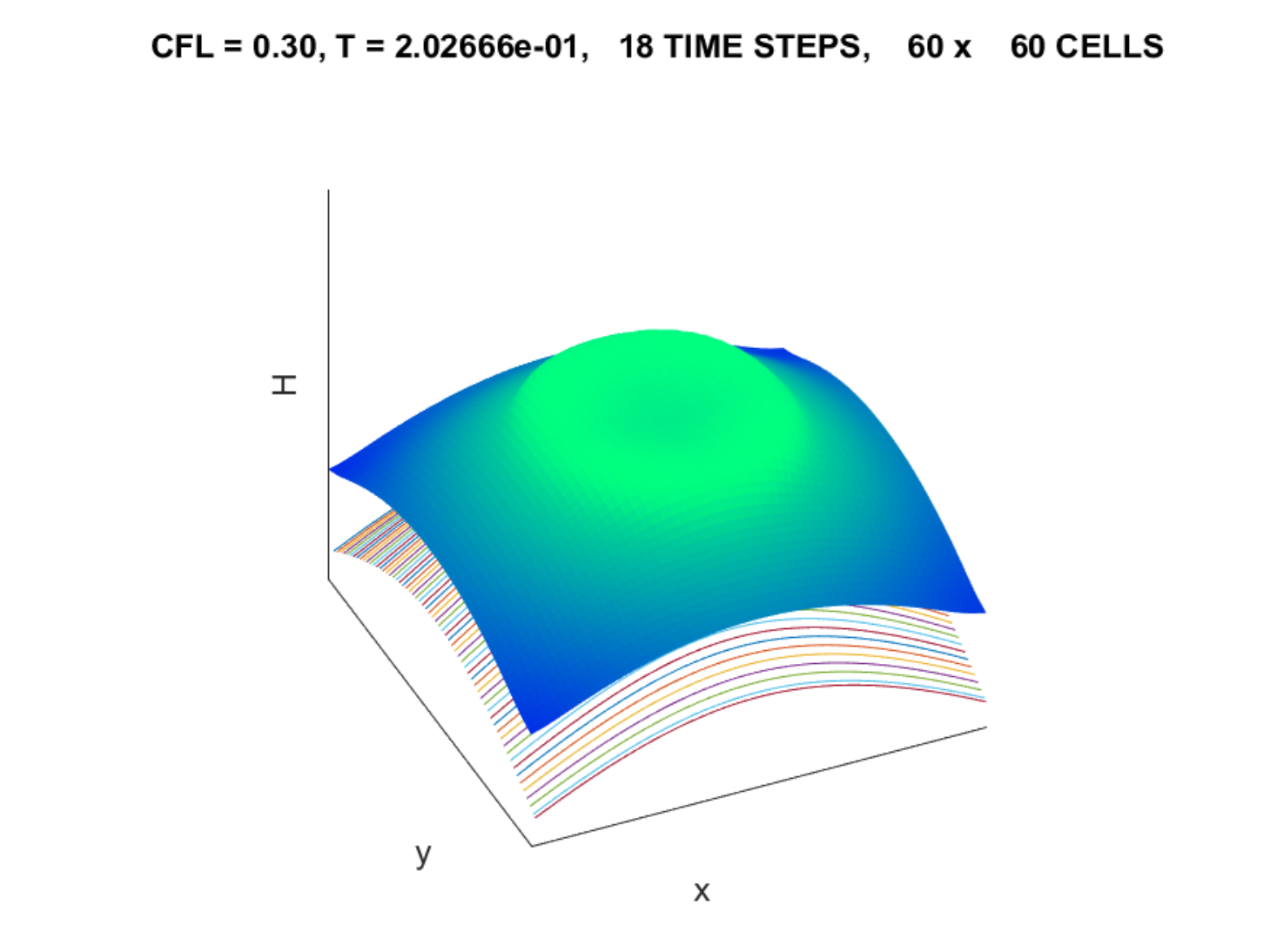}
	\includegraphics[width=0.38\textwidth]{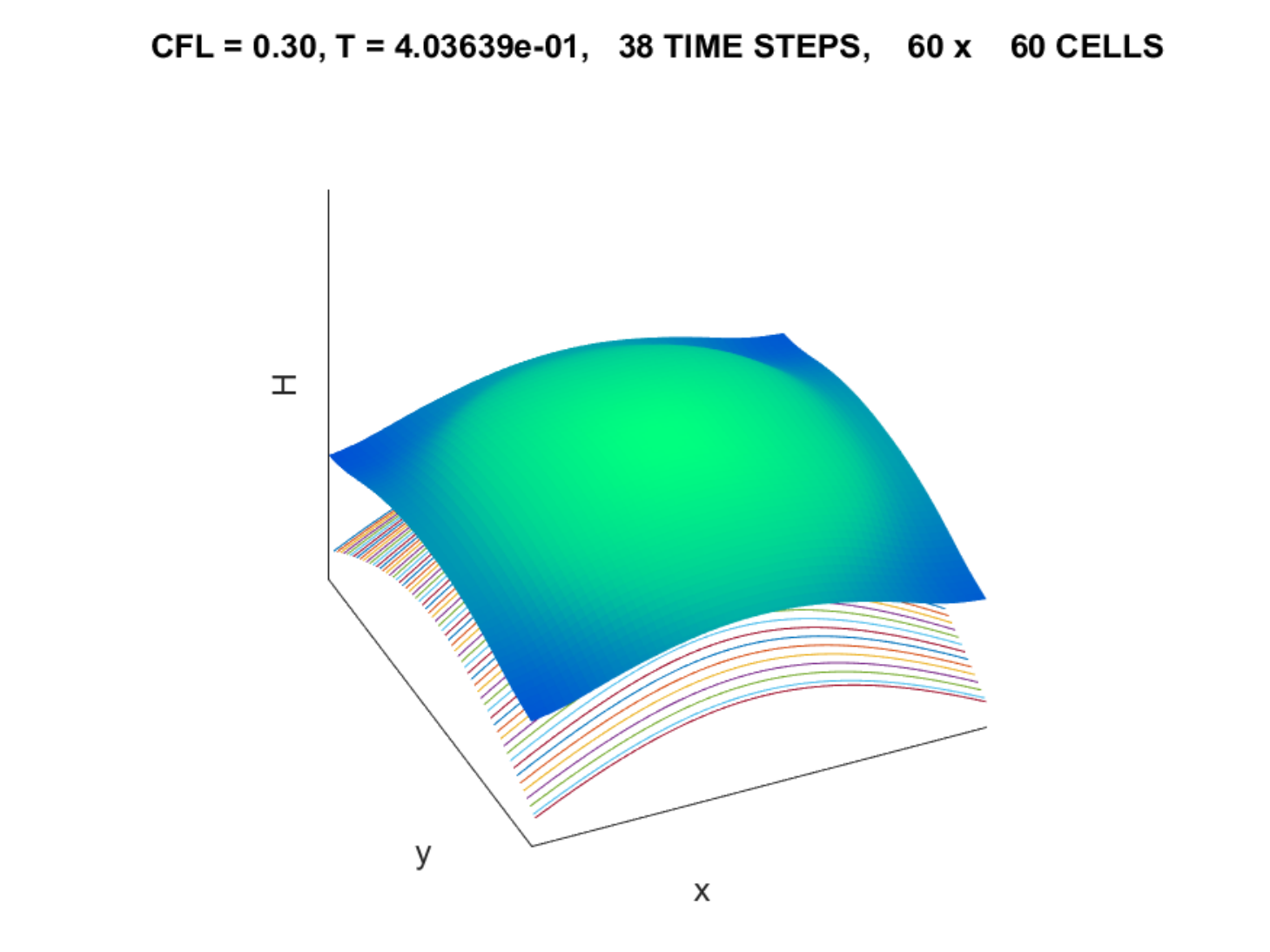}
	\includegraphics[width=0.38\textwidth]{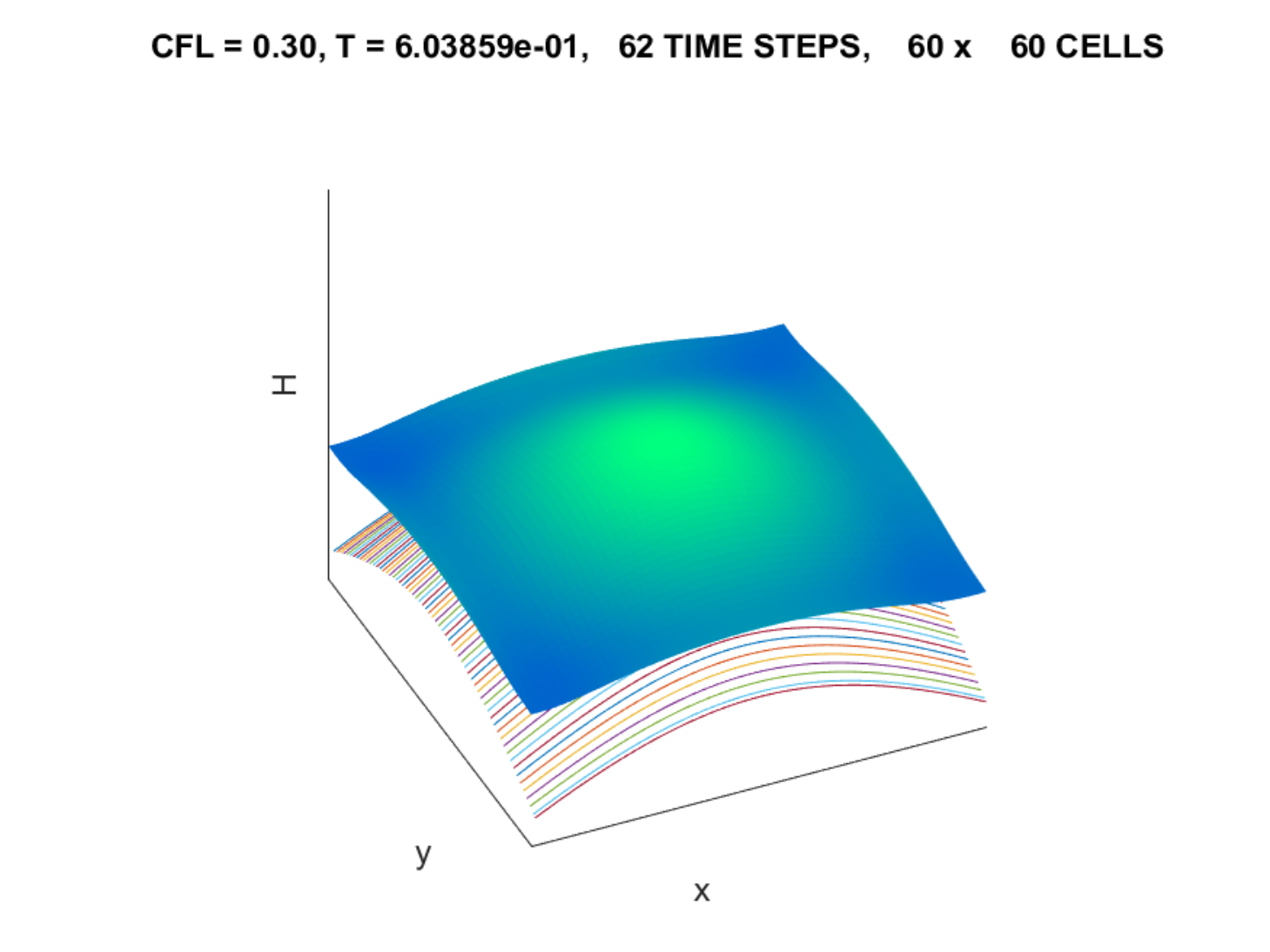}
	\caption{Hyperboloid-central-bump: evolution of the water
          depth at times $t=0.0, 0.2, 0.4, 0.6$~s.}
	\label{centralbump}
\end{figure}

\begin{figure}
	\centering
	\includegraphics[width=0.38\textwidth]{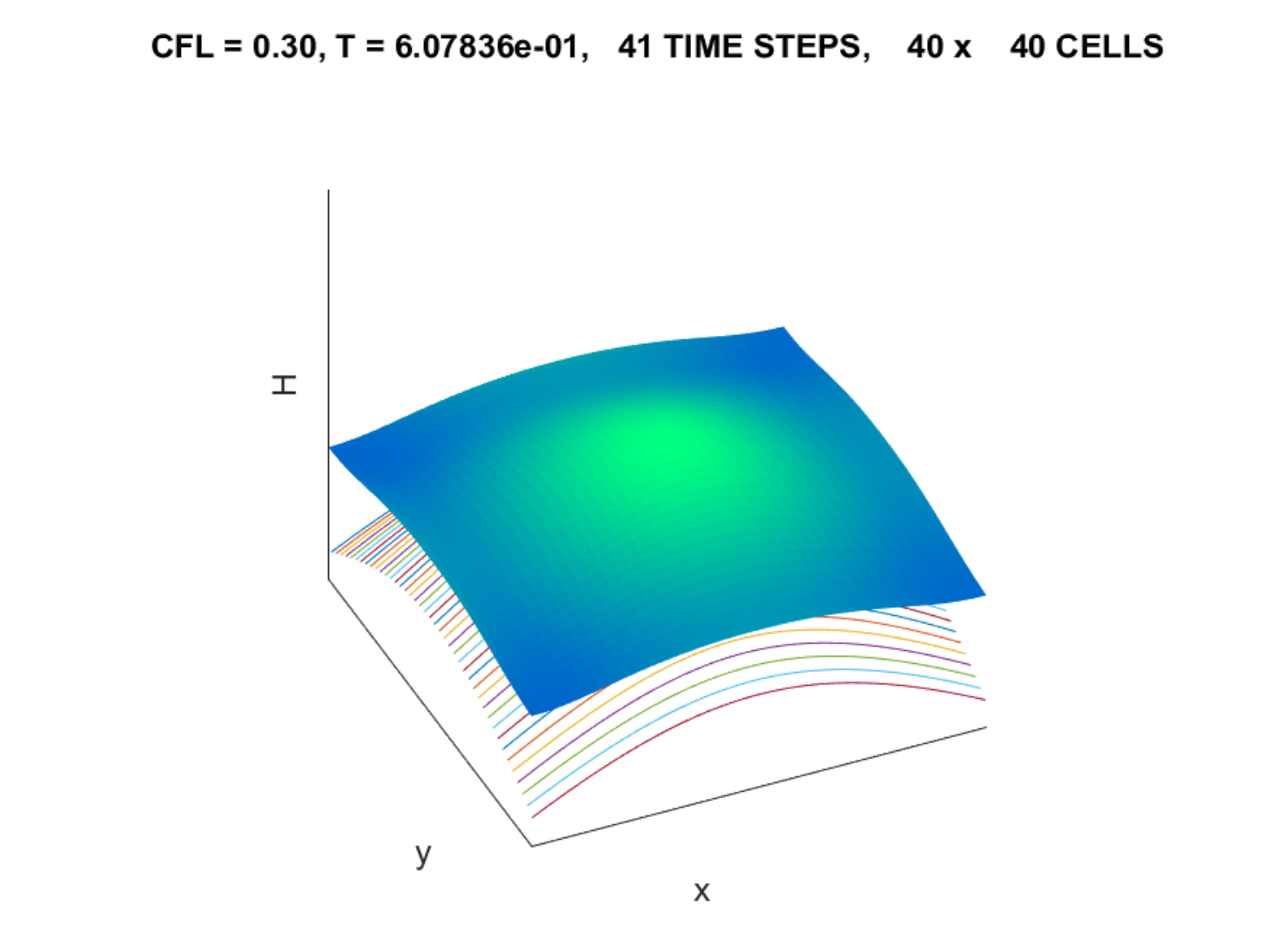}
	\includegraphics[width=0.38\textwidth]{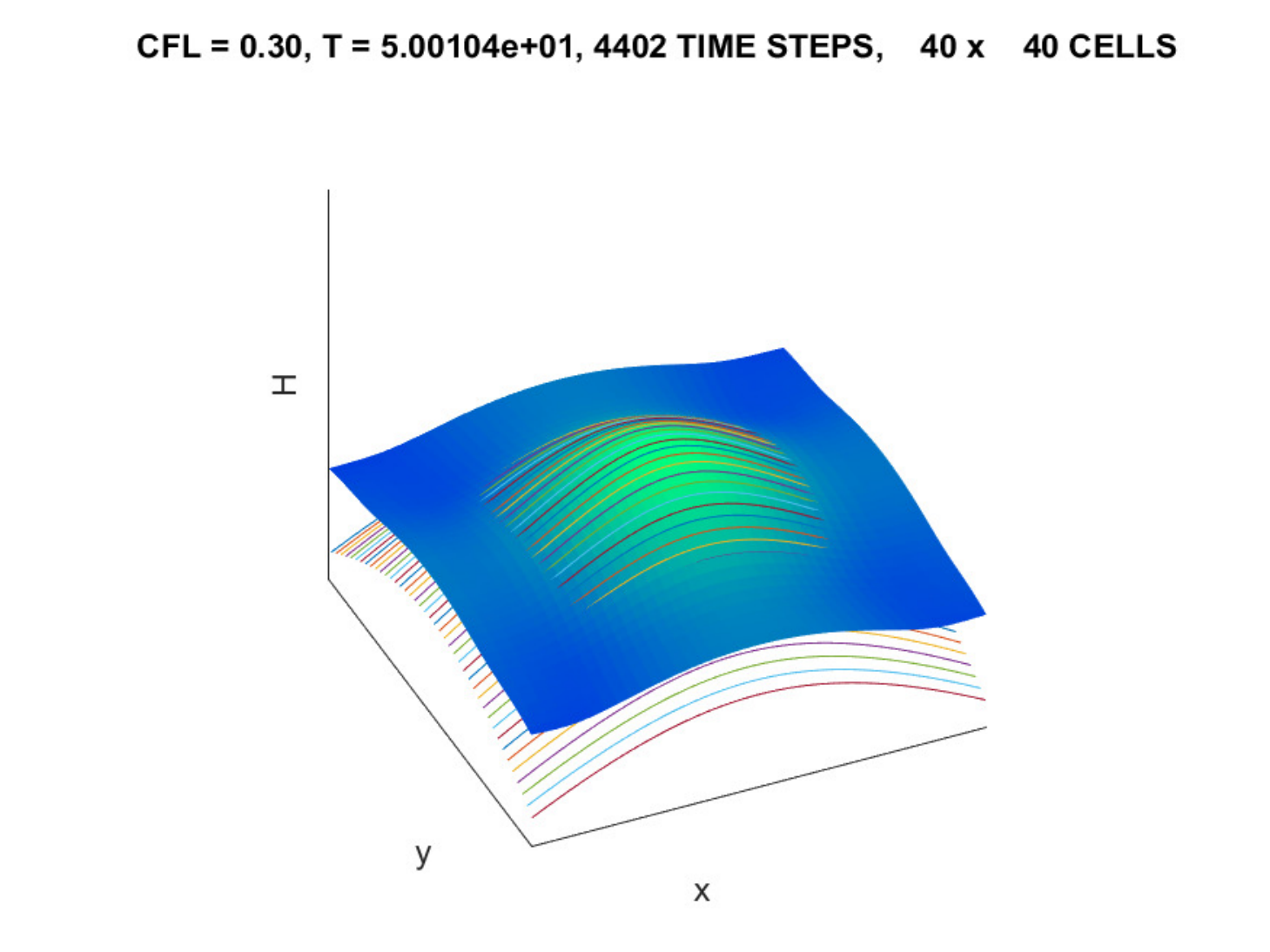}
	\includegraphics[width=0.38\textwidth]{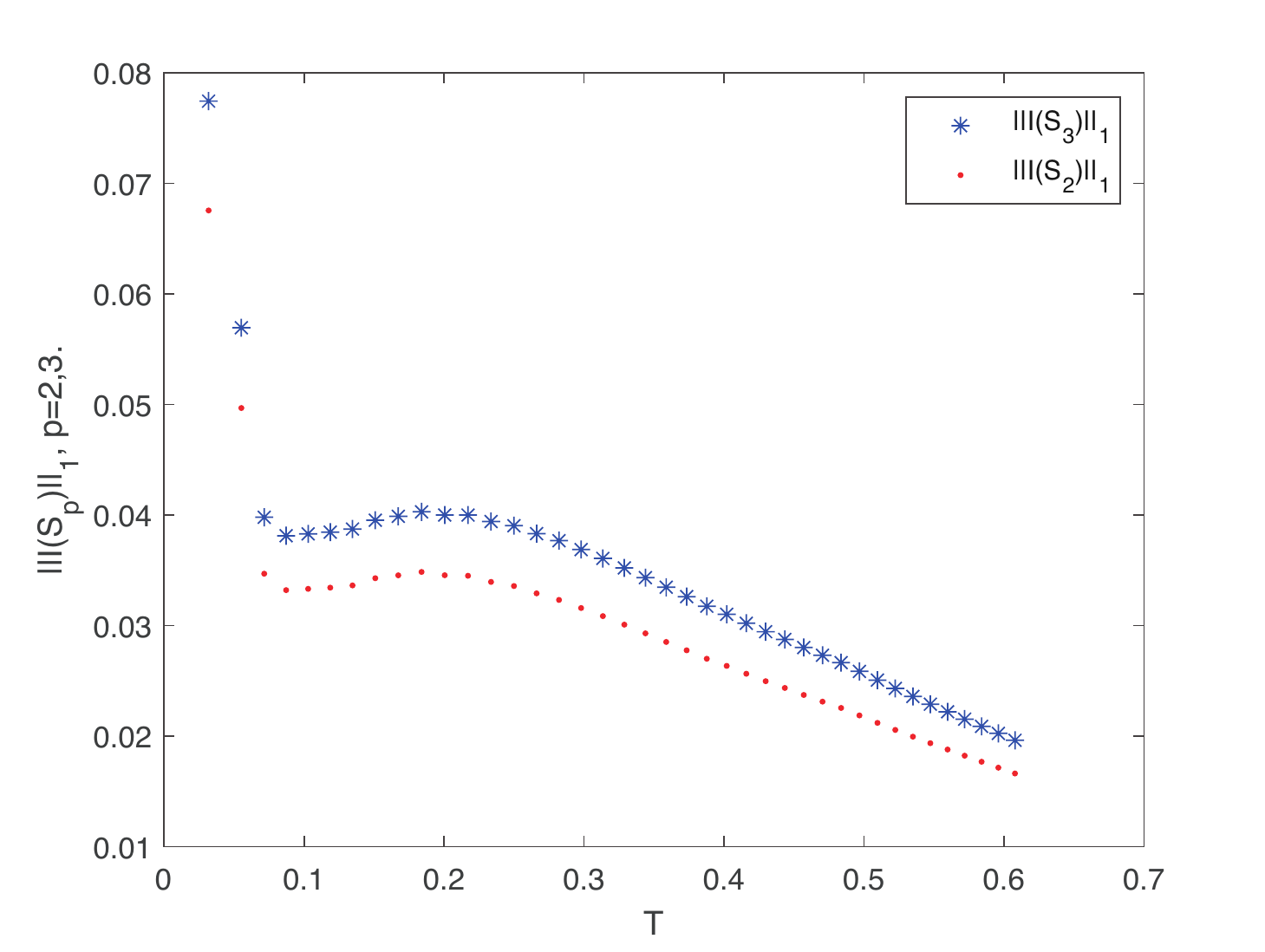}
	\includegraphics[width=0.38\textwidth]{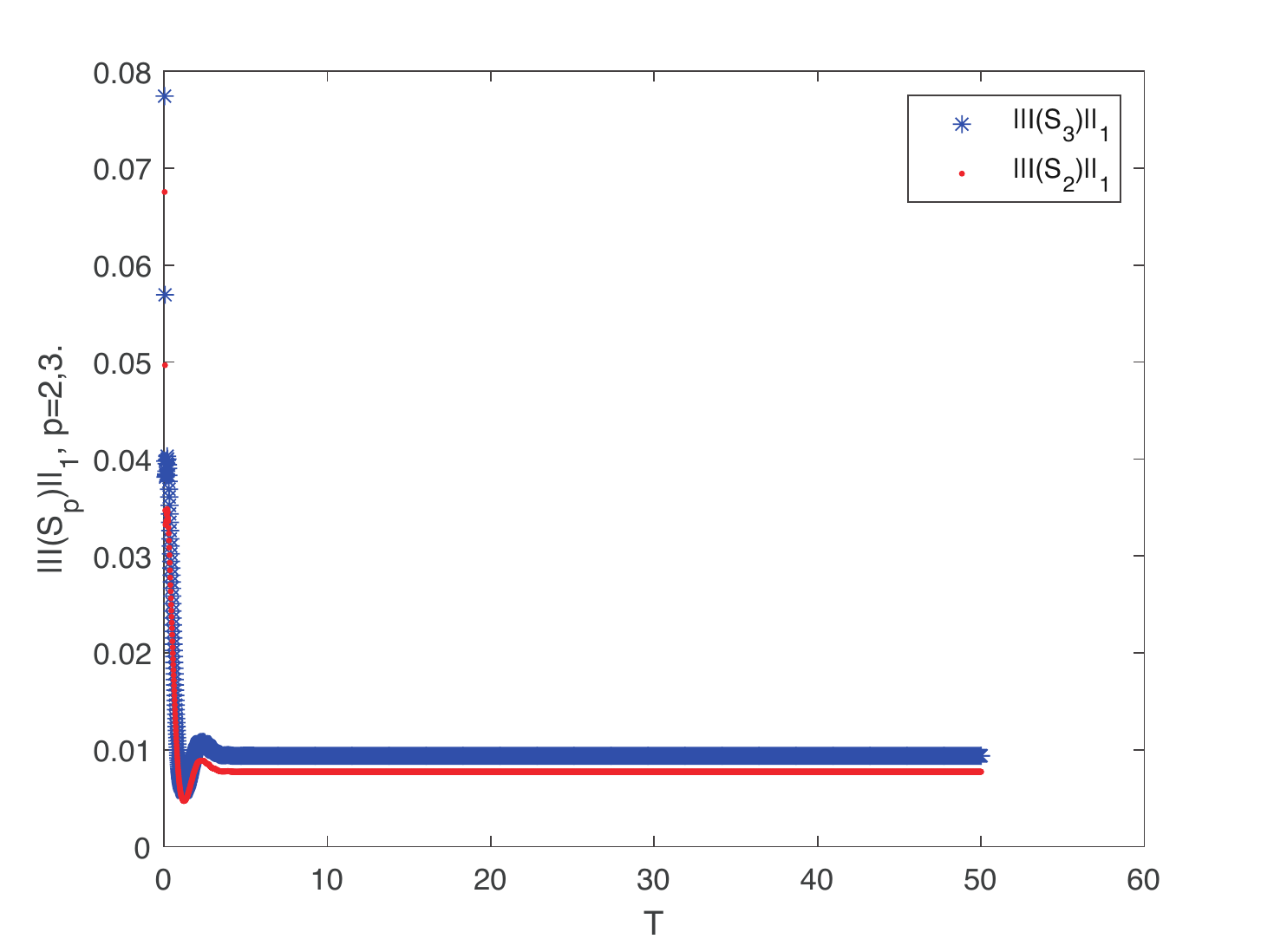}
	\caption{Hyperboloid-central-bump, 40$\times$40 cells: water depth (top panels) and
          plots of the $L^1$-norm of $I(\Source^2)$ and $I(\Source^3)$ (bottom panels), at times
          $t=0.6$~s (left) and $t=50$~s (right). Computational time for $t=50$: 6.46~s.} 
	\label{CI7b}
\end{figure}

The time evolution of the normal water depth at different times
($t=$0.0, 0.2, 0.4, 0.6~s) is shown in figure~\ref{centralbump} on a
grid formed by $60\times 60$ cells.
Steady-state results
are reported in
figure \ref{CI7b} for a $40\times 40$ cells grid.
Again the results clearly show the effectiveness of the proposed scheme.

\paragraph{Fully 3D surface}
The final surface
is designed to highlight the behavior of the proposed scheme when
positive and negative curvatures are present and to test its
effectiveness in handling complex wetting and drying patterns.
The height function is given by:
\begin{equation*}
  \BSM(x,y) = -\frac{1}{500}(y)^3 - \frac{1}{100}y(x)^2
\end{equation*}
on a domain $\SubsetU=[-4,4]\times[-10,10]\subset\REAL^2$. The dam is
 located at $y=-8.5$~m and the initial condition for the depth
is set to $2.0$~m of water upstream and $1.0$~m of water
downstream. This is the only example where we used an outlet boundary
condition at $y=10$~m and no-flow conditions for the other boundaries.

\begin{figure}
  \centering
  \includegraphics[width=0.38\textwidth]{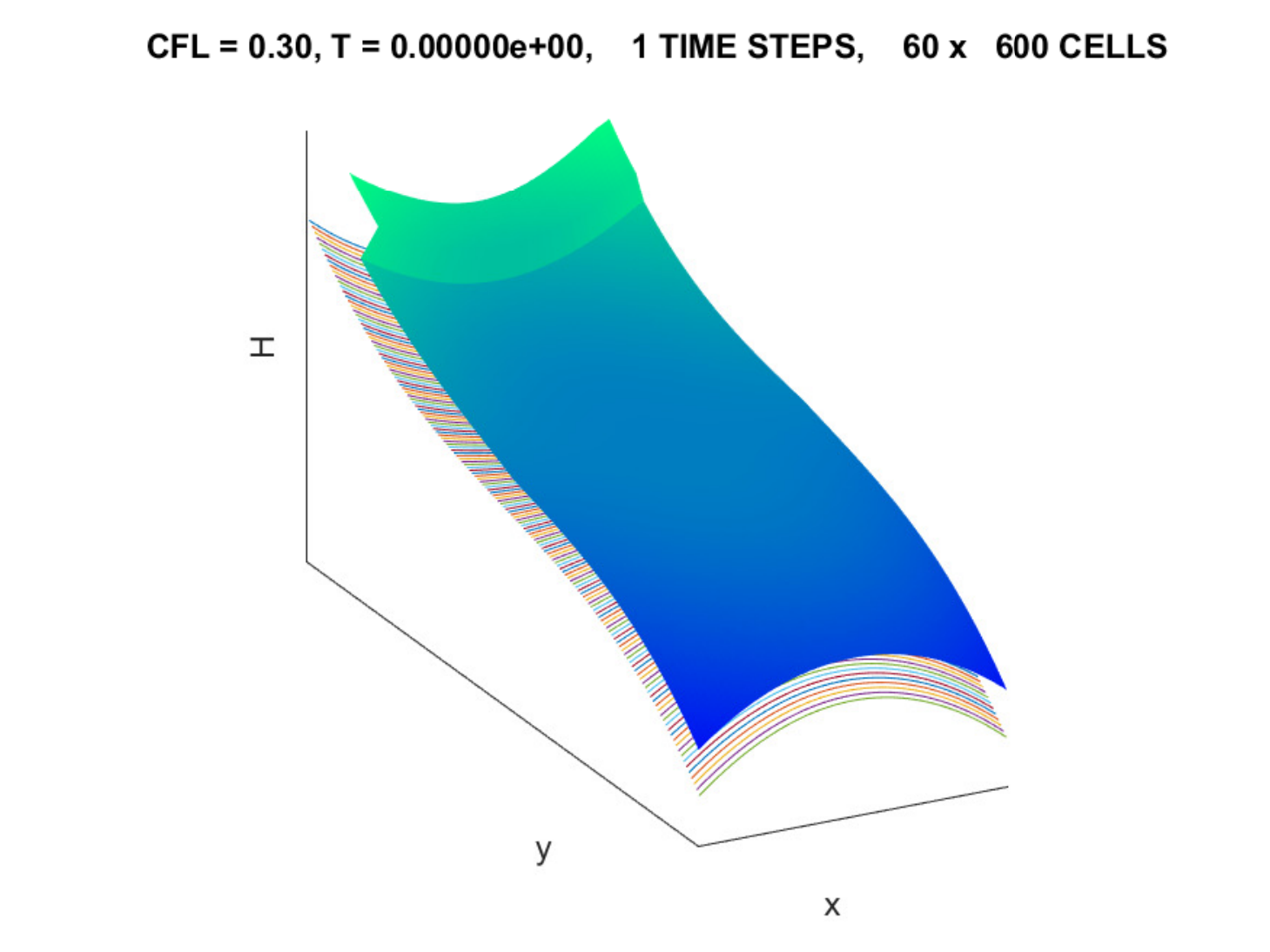}
  \includegraphics[width=0.38\textwidth]{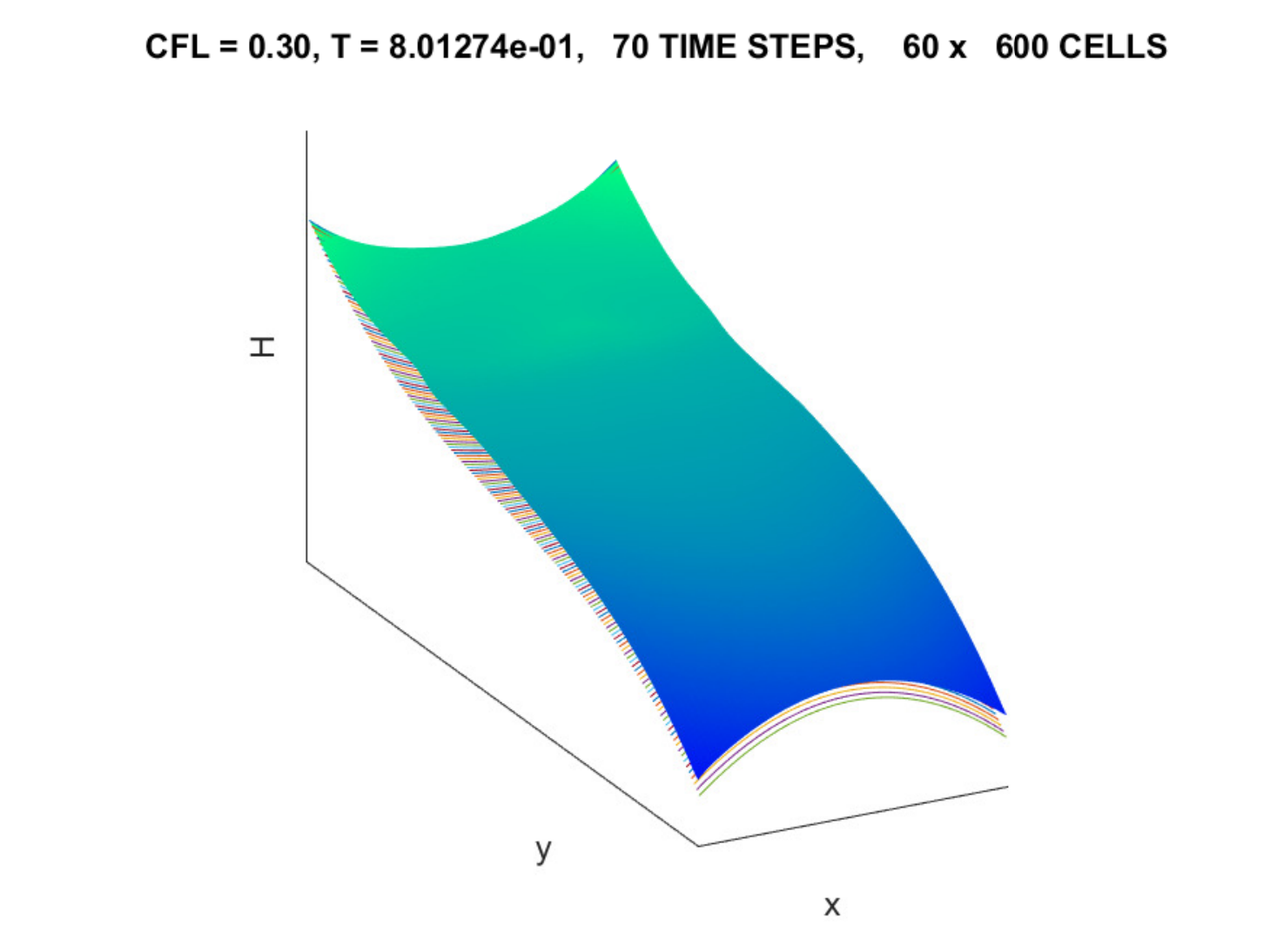}
  \includegraphics[width=0.38\textwidth]{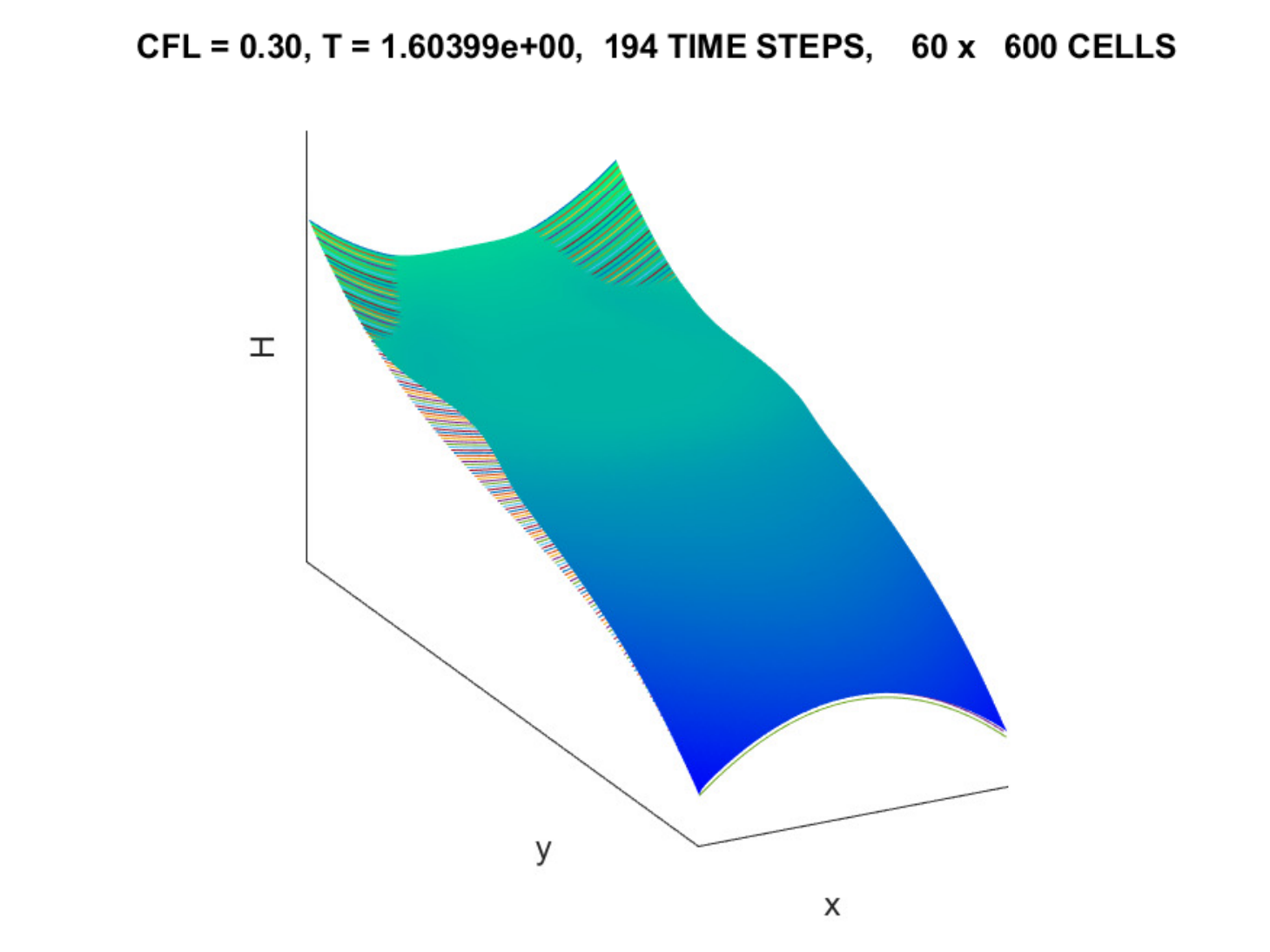}
  \includegraphics[width=0.38\textwidth]{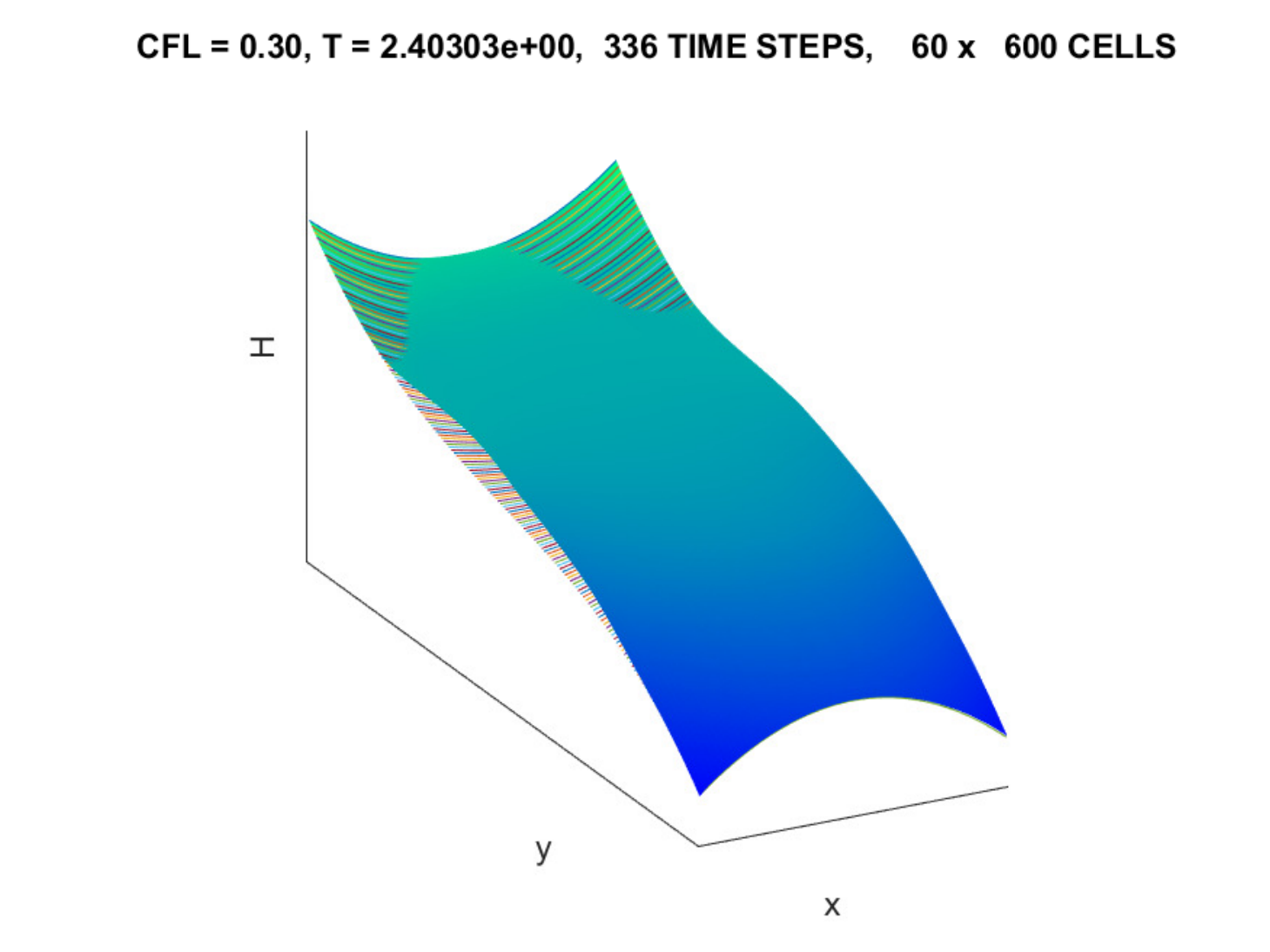}
  \caption{Fully 3D surface: evolution of the water depth at times
    $T=$0.0, 0.8, 1.6, 2.4~s.}
  \label{Fully3D}
\end{figure}

\begin{figure}
	\centering
	\includegraphics[width=0.38\textwidth]{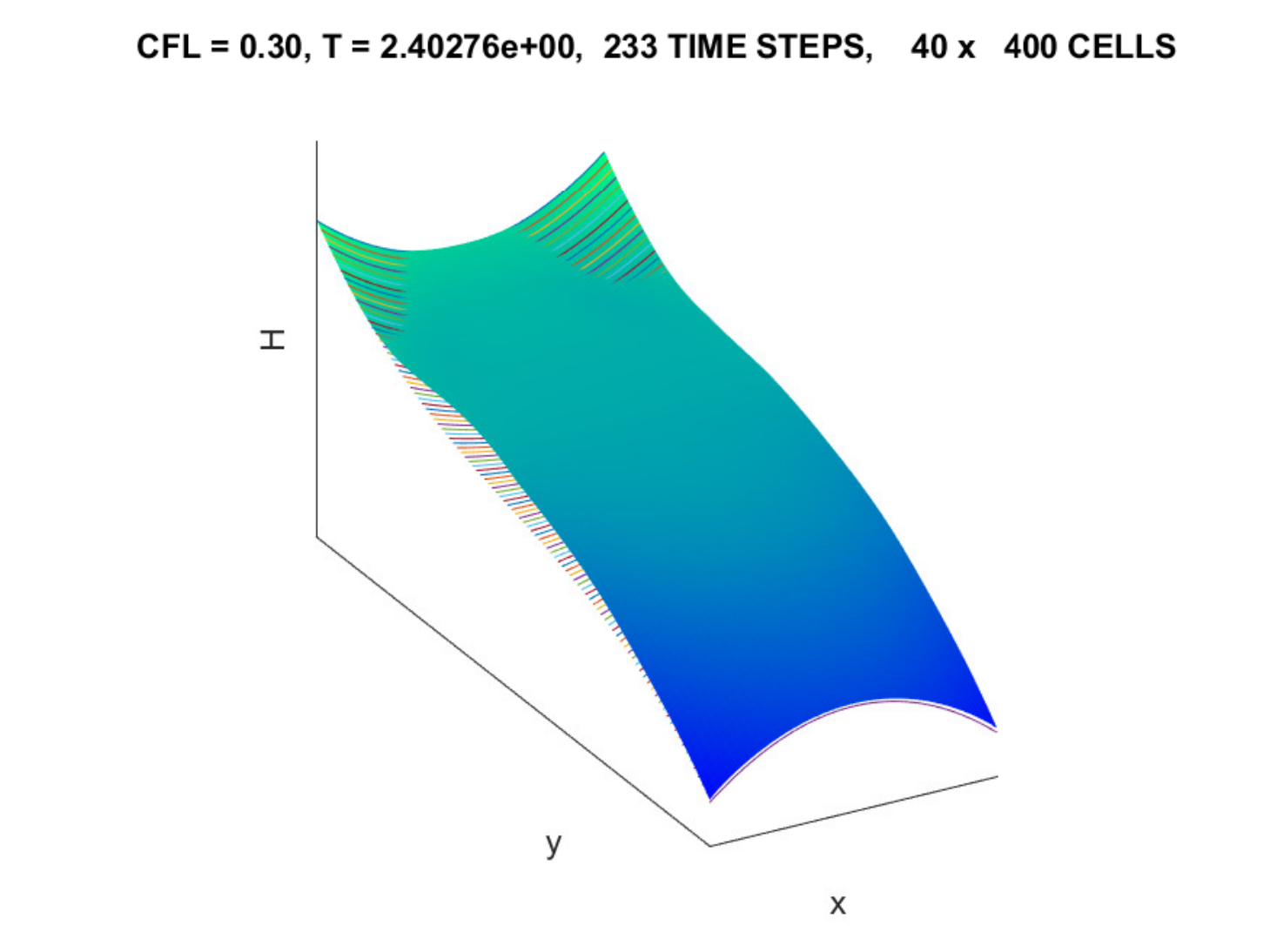}
	\includegraphics[width=0.38\textwidth]{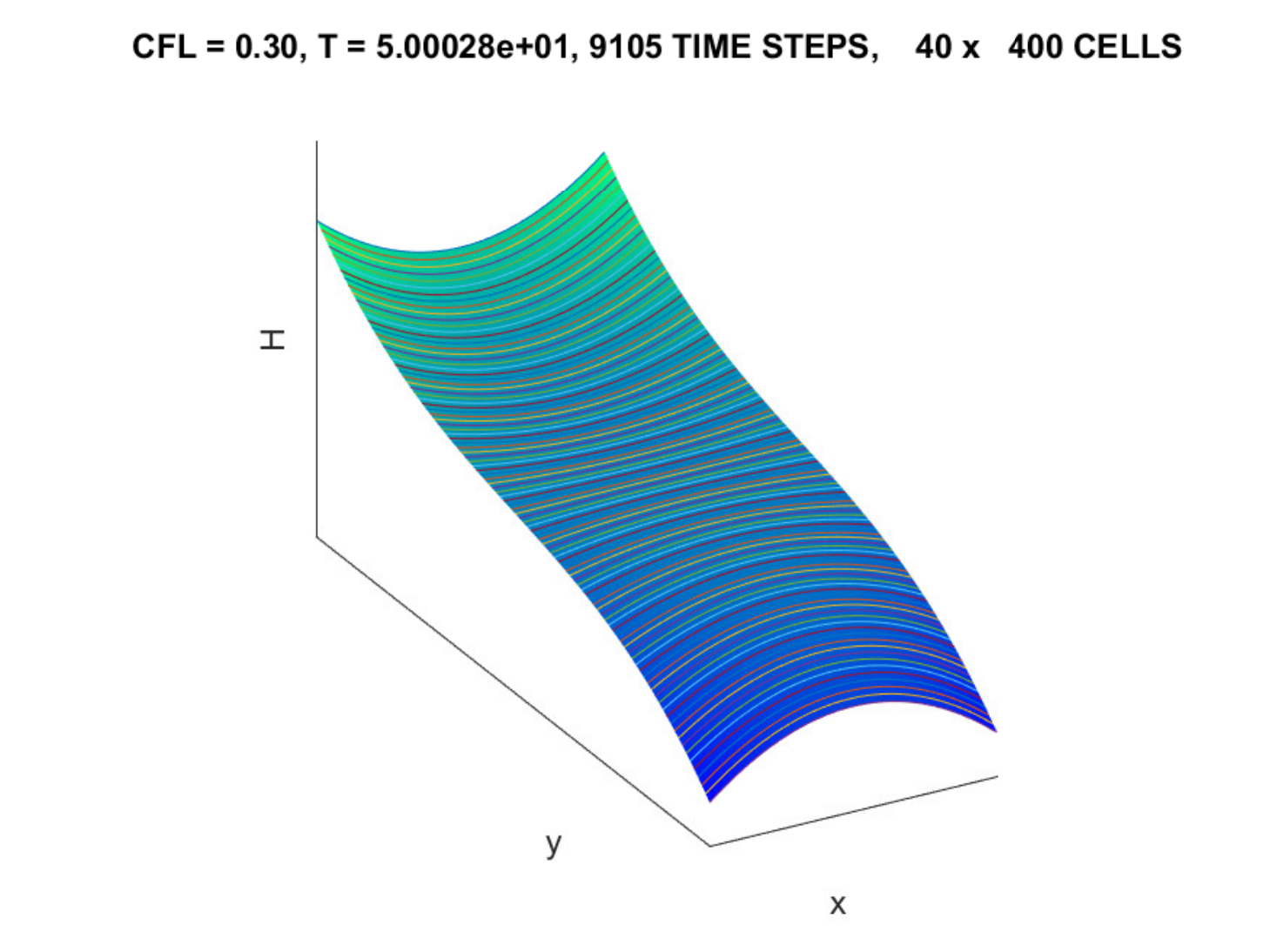}
	\includegraphics[width=0.38\textwidth]{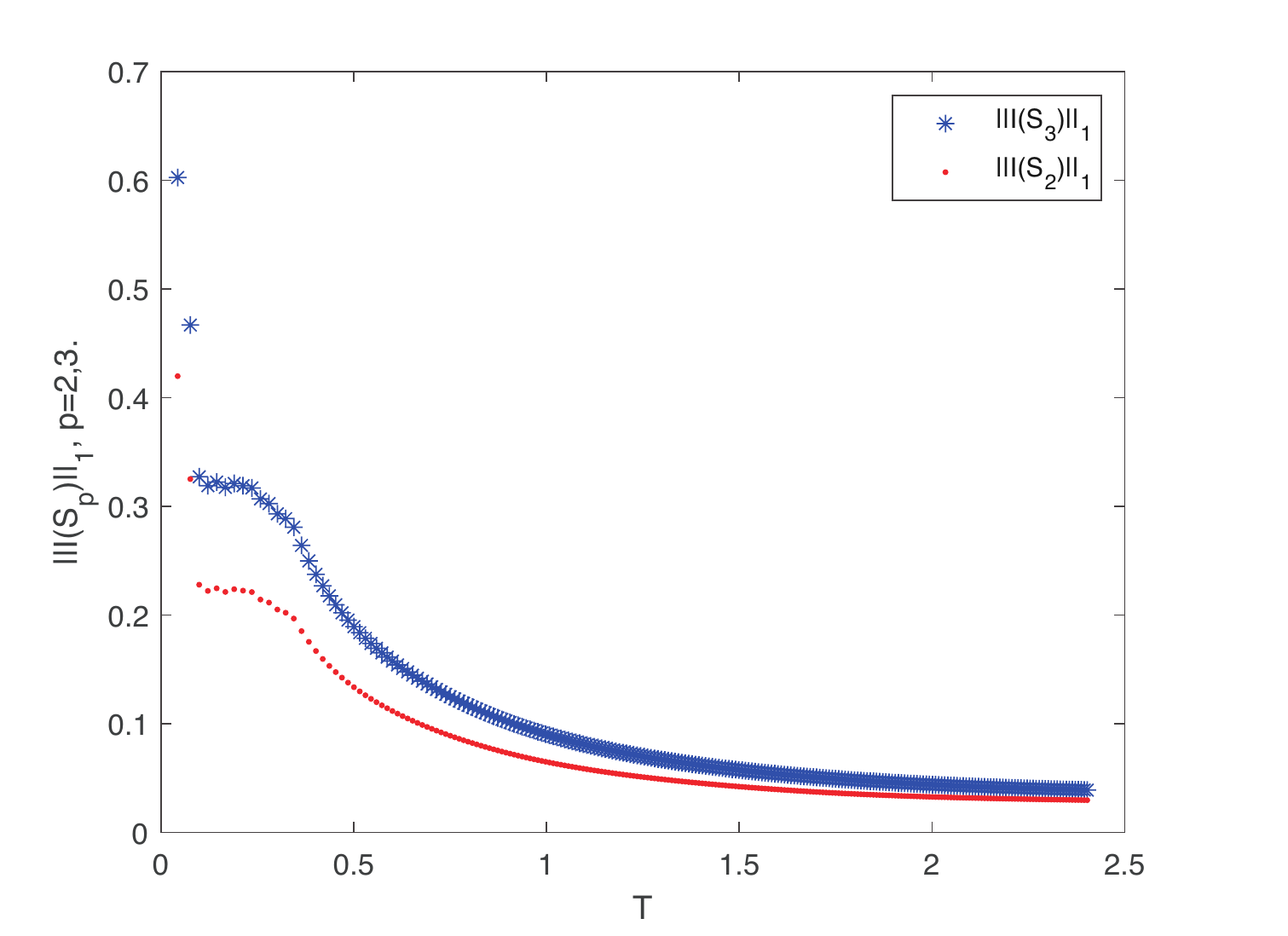}
	\includegraphics[width=0.38\textwidth]{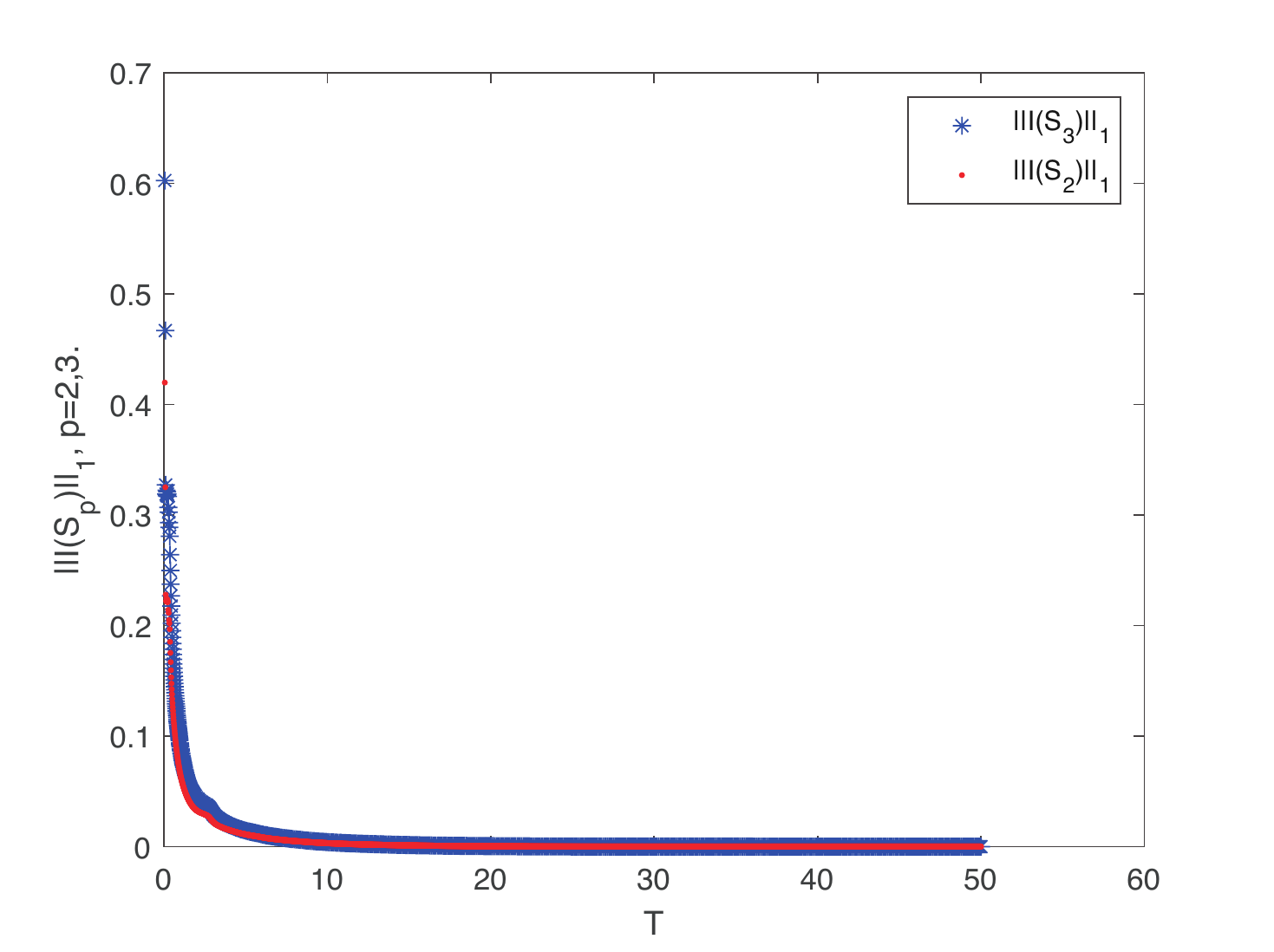}
	\caption{Fully 3D surface, 40$\times$400 cells: water depth (top panels) and
          plots of the $L^1$-norm of $I(\Source^2)$ and $I(\Source^3)$ (bottom panels), at times
          $t=0.6$~s (left) and $t=50$~s (right). Computational time for $t=50$: 45.65~s.} 
	\label{CI5}
\end{figure}

Figure~\ref{Fully3D} shows the evolution of the water depth from the
initial solution, $t=0.0$~s, and at times $t=$0.8, 1.6, 2.4~s, on a
grid that considers $60\times 600$ cells.
 The long-time numerical tests shown in figure \ref{CI5} verify the
 well-balance property of our approach.

\paragraph{Discontinuous initial conditions}

Finally,we present a test case in which strongly discontinuous
initial conditions are present.
The domain is a rectangle $[0,3]\times[0,10]\subset\REAL^2$ and the
height function describes a sloping plane, as in the test case 1.  We
define a dam-break initial condition with a square discontinuity in
the central portion of the domain centered at the point $(1.5,5)$,
with support in $\{|x - 1.5| <0.5\}\cup\{|y - 5| < 0.5\}$, and shape
given by the graph of the function $(x-1.5)^2 + (y-5)^2 - 0.8$.

\begin{figure}
	\centering
	\includegraphics[width=0.3\textwidth]{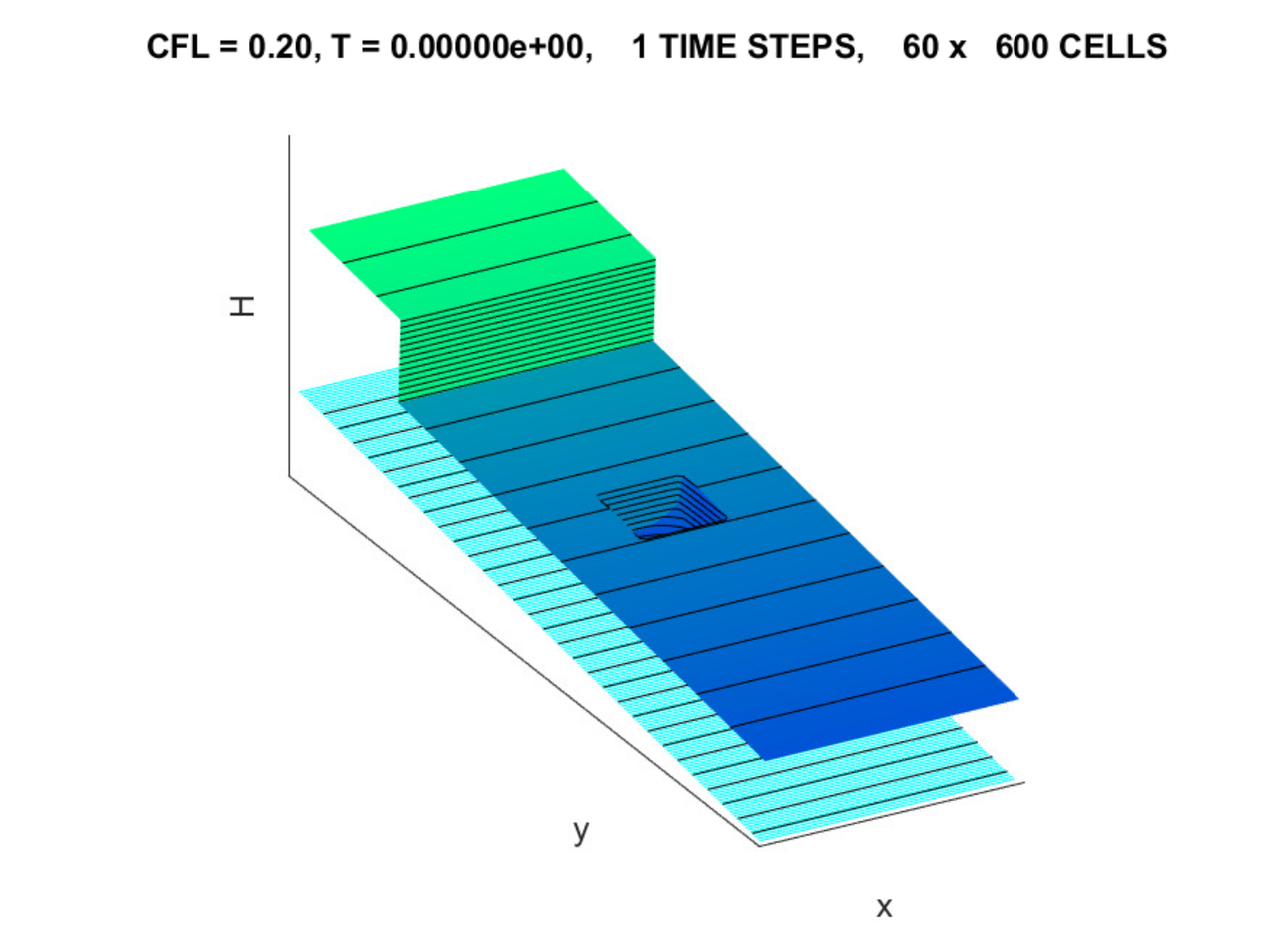}
	\includegraphics[width=0.3\textwidth]{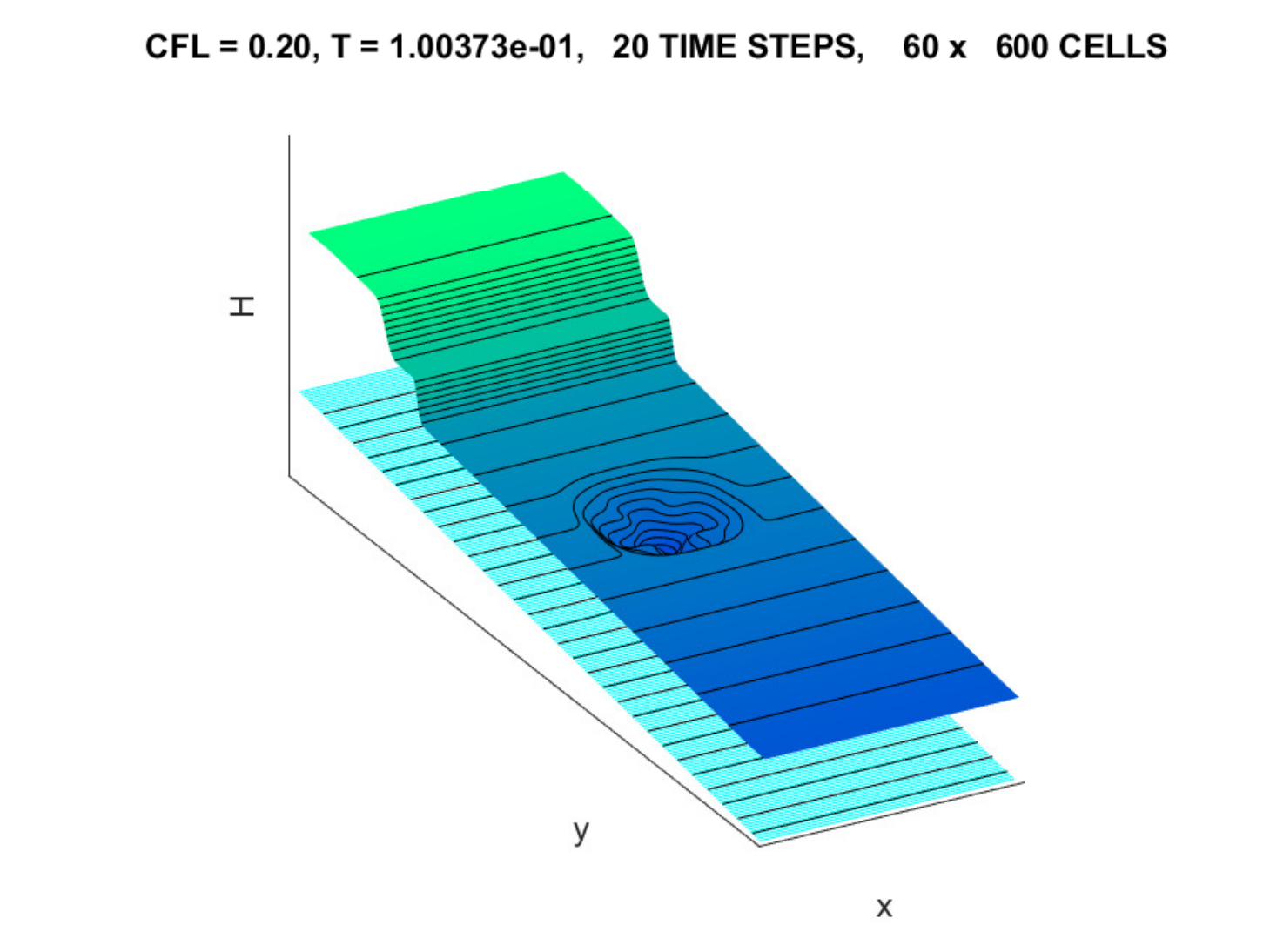}
	\includegraphics[width=0.3\textwidth]{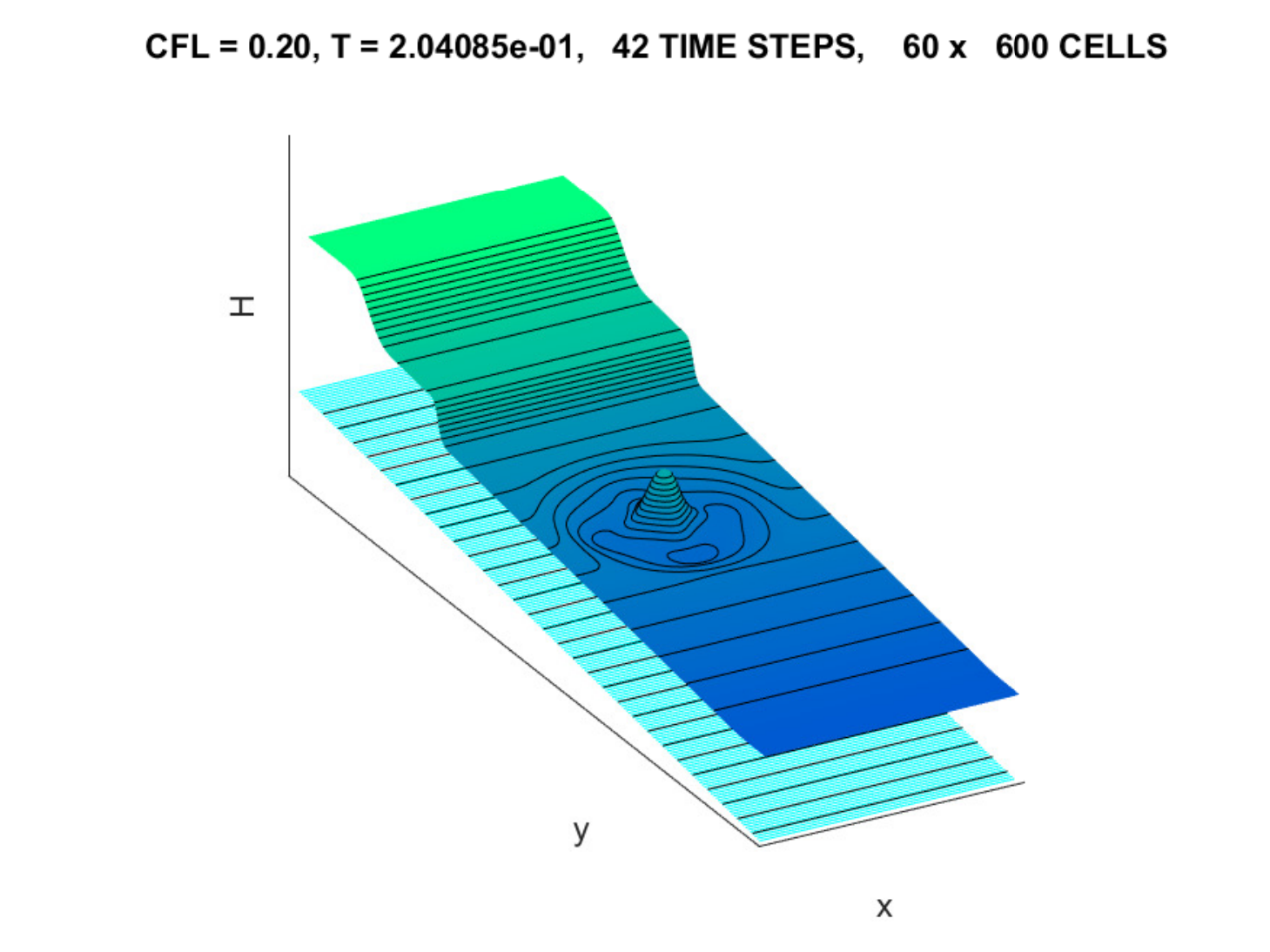}
	\includegraphics[width=0.3\textwidth]{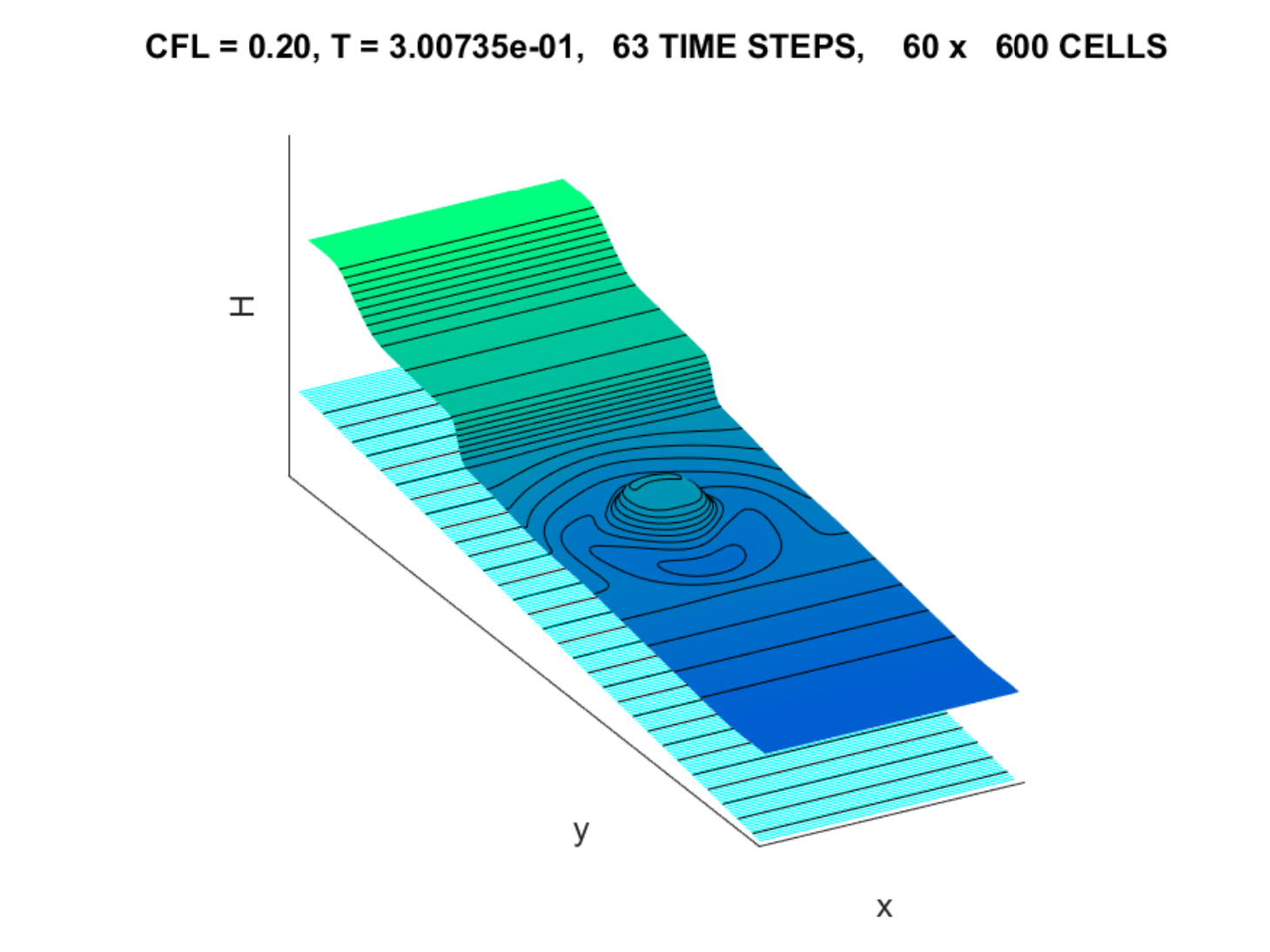}
	\includegraphics[width=0.3\textwidth]{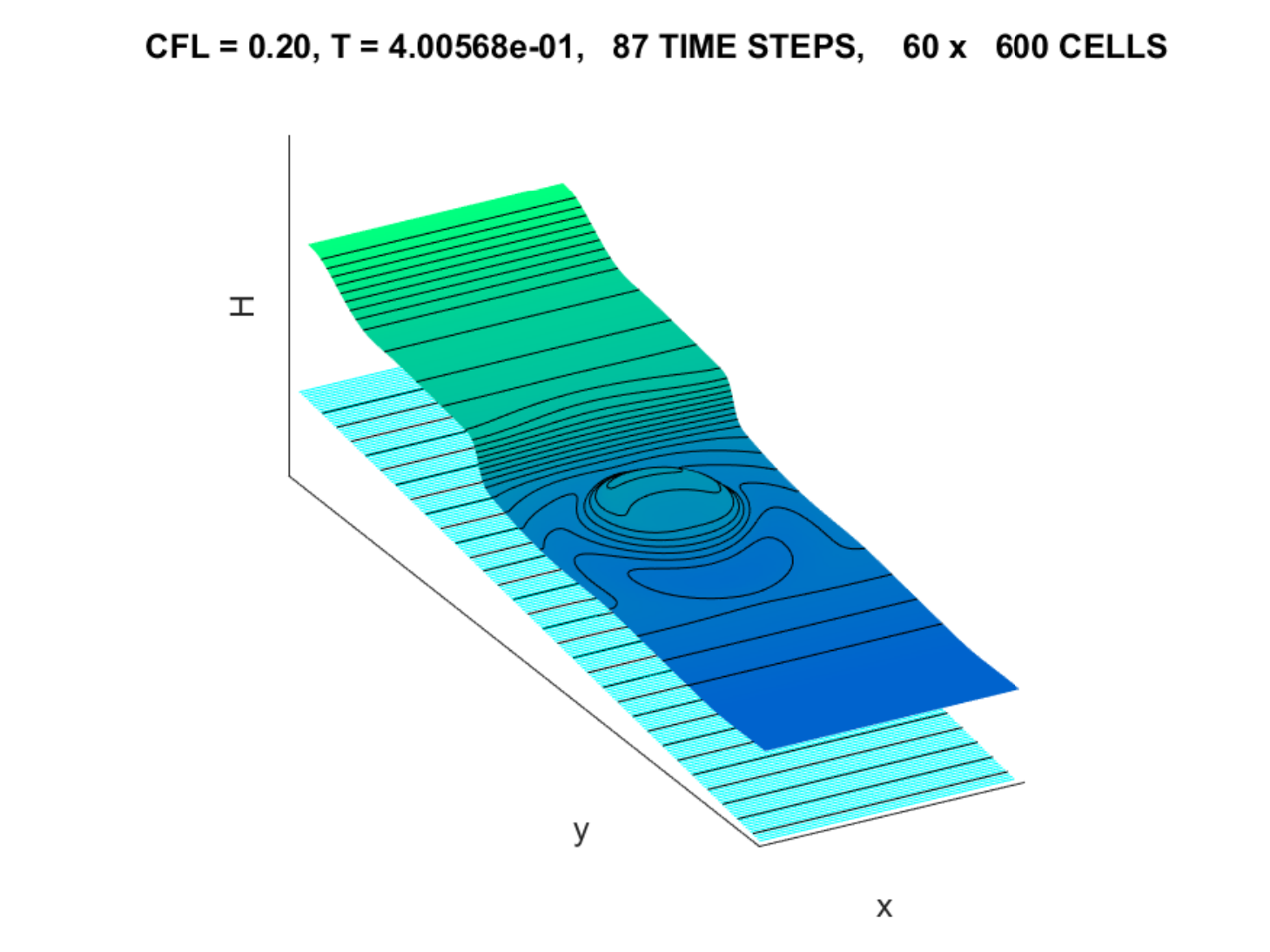}
	\includegraphics[width=0.3\textwidth]{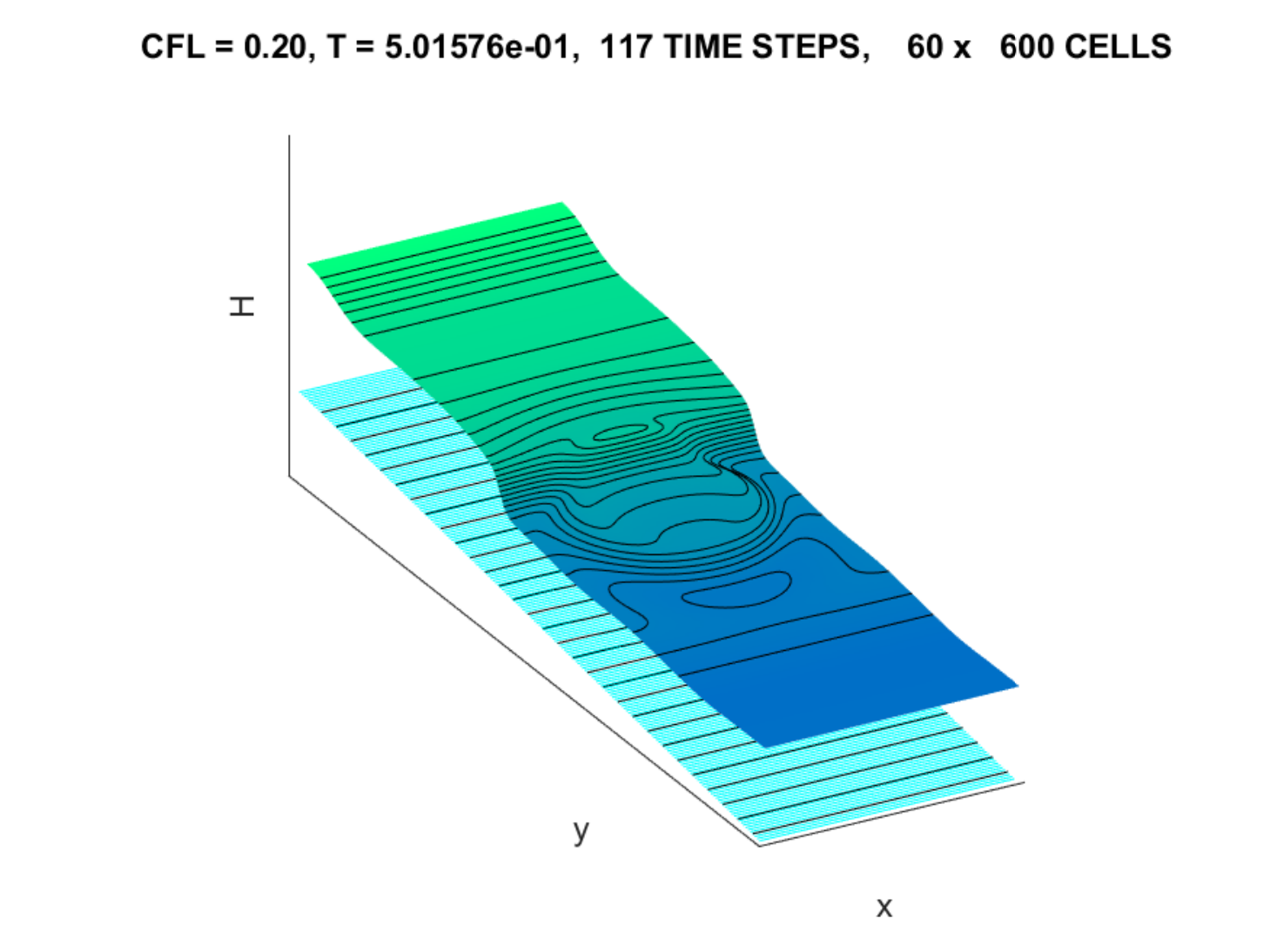}
	\includegraphics[width=0.3\textwidth]{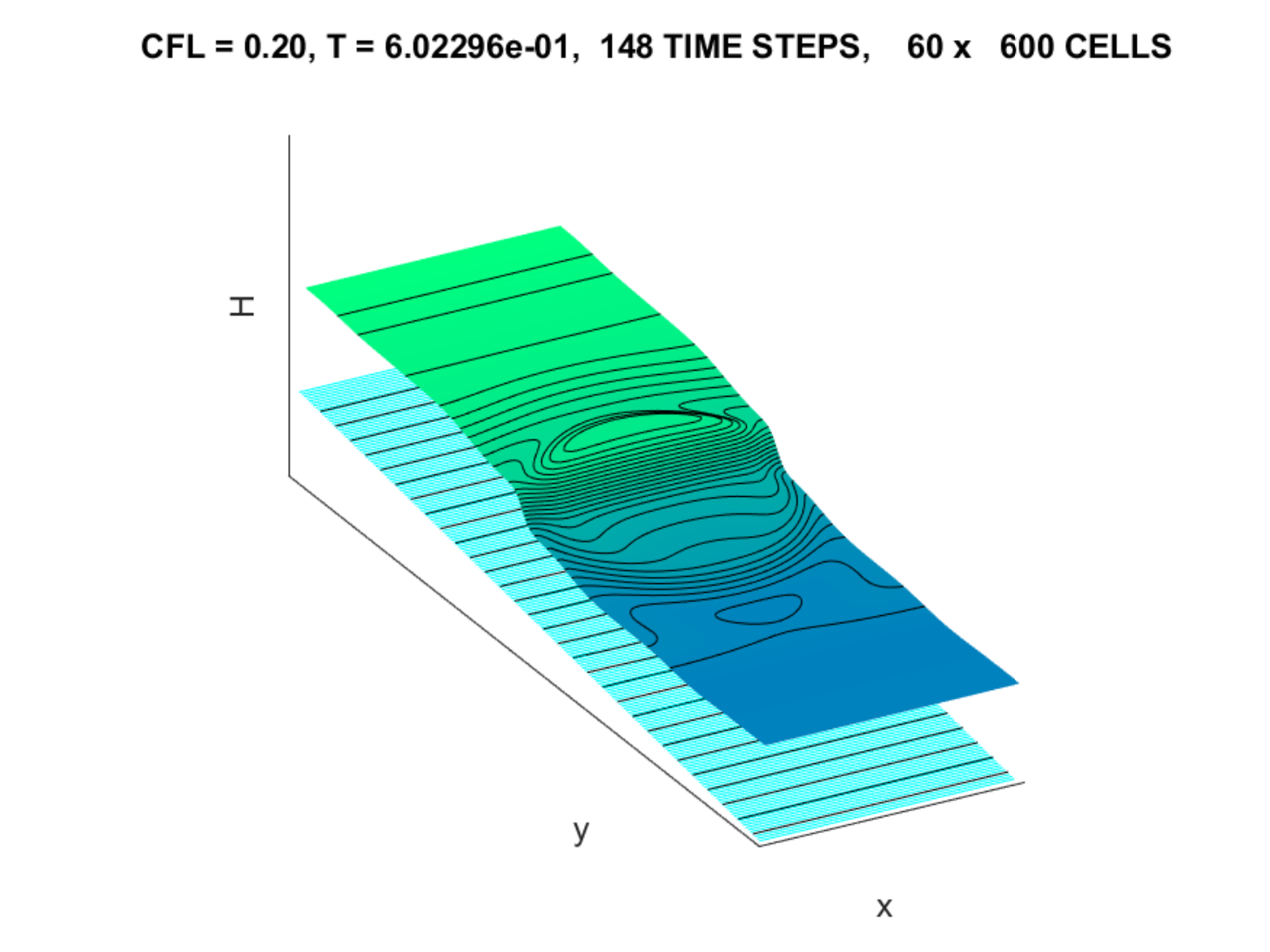}
	\includegraphics[width=0.3\textwidth]{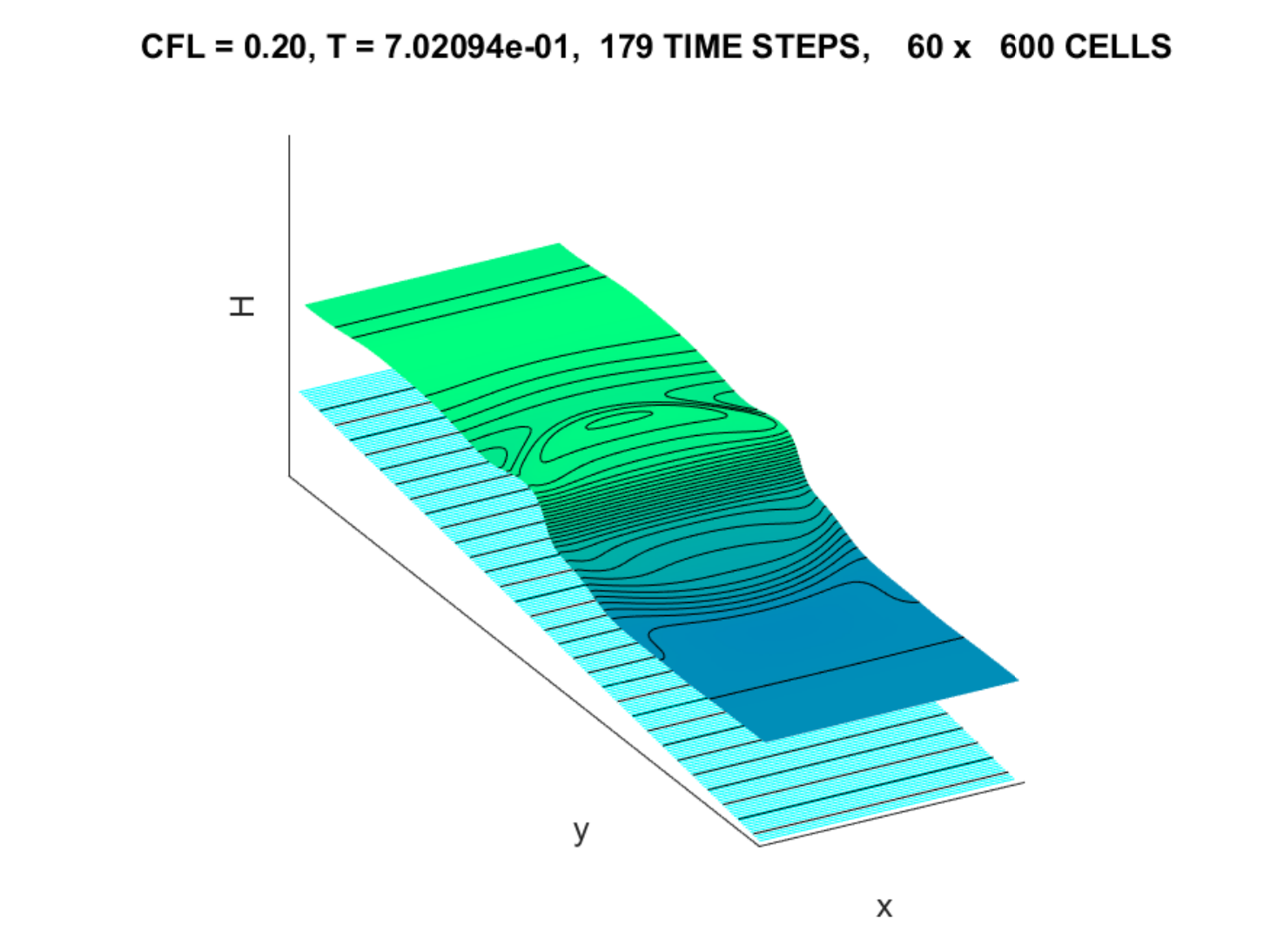}
	\includegraphics[width=0.3\textwidth]{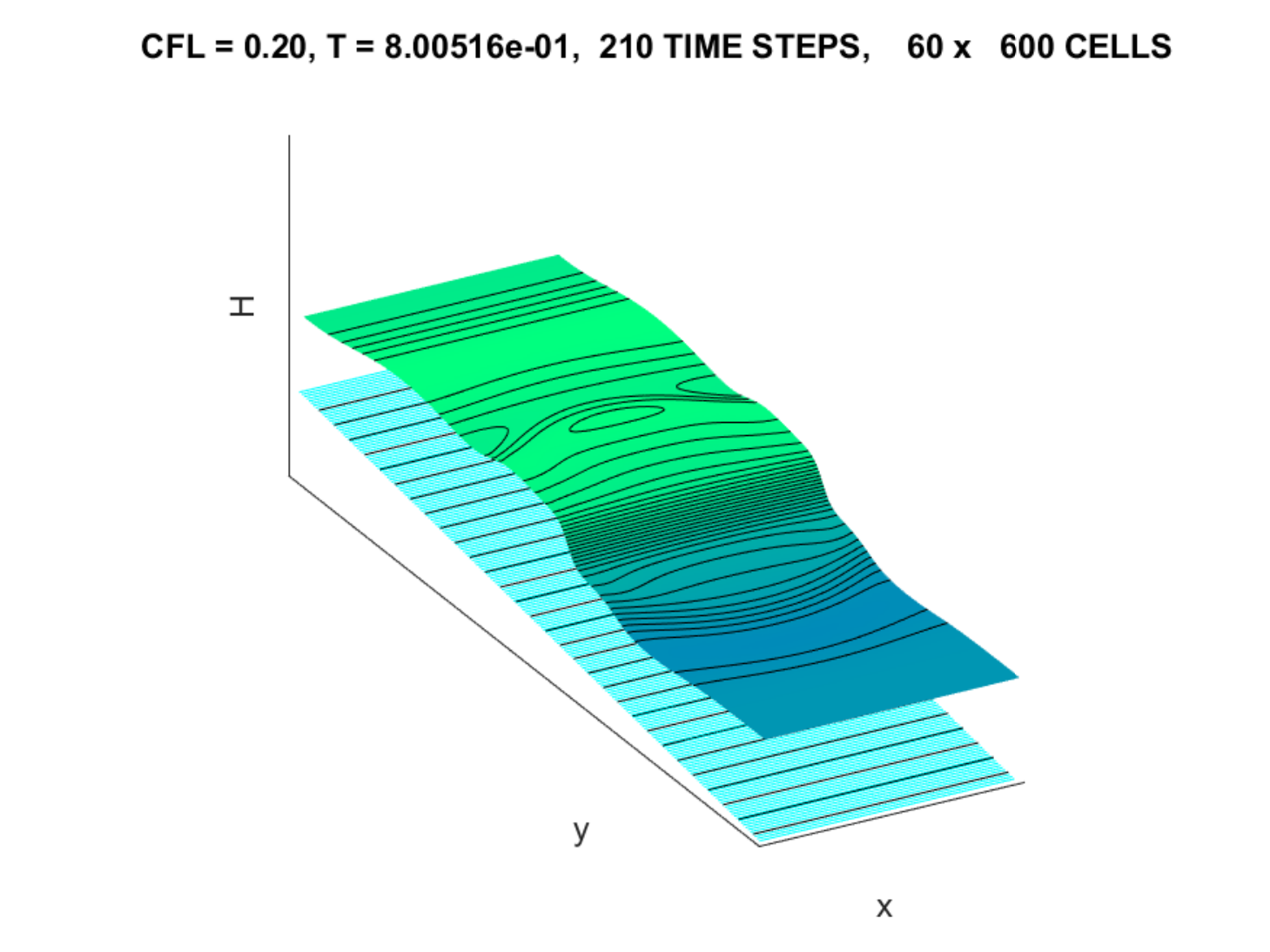}
	\includegraphics[width=0.3\textwidth]{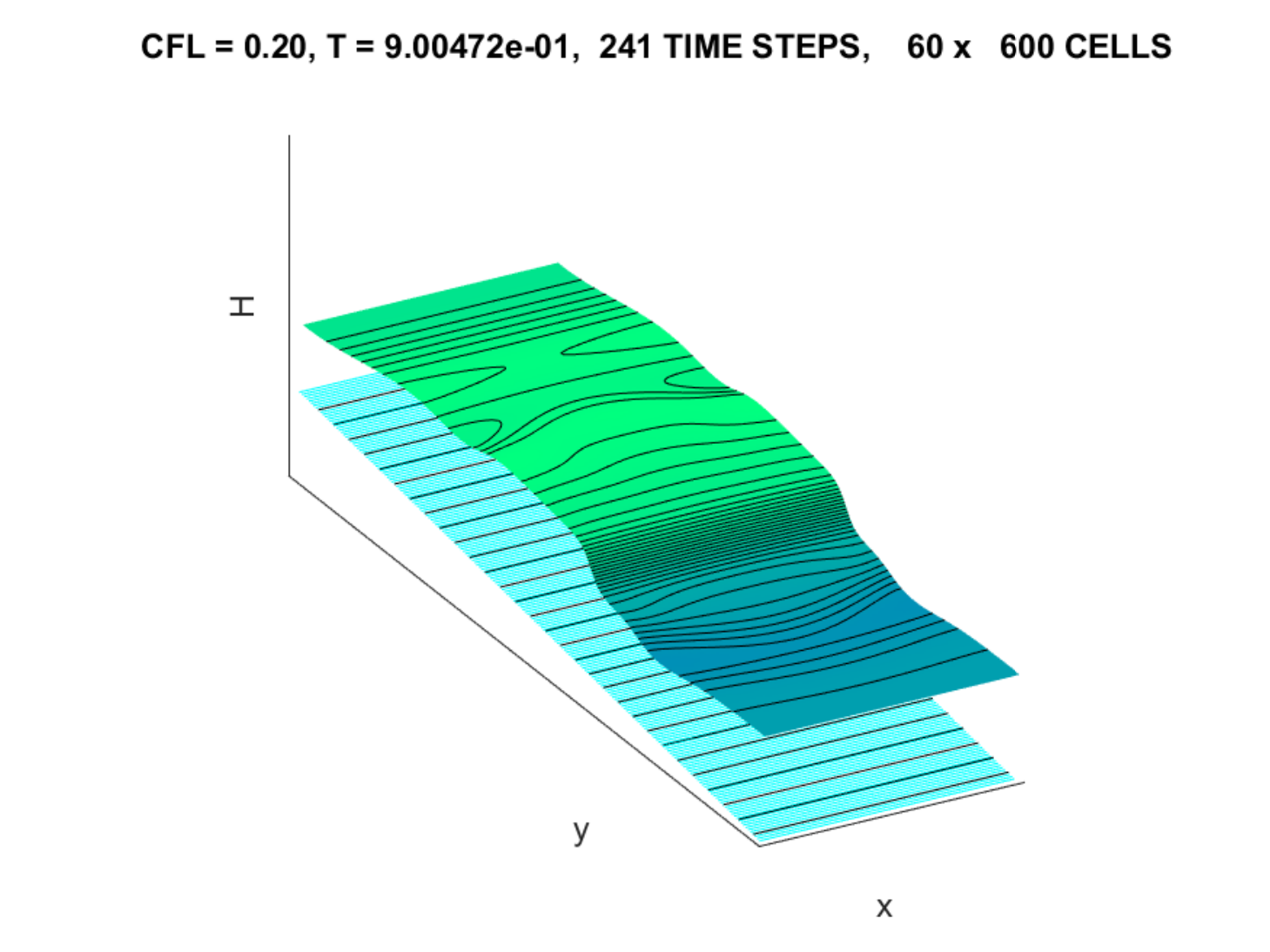}
	\includegraphics[width=0.3\textwidth]{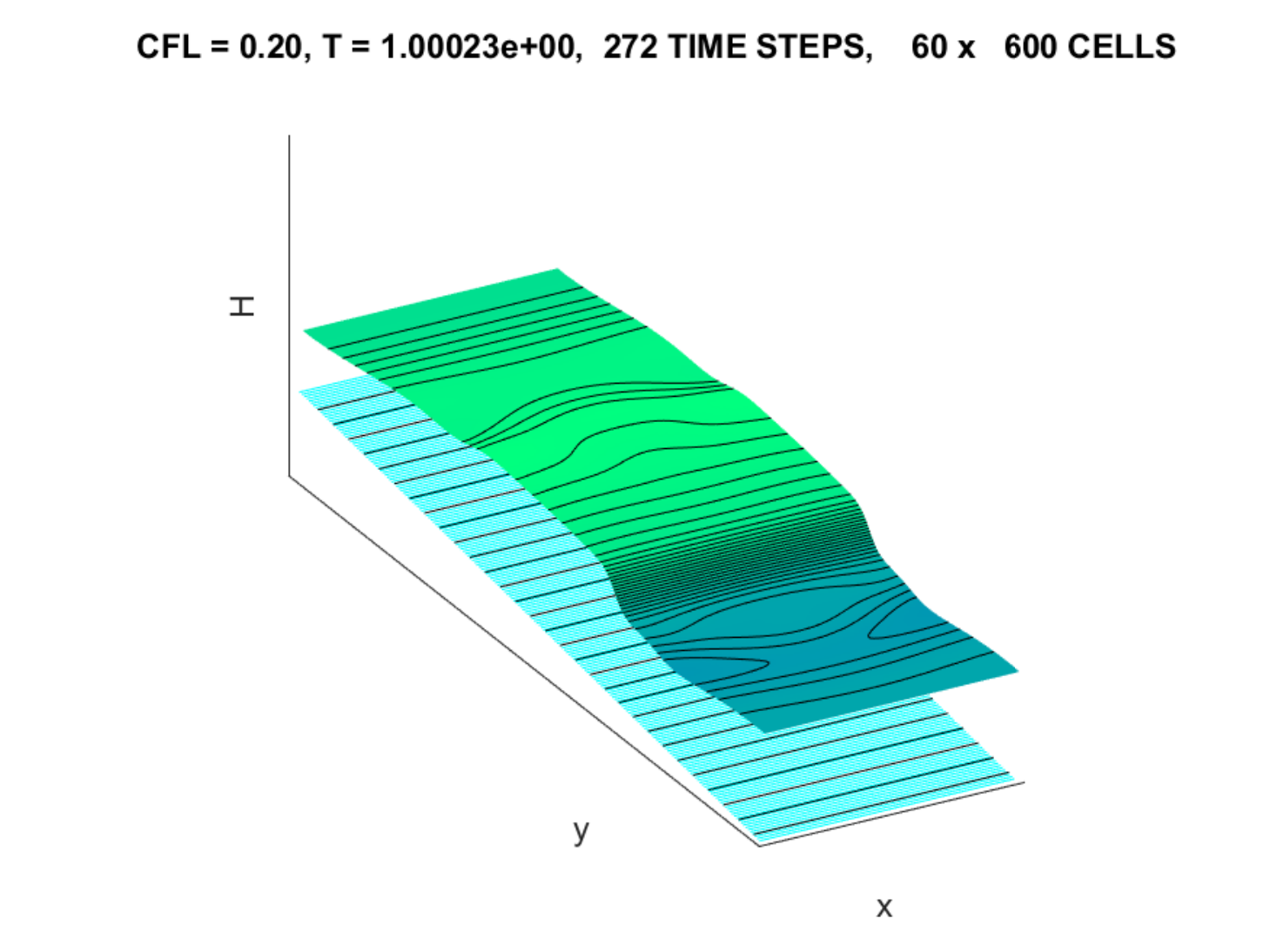}
	\includegraphics[width=0.3\textwidth]{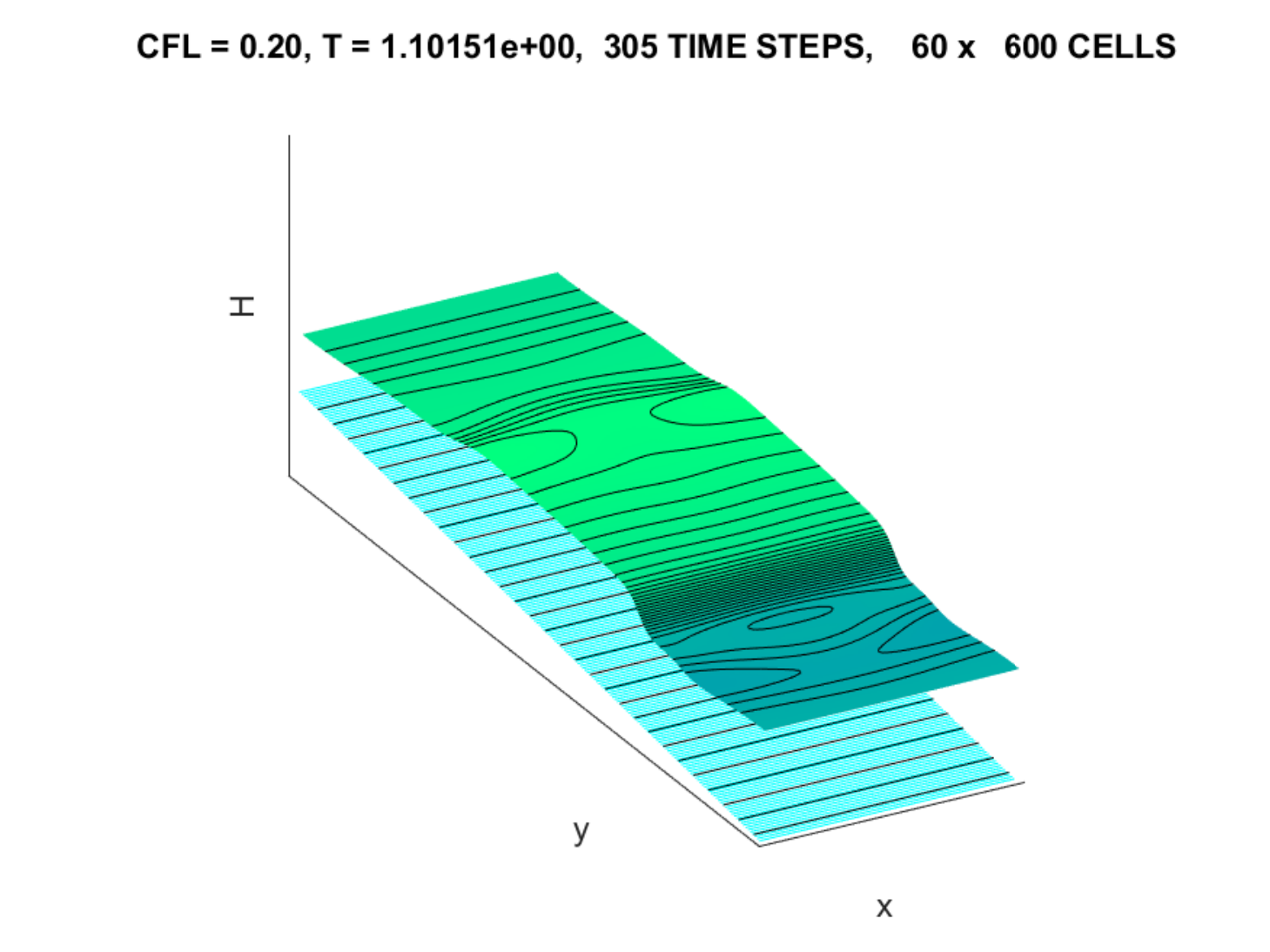}
	\caption{Sloping plane with discontinuous initial data:
          dam-break problem on a sloping plane with initial condition
          that are discontinuous in the central portion of the
          domain. The different frames show the evolution of the water
          depth from the initial solution, $t=0.0$~s (top left corner)
          to a final time of $t=1.10$~s (bottom right corner).}
	\label{discID}
\end{figure}

\begin{figure}
	\centering
	\includegraphics[width=0.38\textwidth]{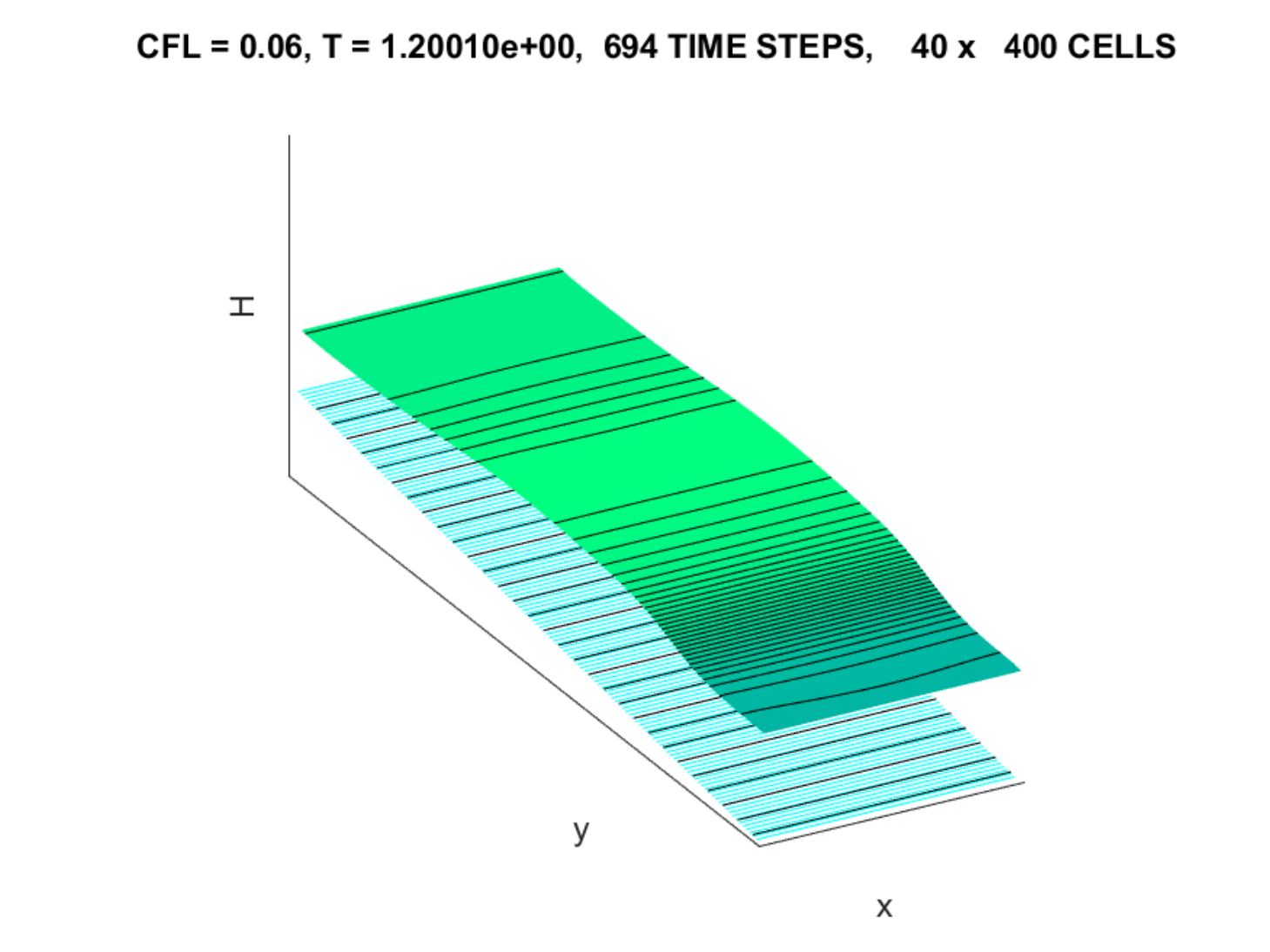}
	\includegraphics[width=0.38\textwidth]{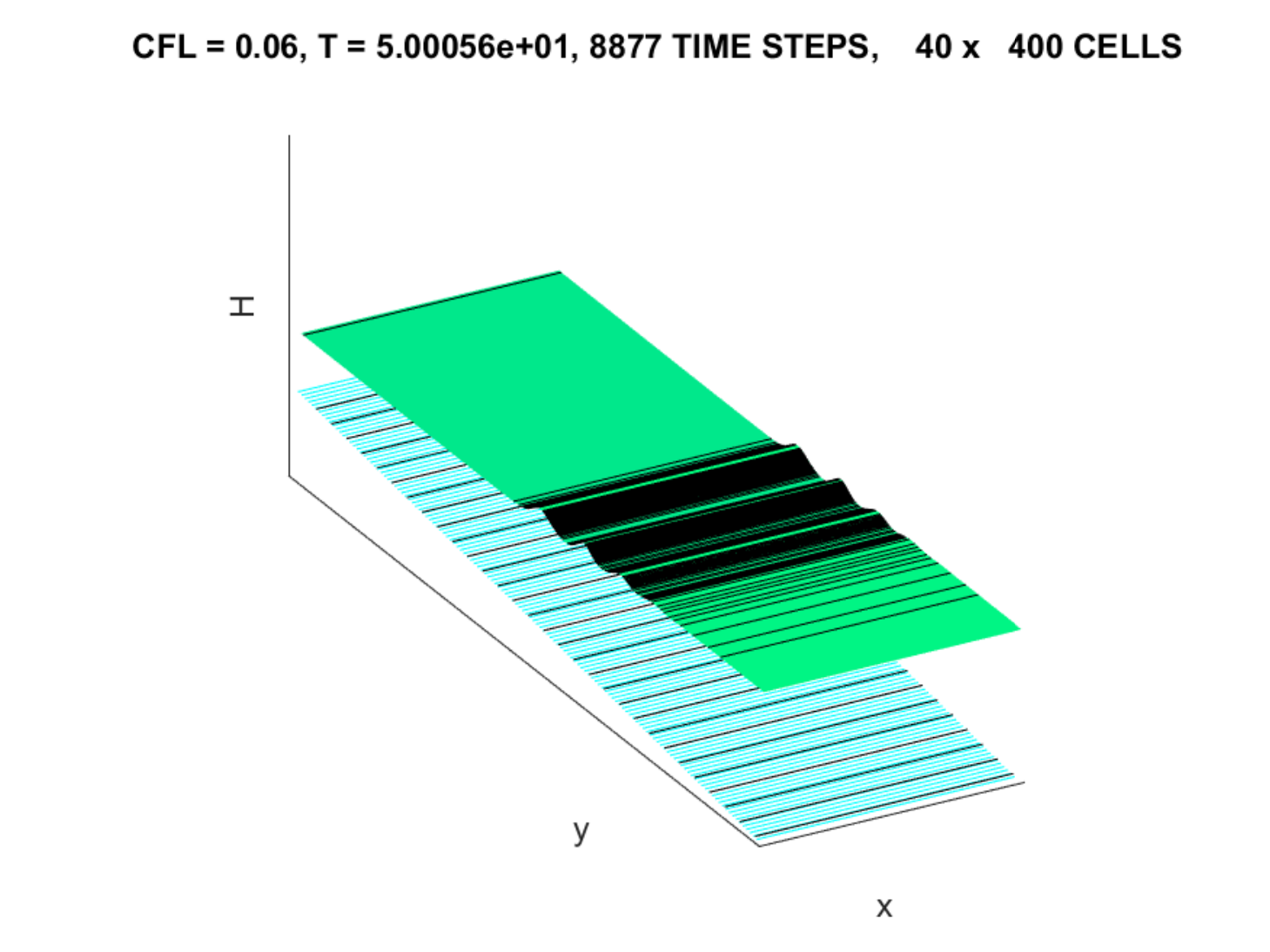}
	\includegraphics[width=0.38\textwidth]{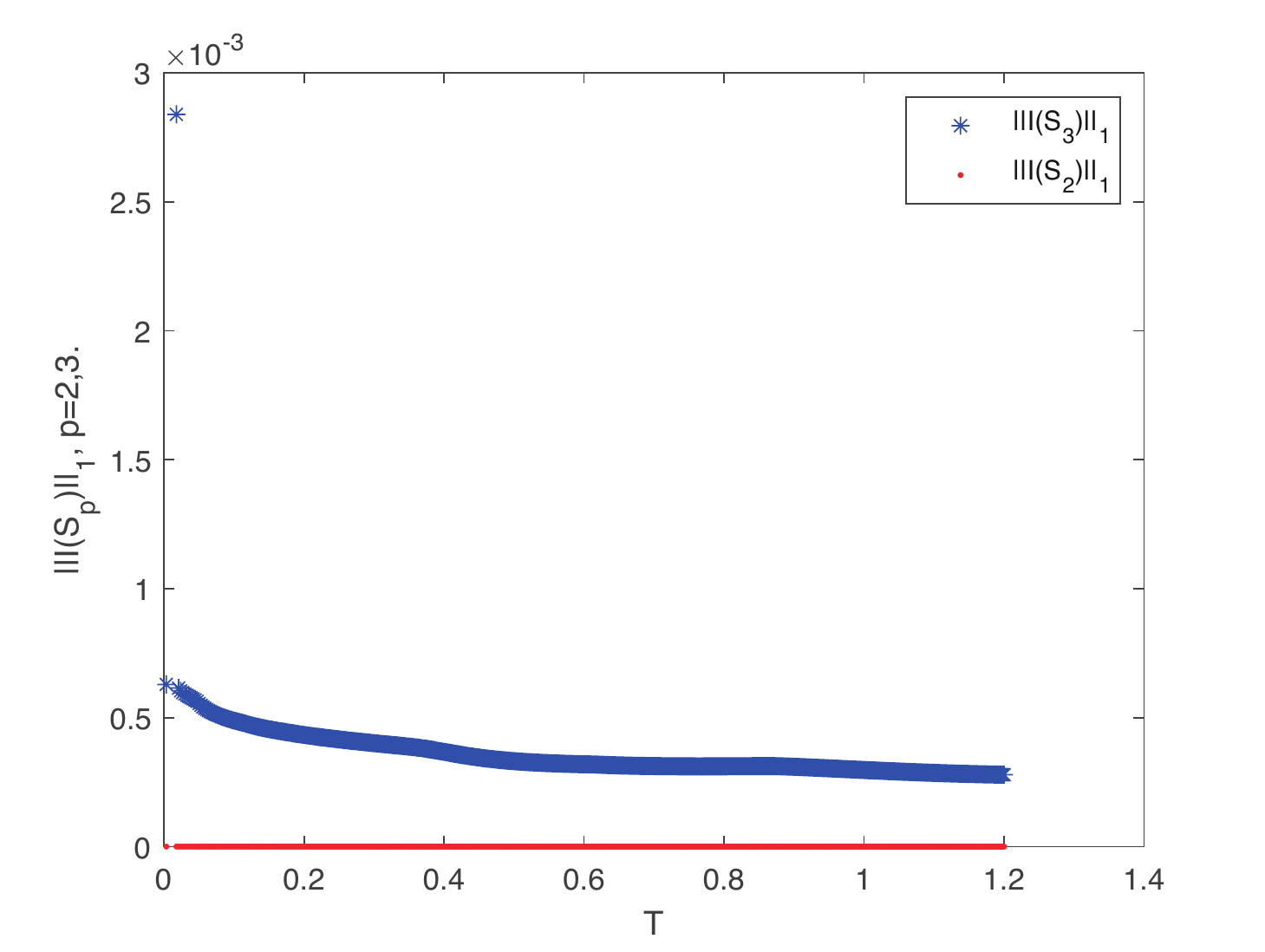}
	\includegraphics[width=0.38\textwidth]{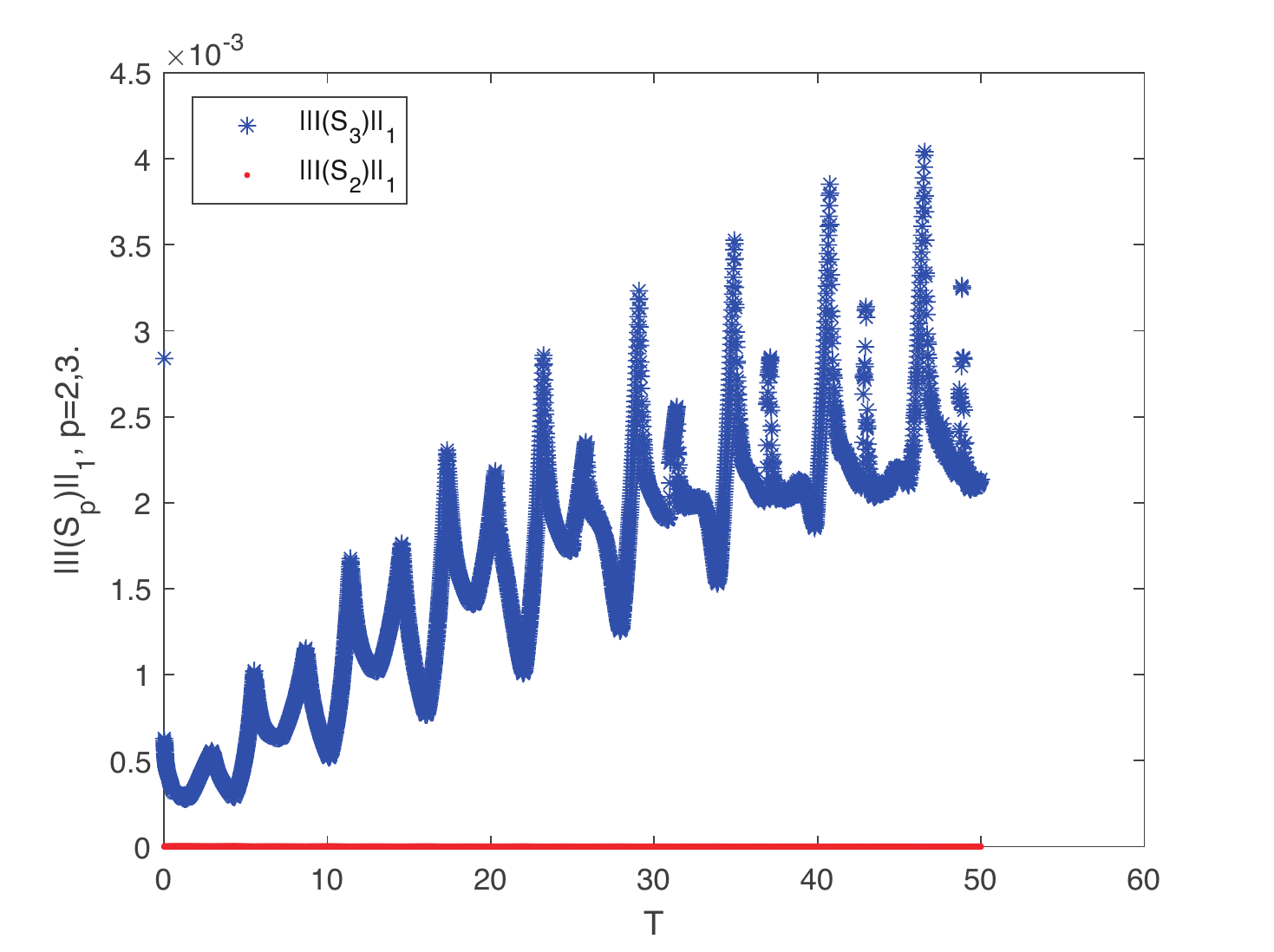}
	\caption{Sloping plane with discontinuous initial datum,
          40$\times$400 cells: water depth (top panels) and plots of
          the $L^1$-norm of $I(\Source^2)$ and $I(\Source^3)$ (bottom
          panels), at times $t=1.2$~s (left) and $t=50$~s
          (right). Computational time for $t=50$: 47.04~s.}
	\label{CI9}
\end{figure}

Figure~\ref{discID} shows the evolution of the normal water depth at
different times, until the final time $t=1.10$~s.
Again figure~\ref{CI9} reports the long-time simulations that show the
robustness over the treatment of the lake-at-rest conditions.

\section{Concluding remarks}
\label{sec:conclusion}

In this work, we present and discuss a new geometrically
intrinsic Lagrangian-Eulerian scheme able to handle non-autonomous
fluxes that arise from the spatially varying bottom topography.
The scheme is based on the concept of no-flow surfaces that define
the integral tube along which the initial value problem is solved.
The integral tube is defined by re-writing the ISWE as a time-space
divergence of the fluxes that takes into account the intrinsic
geometric features of the SW equations, written on a local reference
system anchored on the bottom surface.
A first Lagrangian evolution step that evolves the solution from time
$t^n$ to $t^{n+1}$ along the integral tube is followed by a second
Eulerian projection step where the current numerical approximation is
averaged over the original grid.
The scheme is efficient in terms of computational and memory cost
and, at the same time, is easy to implement since no (local) Riemann
problems are solved.  Moreover, it allows to handle seamlessly the
non-autonomous fluxes arising in the ISWE formulation.  We show
monotonicity and stability of the scheme, provided by a new weak
CFL-type constraint. These features are significant and ensure the
simplicity and robustness of this class of \textit{no-flow surface}
Lagrangian-Eulerian schemes.

This new scheme is also applicable to scalar equations and systems of
hyperbolic conservation laws as well as for modeling transport
problems in the presence of a general bottom topography with
non-negligible slopes and curvatures, such as a mountain landscape.
Numerical examples on mountain-like topographies are used to verify
the theoretical developments and illustrate the capabilities of the
proposed approach in a range of two-dimensional shallow water
equations. The results show the robusteness of the Lagrangian-Eulerian
formulation on Cartesian grids for relatively mild and slowly varying
curvatures verifying the ISWE hypothesis.

\subsection*{Acknowledgments}
 E. Abreu gratefully acknowledges the financial support of the the
 National Council for Scientific and Technological Development -
 Brazil (CNPq) (Grant 306385/2019-8). E. Bachini wishes to thank the
 German Research Foundation (DFG) for financial support within the
 Research Unit “Vector- and Tensor-Valued Surface PDEs” (FOR 3013)
 with project no. VO 899/22-1.

\appendix
\section{Geometrically intrinsic differential operators}
\label{app1}

We recall here the definition of the intrinsic differential
operators taking into account the metric tensor
(see~\cite{art:BP20}).
Let $(\Time,\xcl,\ycl)$ be the coordinate set of the LCS and
$\STFirst = \left\{\metrcoef{i}\right\}$ the
associated metric tensor.  Let
$\scalFun\From\Omega\To\REAL$ be a scalar function,
$\vectvel\From\Omega\To\REAL^3$ a contravariant vector field given
by
$\vectvel=\velcompContr[1]{}\normalvec +
\velcompContr[2]{}\vecBaseCCcv[1] +
\velcompContr[3]{}\vecBaseCCcv[2]$, and
$\tens\From\Omega\To\REAL^{3\times 3}$ a rank-2 contravariant
symmetric tensor
given by $\tens=\left\{\tenscomp{}\right\}$. Then, the differential
operators in the LCS are given by the following expressions:
\begin{itemize}
\item the gradient of $\scalFun$ is:
  \begin{equation*}
    \label{eq:gradf}
    \GradSurfST \scalFun = 
    \STFirst[]^{-1}
    \Grad \scalFun = 
    \left( 
    \DerPar{\scalFun}{\Time},\;
    \frac{1}{\metrcoefH{1}^2} \DerPar{\scalFun}{\xcl},\;
    \frac{1}{\metrcoefH{2}^2} \DerPar{\scalFun}{\ycl},\;
    \right) \, ;
  \end{equation*}    
\item the divergence of $\vectvel$ is:
  \begin{multline}
    \label{eq:divVect3D}
    \DivSurfST \vectvel 
    = \frac{1}{\sqrt{\DET{\First}}} \Div 
    \left( \sqrt{\DET{\First}} \vectvel \right)
    =\frac{1}{\metrcoefH{1}\metrcoefH{2}} 
    \left(
    \DerPar{\left(\metrcoefH{1}\metrcoefH{2}\velcomp[1]{}\right)}{\Time} +
    \DerPar{\left(\metrcoefH{1}\metrcoefH{2}\velcomp[2]{}\right)}{\xcl} +
    \DerPar{\left(\metrcoefH{1}\metrcoefH{2}\velcomp[3]{}\right)}{\ycl}
    \right) \, ;
  \end{multline}
\item the $j$-th component of the divergence $\tens$ is:
  \begin{multline}
    \label{eq:divTen3D}
    \left(\DivSurfST\tens\right)^{j} 
    = {\GradSurfST}_{i}\tenscomp[ij]{}
    = \frac{1}{\sqrt{\DET{\First}}} \partial_i 
    \left( \sqrt{\DET{\First}} \tenscomp[ij]{} \right)
    +\ChristSymb{ik}{j}\tenscomp[ik]{} =
    \\
    = \DivSurfST\tenscomp[(\cdot j)]{}
    + \frac{1}{\sqrt{\metrcoef{j}}}
    \left(
    2 \tenscomp[1j]{}\DerPar{\sqrt{\metrcoef{j}}}{\Time}
    - \tenscomp[11]{}\frac{\sqrt{\metrcoef{1}}}{\sqrt{\metrcoef{j}}}
    \DerPar{\sqrt{\metrcoef{1}}}{\svcomp[j]{}}
    \right) 
    + \frac{1}{\sqrt{\metrcoef{j}}}
    \left(
    2 \tenscomp[2j]{}\DerPar{\sqrt{\metrcoef{j}}}{\xcl}
    - \tenscomp[22]{}\frac{\sqrt{\metrcoef{2}}}{\sqrt{\metrcoef{j}}}
    \DerPar{\sqrt{\metrcoef{2}}}{\svcomp[j]{}}
    \right)\\
    + \frac{1}{\sqrt{\metrcoef{j}}}
    \left(
    2 \tenscomp[3j]{}\DerPar{\sqrt{\metrcoef{j}}}{\ycl}
    - \tenscomp[33]{}\frac{\sqrt{\metrcoef{3}}}{\sqrt{\metrcoef{j}}}
    \DerPar{\sqrt{\metrcoef{3}}}{\svcomp[j]{}}
    \right) \qquad \mbox{ (if $j=1$ the $\svcomp[j]=\Time$)}\\[1em]
    = \DivSurfST\tenscomp[(\cdot j)]{}
    + \frac{1}{\sqrt{\metrcoef{j}}}
    \left(
    2 \tenscomp[2j]{}\DerPar{\sqrt{\metrcoef{j}}}{\xcl}
    - \tenscomp[22]{}\frac{\sqrt{\metrcoef{2}}}{\sqrt{\metrcoef{j}}}
    \DerPar{\sqrt{\metrcoef{2}}}{\svcomp[j]{}}
    \right)
    + \frac{1}{\sqrt{\metrcoef{j}}}
    \left(
    2 \tenscomp[3j]{}\DerPar{\sqrt{\metrcoef{j}}}{\ycl}
    - \tenscomp[33]{}\frac{\sqrt{\metrcoef{3}}}{\sqrt{\metrcoef{j}}}
    \DerPar{\sqrt{\metrcoef{3}}}{\svcomp[j]{}}
    \right)\,,
  \end{multline}
  where $\DivSurfST\tenscomp[(\cdot j)]{}$ identifies the divergence of the
  $j$-th column of $\tens$, and $\ChristSymb{ij}{k}$ denote the
  Christoffel symbols.
\end{itemize}

By applying this definitions, we can prove the equivalence between
formulation \eqref{eq1swe} and \eqref{eq:eqSTdiv}.  Using
equation~\eqref{eq:divTen3D}, the divergence of the tensor $\Flux$ is
given by:
\begin{align*}
  \mbox{ROW 1: } &   \left(\DivSurfST\Flux\right)^{1}
                   = \DivSurfST\Flux^{(\cdot 1)} \\
  \mbox{ROW 2: } &   \left(\DivSurfST\Flux\right)^{2}
                   = \DivSurfST\Flux^{(\cdot 2)}
                   + \frac{1}{\metrcoefH{1}}
                    F^{21}\DerPar{\metrcoefH{1}}{\xcl}
                   + \frac{1}{\metrcoefH{1}}
                   \left(
                   2 F^{31}\DerPar{\metrcoefH{1}}{\ycl}
                   - F^{32}\frac{\metrcoefH{2}}{\metrcoefH{1}}
                   \DerPar{\metrcoefH{2}}{\xcl}
                   \right) \\
    \mbox{ROW 3: } &   \left(\DivSurfST\Flux\right)^{3}
                   = \DivSurfST\Flux^{(\cdot 3)}
                   + \frac{1}{\metrcoefH{2}}
                     \left(
                     2 F^{22}\DerPar{\metrcoefH{2}}{\xcl}
                     - F^{21}\frac{\metrcoefH{1}}{\metrcoefH{2}}
                     \DerPar{\metrcoefH{1}}{\ycl}
                     \right)
                     + \frac{1}{\metrcoefH{2}}
                      F^{32}\DerPar{\metrcoefH{2}}{\ycl}
\end{align*}
Applying the definition of the divergence of a vector given
in~\eqref{eq:divVect3D} to ROW 1 we obtain:
\begin{equation*}
  \DerPar{\Depth}{\Time}+
                   \frac{1}{\metrcoefH{1}\metrcoefH{2}} 
                   \left(
                   \DerPar{\left(\metrcoefH{1}\metrcoefH{2}\Qcomp[1]{}\right)}{\xcl} +
                   \DerPar{\left(\metrcoefH{1}\metrcoefH{2}\Qcomp[2]{}\right)}{\ycl}
                   \right) = \DerPar{\Depth}{\Time}+\DivSurf\Qvect\,,
\end{equation*}
while for ROWS 2 and 3 we need
to apply the expression for the divergence of a
tensor~\eqref{eq:eqSTdiv}, obtaining:
\begin{equation*}
  \DerPar{\Qvect}{\Time} + \DivSurf \FluxSWE\,.
\end{equation*}
These two equations provide the ISWE model~\eqref{eq1swe} in absence
of source term.

\bibliography{strings,biblio}
\bibliographystyle{elsarticle-num}

\end{document}